\documentclass[12pt]{article}
\usepackage{amsmath}
\usepackage{amssymb}
\usepackage{graphicx}
\usepackage[cp1251]{inputenc}
\usepackage[russian]{babel}
\newcommand{\il}[2]{\int\limits_{#1}^{#2}}
\newcommand{\ilp}[1]{\int\limits_{#1}^{+\infty}}

\newcommand{\ph}{\phantom{a}}
\newcommand{\phh}{\phantom{aaa}}
\newcommand{\ilpp}{\ilp{t_0}}
\newcommand{\sist}[2]{\left\{
\begin{array}{l}
{#1}\\
\ph\\
{#2}
\end{array}
\right.}
\newcommand{\mb}[1]{\ph\mbox{#1}\ph}
\newcommand{\liml}[1]{\lim\limits_{#1}}
\newcommand{\gr}[9]{
\begin{picture}(10,10)(#1,#2)
\qbezier(0,#3)(10,#4)(20,#5)
\qbezier(20,#5)(30,#6)(40,#7)
\qbezier(40,#7)(50,#8)(60,#9)
\end{picture}}
\newcommand{\graph}[9]{
\begin{picture}(10,10)
\qbezier(0,#1)(10,#2)(20,#3)
\qbezier(20,#3)(30,#4)(40,#5)
\qbezier(40,#5)(50,#6)(60,#7)
\qbezier(60,#7)(70,#8)(80,#9)
\end{picture}}
\newcommand{\grik}[5]{
\begin{picture}(10,10)(#1,#2)
\qbezier(0,#3)(10,#4)(20,#5)
\end{picture}}

\begin{document}


\vskip 20pt

\centerline{\bf The behavior of solutions of the systems of two first order}

 \centerline{\bf linear ordinary  differential equations}

\vskip 20 pt

\centerline{\bf G. A. Grigorian}

\vskip 10 pt

\centerline{0019 Armenia c. Yerevan, str. M. Bagramian 24/5}
\centerline{Institute of Mathematics NAS of Armenia}
\centerline{E - mail: mathohys2@instmath.sci.am}

\vskip 20 pt

\noindent
Abstract. The Riccati equation method is used for study the behavior of  solutions of the systems of two linear first order ordinary differential equations. All types of oscillation and regularity of these system are revealed. A generalization of Leighton's theorem is obtained. Three new principles for the second order linear differential equations are derived. Stability and non conjugation criteria are proved for the mentioned systems, as well as estimates are obtained for the solutions of the last ones.

\vskip 20 pt


\vskip 20 pt

\centerline{\bf \S 1. Introduction}

\vskip 20 pt

Let  $a_{jk}(t) \phantom{a} (j,k =1,2)$ be real valued continuous functions on $[t_0;+\infty)$. Consider the system of equations
$$
\left\{
\begin{array}{l}
\phi' = a_{11}(t) \phi + a_{12}(t) \psi;\\
\phantom{a}\\
\psi' = a_{21}(t) \phi + a_{22}(t) \psi, \phantom{a} t\ge t_0.
\end{array}
\right. \eqno (1.1)
$$
Study of  the asymptotic behavior problem (oscillation, non oscillation, non \linebreak conjugation, rate of growth) of solutions  linear systems of ordinary differential equations as well as  stability  problem of the last ones, in particular of the system (1.1), are  important problems  of the qualitative theory of differential equations and many works are devoted to them (see [1] and cited works therein, [2], [3],  [4],  [5], [6], [7]). Let  $p(t), \phantom{a} q(t)$  and $r(t)$  be real valued continuous functions on  $[t_0;+\infty)$, and let  $p(t) > 0, \phantom{a} t\ge t_0$.  Along with the system (1.1) consider the equation
$$
(p(t) \phi')' + q(t) \phi' + r(t) \phi = 0, \phantom{aaa} t\ge t_0. \eqno (1.2)
$$
The substitution
$$
p(t)\phi' = \psi   \eqno (1.3)
$$
in this equation reduces it to the system
$$
\left\{
\begin{array}{l}
\phi' = \phantom{aaaaaaaaaa} \frac{1}{p(t)} \psi;\\
\phantom{a}\\
\psi' = - r(t) \phi  - \frac{q(t)}{p(t)}(t) \psi, \phantom{a} t\ge t_0,
\end{array}
\right. \eqno (1.4)
$$
which is a particular case of the system (1.1).  For  $p(t) \equiv 1, \phantom{a} q(t)\equiv 0$  Eq. (1.2) takes the forme
$$
\phi'' + r(t) \phi = 0, \phantom{aaa} t\ge t_0. \eqno (1.5)
$$
It is well known (see for example [8]), that by using different transformations Eq. (1.2) can be reduced to the Eq. (1.5). One can show that the system (1.1) can be reduced to Eq. (1.5), if (for example) $a_{12}(t) \ne 0, \phantom{a} t\ge t_0$. There exist also other conditions for which the system (1.1) can be reduced to Eq.  (1.5). Of course the reduction of the system (1.1) to Eq. (1.5), if it is possible to carry it out (until now, it is not known whether this can always be done), can be very useful for study of different qualitative characteristics of the system (1.1). However this method not always can help to solve the assigned problem. One of effective methods of qualitative investigation of Eq. (1.5), as well as of the system  (1.1) is the Riccati equation method. In this work we use this method for the study of the behavior of  solutions of the system of two linear first order ordinary differential equations. We reveal  all types of oscillation and regularity of this system. We obtain a  generalization of Leighton's oscillation theorem. We derive three new principles for the second order linear ordinary differential equations. We prove some stability and non conjugation criteria for the mentioned above system. We also obtain  estimates  for the solutions of the last one.

\vskip 20pt

\centerline{\bf \S 2. Auxiliary propositions}

\vskip 20pt

Let  $a(t), \phantom{a} b(t)$    and   $c(t)$  be real valued continuous function on  $[t_0;+\infty)$. Consider the Riccati equation
$$
x' + a(t) x^2 + b(t) x + c(t) = 0, \phantom{aaa} t\ge t_0. \eqno (2.1)
$$
In this paragraph we study some important properties of global solutions (existing on $[t_1;+\infty)$  for some $t_1 \ge t_0$) of this equation which will be used further for the study of asymptotic properties of solutions of the system (1.1). Along with Eq. (2.1) consider the system of equations
$$
\left\{
\begin{array}{l}
\phi'= \phantom{aaaaaaaa} a(t) \psi;\\
\phantom{aa}\\
\psi' = - c(t) \phi - b(t) \psi, \phantom{aaa} t\ge t_0.
\end{array}
\right.
\eqno (2.2)
$$
The solutions  $x(t)$  of Eq. (2.1), existiong on some interval  $[t_1;t_2) \phantom{a} (t_0 \le t_1 < t_2 \le +\infty)$, are connected with the solutions $(\phi(t), \psi(t))$  of the system  (2.2) by the equalities  (see [9], pp. 153 - 154)
$$
\phi(t) = \phi(t_1) \exp\biggl\{\int\limits_{t_1}^t a(\tau) x(\tau) d\tau\biggr\}, \phantom{a} \phi(t_1) \ne 0, \phantom{a} \psi(t) = x(t) \phi(t).  \eqno (2.3)
$$
In this paragraph we will take that all solutions of equations and systems of equations are real valued. For brevity we introduce the denotations:
$$
J_u(t_1;t) \equiv\exp\biggl\{\int\limits_{t_1}^t a(\tau) u(\tau) d\tau\biggr\}, \phantom{a} J_u(t) \equiv J_u(t_0;t),
$$
$$
  I_{u,v}^+(t_1;t) \equiv \int\limits_{t_1}^t u(\tau) J_{-v}(t_1;\tau) d\tau, \phantom{a} I_{u,v}^+(t_1) \equiv \int\limits_{t_1}^{+\infty} u(\tau) J_{-v}(t_1;\tau) d\tau, \phantom{a}
$$
$$
I_{u,v}^- (t_1;t) \equiv \int\limits_{t_1}^t J_{-u}(\tau; t) v(\tau) d\tau, \phantom{a} t_1, t \ge t_0,
$$
where  $u(t)$  and $v(t)$  are arbitrary continuous functions on  $[t_0;+\infty)$.
Rewrite Eq. (2.1) in the form:
$$
x' + h_x(t) x + b(t) = 0, \phantom{aaa} t\ge t_0,
$$
where  $h_x(t) \equiv a(t) x + b(t), \phantom{a} t\ge t_0.$  Then by virtue of the Cauchy formula Eq. (2.1) is equivalent to the following integral equation
$$
x = J_{-h_x}(t_1;t)\biggl[x(t_1) - \int\limits_{t_1}^t J_b(t;\tau) c(\tau) \phi_{x}(t_1;\tau) d\tau\biggr], \phantom{aaa} t\ge t_0, \eqno (2.4)
$$
where  $\phi_x(t_1;t) \equiv \exp\biggl\{\int\limits_{t_1}^t a(\tau) x(\tau) d\tau\biggr\}, \phantom{a} t_1,t \ge t_0.$  Let  $a_1(t), \phantom{a} b_1(t)$ and $c_1(t)$ be real valued continuous function on  $[t_0;+\infty)$. Along with Eq. (2.1) consider the equation
$$
x' + a_1(t) x^2 +  b_1(t) x + c_1(t) = 0, \phantom{aaa} t\ge t_0, \eqno (2.5)
$$
and differential inequality
$$
\eta' + a(t) \eta^2 +  b(t) \eta  + c(t) \ge 0, \phantom{aaa} t\ge t_0. \eqno (2.6)
$$
Note that for $a(t) \ge 0, \phantom{a} t\ge t_0$, any solution of the linear equation
$$
\eta'+  b(t) \eta  + c(t) = 0, \phantom{aaa} t\ge t_0,
$$
is a solution of (2.6). Therefore for any initial condition  $\eta_{(0)}$ inequality (2.6) has a solution   $\eta_0(t)$   on  $[t_0;+\infty)$  with $\eta_0(t_0) = \eta_{(0)}$.

{\bf Theorem 2.1}. {\it Let Eq. (2.5) has a solution  $x_1(t)$  on  $[t_0;\tau_0)\phantom{a} (\tau_0 \le +\infty)$  and let the following condition be satisfied:
$$
a(t) \ge 0, \phantom{a} \int\limits_{t_0}^t \exp\bigg\{\int\limits_{t_0}\bigl[a(\xi)\bigl(\eta_0(\xi) + x_1(\xi)\bigr) + b(\xi)\bigr]d\xi\biggr\}\times
$$
$$
\times\bigl[\bigl(a_1(\tau) - a(\tau)) x_1^2(\tau) + \bigl(b_1(\tau) - b(\tau)) x_1(\tau) + c_1(\tau) - c(\tau)\bigr] d\tau \ge 0, \phantom{aaa} t\in [t_0;\tau_0),
$$
where $\eta_0(t)$ is a solution of (2.6) on $[t_0;\tau_0)$  with  $\eta_0(t_0) \ge x_1(t_0)$.
Then for each $x_{(0)} \ge x_1(t_0)$  Eq. (2.1) has a solution $x_0(t)$  on  $[t_0;\tau_0)$,  and
$x_0(t) \ge x_1(t), \phantom{a} t\in [t_0;\tau_0)$.}

See proof in [10].

Let  $t_1 \ge t_0$.

{\bf Definition 2.1.} {\it  A solution of Eq. (2.1) is called $t_1$-regular, if it exists on $[t_1;+\infty)$. Eq. (2.1) is called regular if it has a
 $t_1$-regular solution for some  $t_1 \ge t_0$.}

{\bf Definition 2.2.} {\it A $t_1$-regular solution  $x(t)$ of Eq. (2.1) is called   $t_1$-normal, if there exists a neighborhood $U_x(t_1)$ of the point $x(t_1)$ such that each solution  $\widetilde{x}(t)$ of Eq.  (2.1) with  $\widetilde{x}(t_1) \in U_x(t_1)$  is  $t_1$-regular. Otherwise $x(t)$  is called    $t_1$-extremal.}

{\bf Remark 2.1.} {\it From the results of work  [11] it follows that for some $t_1 \ge t_0$ the regular equation (2.1) can have: the unique  $t_1$-regular solution (then it is  $t_1$-extremal); no    $t_1$-extremal solution (then its all $t_1$-regular solutions are $t_1$-normal);  the unique $t_1$-extremal solution (and all  other $t_1$-regular solutions are $t_1$-normal);  two  $t_1$- extremal solutions (all other  $t_1$-regular solutions are $t_1$-normal).}

In what follow we will assume that the functions  $a(t)$  and  $c(t)$ have unbounded supports (the case when one of these functions has a bounded support is trivial). For arbitrary continuous   function $u(t)$  on  $[t_0;+\infty)$ denote:
$$
\mu_u(t_1;t) \equiv \int\limits_{t_1}^ta(\tau)\exp\biggl\{ - \int\limits_{t_1}^\tau \bigl[2a(\xi) u(\xi) + b(\xi)] d\xi\biggr\} d\tau,\phantom{aaaaaaaaaaaaaaaaaaaaaaa}
$$
$$
\phantom{aaaaaaaaa}\nu_u(t) \equiv \int\limits_t^{+\infty}a(\tau)\exp\biggl\{ - \int\limits_t^\tau \bigl[2a(\xi) u(\xi) + b(\xi)] d\xi\biggr\} d\tau, \phantom{aaa} t_1, t\ge t_0.
$$

{\bf Theorem 2.2.} {\it Let for some  $t_1$-regular solution $x_0(t)$  of Eq. (2.1) the integral   $\nu_{x_0}(t_1)$ be convergent. Then the following assertions are valid.

\noindent
A) For any $t\ge t_1$  and for all  $t_1$-normal solutions $x(t)$  of Eq. (2.1) and only for them the integrals  $\nu_x(t)$  converge.

\noindent
B) In order that Eq. (2.1) have  $t_1$-extremal solution it is necessary and sufficient that
$$
\nu_{x_0}(t) \ne 0, \phantom{aaa} t\ge t_1.
$$
Under this condition  the unique  $t_1$-extremal solution  $x_*(t)$ is defined by formula
$$
x_*(t) = x_0(t) - \frac{1}{\nu_{x_0}(t)}, \phantom{aaa} t\ge t_1, \eqno (2.7)
$$
and
$$
\nu_{x_*}(t) = +\infty,  \phantom{a} t\ge t_1,\phantom{a}  \mbox{ or } \phantom{a}\nu_{x_*}(t) = -\infty,  \phantom{a} t\ge t_1, \eqno (2.8)
$$
$$
\int\limits_t^{+\infty} a(\tau)\bigl[x_1(\tau) - x_2(\tau)\bigr] d\tau = \ln \biggl[\frac{x_*(t) - x_1(t)}{x_*(t) - x_2(t)}\biggr], \phantom{aaa} t\ge t_1, \eqno (2.9)
$$
$$
\int\limits_{t_1}^{+\infty} a(\tau)\bigl[x_1(\tau) - x_2(\tau)\bigr] d\tau =  - \infty. \eqno (2.10)
$$
}

Proof. All assertions of this theorem except (2.8) and (2.10), are proved in [11]. Let as prove  (2.8). We will use the equalities (see  [11]):
$$
\mu_{x_*}(t_2;t) = \frac{\mu_{x_0}(t_2;t)}{1 + \lambda_*(t_2)\mu_{x_0}(t_2;t)},  \eqno (2.11)
$$
$$
x_1(t) = x_2(t) + \frac{\lambda_{12}(t_2)\exp\biggl\{- \int\limits_{t_1}^t [2 a(\tau) x_2(\tau) + b(\tau)] d\tau \biggr\}}{1 + \lambda_{12}(t_2)\mu_{x_2}(t_2;t)}, \phantom{a} t\ge t_2 \ge t_1, \eqno (2.12)
$$
where  $\lambda_*(t_2)\equiv x_*(t_2) - x_0(t_2), \phantom{a} \lambda_{12}(t_2)\equiv x_1(t_2) - x_2(t_2), \phantom{a} t\ge t_2 \ge t_1, \phantom{a} x_1(t)$  and  $x_2(t)$  are arbitrary  $t_1$-regular solutions of Eq. (2.1). From (2.12) it follows that $\mu_{x_0}(t_2;t)$  is bounded by  $t$  on  $[t_2;+\infty)$. Then since obviously
$$
\nu_{x_0}(t_2) = \lim \limits_{t\to +\infty} \mu_{x_0}(t_2;t) \ne 0,\phantom{aaa} t_2 \ge t_1, \eqno (2.13)
$$
necessarily
$$
\lim\limits_{t\to +\infty} [1 + \lambda_*(t_2)\mu_{x_0}(t_2;t)] = 0, \phantom{aaa} t_2 \ge t_1. \eqno (2.14)
$$
From here from (2.11) and  (2.13) it follows (2.8). Let us prove  (2.10). By (2.12) we have:
$$
x_*(t) - x_0(t) = \frac{\lambda_{*}(t_1)\exp\biggl\{- \int\limits_{t_1}^t [2 a(\tau) x_0(\tau) + b(\tau)] d\tau \biggr\}}{1 + \lambda_{*}(t_1)\mu_{x_0}(t_1;t)}, \phantom{a} t \ge t_1.
$$
Multiplying both sides of this equality on  $a(t)$  and integrating from $t_1$  to  $t$ we will get
$$
\int\limits_{t_1}^t a(\tau)[x_*(\tau) - x_0(\tau)] d\tau = \ln \bigl[1 + \lambda_{*}(t_1)\mu_{x_0}(t_1;t)\bigr], \phantom{aaa} t\ge t_1.
$$
Passing to the limit in this equality when  $t\to +\infty$  and taking into account (2.14) we obtain  (2.10). The theorem is proved.

{\bf Corollary 2.1.} {\it If for some $t_1$-regular solution  $x_*(t)$
the equality $\nu_{x_*}(t_1) = \pm\infty$ is fulfilled, then  $x_*(t)$ is the unique $t_1$-extremal solution of Eq. (2.1), and Eq. (2.1) has $t_1$-normal solutions,
and for each $t\ge t_1$  and for all  $t_1$-normal solutions $x(t)$  of Eq. (2.1) the integrals  $\nu_x(t)$    converge; for every normal solutions  $x_0(t), \phantom{a} x_1(t), \phantom{a} x_2(t)$
of Eq. (2.1) and for  $x_*(t)$  the correlations (2.7) - (2.10) are satisfied.}

Proof. Let  $\nu_{x_*}(t_1) = + \infty$  (the proof in the case  $\nu_{x_*}(t_1) = - \infty$ by analogy).  Then
$$
\mu_{x_*}(t_1;t) > 1, \phantom{aaa} t\ge T, \eqno (2.15)
$$
for some  $T > t_1$.  Let  $\overline{\mu}_{x_*}\equiv \max\limits_{t\in [t_0;T]}|\mu_{x_*}(t_1;t)|$, and let $x_0(t)$  be a solution to Eq. (2.1) with  $x_*(t_1) < x_0(t_1) < x_*(t_1) + \frac{1}{\max\{1, \overline{\mu}_{x_*}\}}$. Then taking into account  (2.15) we will have:
\pagebreak
$$
1 + \lambda_*(t_1)\mu_{x_*}(t_1;t) > 0, \phantom{aaa} t\ge t_1,
$$
where  $\lambda_*(t_1)\equiv x_)(t_1) - x_*(t_1) > 0$. Hence  (see  [11]) by  (2.11)  $x_0(t)$  is a  $t_1$-regular solution of Eq. (2.1). Show that the integrals   $\nu_{x_*}(t)$
converge for all  $t\ge t_1$  and   $\nu_{x_*}(t) \ne 0, \phantom{a}$.
We use the equality   (see  [11])
$$
\mu_{x_0}(t_2;t) = \frac{\mu_{x_*}(t_2;t)}{1 + \lambda_*(t_2)\mu_{x_*}(t_2;t)}, \phantom{aaa} t\ge t_2\ge t_1,
$$
where  $\lambda_*(t_2) \equiv x_0(t_2) - x_*(t_2) \ne 0.$
For enough large values of  $t > t_2$ we have  $\mu_{x_*}(t_2;t) > 0$. Therefore,
$$
\nu_{x_0}(t_2) = \lim\limits_{t\to +\infty} \mu_{x_0}(t_2;t) = \lim\limits_{t\to +\infty} \frac{1}{\lambda_*(t_2) + \frac{1}{\mu_{x_*}(t_2;t)}} = \frac{1}{\lambda_*(t_2)} \ne 0.
$$
So for  $x_0(t)$ all conditions of Theorem 2.2 are fulfilled. Therefore Eq. (2.1)
has  $t_1$-normal solutions and for every  $t_1$-normal solutions  $x(t)$ of Eq. (2.1)
and for all  $t\ge t_1$ the integrals  $\nu_x(t)$ converge; also for every  $t_1$-normal solutions $x_0(t), \phantom{a} x_1(t), \phantom{a}x_2(t)$ of Eq. (2.1) and for  $x_*(t)$
the correlations  (2.7)  - (2.10) are satisfied. The corollary is proved.

Denote by  $reg(t_1)$ the set of values  $x_{(0)}\in R$,
for which the solution $x(t)$ of Eq. (2.1) with   $x(t_1)=x_{(0)}$ is $t_1$-regular.

{\bf Lemma 2.1}. {\it  Let  $a(t) \ge 0, \phantom{a} t\ge t_0$,
and let Eq. (2.1) has $t_1$-regular solution. Then it has the unique $t_1$-extremal solution  $x_*(t)$,  and $reg(t_1) = [x_*(t_1);+\infty)$.}

See the proof in  [2].

{\bf Lemma 2.2.} {\it let  $a(t) \ge 0, \phantom{a} t\ge t_0; \phantom{a} t_0 \le t_1 < t_2$,  and let  $(t_1;t_2)$ be the maximal existence interval for the solution  $x(t)$ of Eq. (2.1).
Then $\lim\limits_{t\to t_1 + 0} x(t) = +\infty$.}

See the proof in  [10].

{\bf Lemma 2.3}. {\it Let $a(t) \ge 0, \phantom{a} t\ge t_0, \phantom{a} x_0(t)$ be a $t_0$-normal solution of Eq. (2.1), $x_0(t) \ne~0, \phantom{a} t\ge t_0.$   Then for its unique $t_0$-extremal solution $x_*(t)$   the equality
$$
\int\limits_{t_0}^t a(\tau) x_*(\tau) d \tau = - \ln \nu_{x_0}(t_0) +
\ln\biggl[\exp\biggl\{\int\limits_{t_0}^t a(\xi) x_0(\xi) d \xi\biggr\}\times\phantom{aaaaaaaaaaaaa}
$$
$$
\times\int\limits_t^{+\infty}\frac{a(s)x_0(s)}{x_0(t_0)}
\exp\biggl\{\int\limits_{t_0}^s\biggl[\frac{c(\xi)}{x_0(\xi)} - a(\xi)x_0(\xi)\biggr]d\xi\biggr\} d s\biggr], \phantom{aaa} t\ge t_0. \eqno (2.16)
$$
holds}

Proof. By Lemma 2.1 Eq. (2.1) has a  $t_0$-normal solution $x_0(t)$.
Then since  $a(t) \ge 0, \phantom{a} t\ge t_0$ and has unbounded support, the integral  $\nu_{x_0}(t)$
converges for all $t\ge t_0$   and  $\nu_{x_0}(t)\ne 0, \phantom{a} t\ge t_0.$
By virtue of Theorem 2.2  from here it follows that Eq. (2.1) has the unique $t_0$-extremal solution  $x_*(t)$, satisfying the equality  $x_*(t) = x_0(t) - \frac{1}{\nu_{x_0}(t)}, \phantom{a} t\ge t_0$.    From here it follows:
$$
\int\limits_{t_0}^t a(\tau) x_*(\tau) d \tau = \int\limits_{t_0}^t a(\tau) x_0(\tau) d \tau - \int\limits_{t_0}^t \frac{a(\tau)}{\nu_{x_0}(\tau)} d\tau = \ln\bigr [\exp\bigl\{\int\limits_{t_0}^t a(\tau) x_0(\tau) d\tau\bigr\}\bigr] -
$$
$$
-\ln \nu_{x_0}(t_0) + \int\limits_{t_0}^t d\biggl(\ln\biggl[\int\limits_\tau^{+\infty} a(s)\exp\biggl\{- \int\limits_{t_0}^s\bigl[2 a(\xi) x_0(\xi) + b(\xi)\bigr]d\xi\biggr\} d s\biggr]\biggr), \eqno (2.17)
$$
$t\ge t_0$.
On the strength of (2.1) from the condition  $x_0(t) \ne 0, \phantom{a} t\ge t_0$, it follows:
$$
2a(\xi) x_0(\xi) + b(\xi) = - \frac{x_0'(\xi) - a(\xi) x_0^2(\xi) + c(\xi)}{x_0(\xi)}, \phantom{a} \xi \ge t_0.
$$
From here and from (2.17) it follows  (2.16). The lemma is proved.

{\bf Lemma 2.4.} {\it Let  $a(t) \ge 0, \ph c(t) \ge 0, \ph t\ge t_0, \ph I_{a,b}^+(t_0) = +\infty$                 and let Eq. (2.1) has a solution on $[t_1;+\infty)$ for some  $t_1\ge t_0$. Then Eq. (2.1) has a positive solution on $[t_1;+\infty)$.}

See proof in  [5].

{\bf Theorem 2.3}. {\it Let  $a(t) \ge 0, \phantom{a} c(t) \le 0, \phantom{a}  t\ge t_0$.                    Then the following assertions are valid.

\noindent
I$^\circ$). For each $x_{(0)} \ge  \frac{-1}{I_{a,b}^+(t_0)}$ (for $I_{a,b}^+(t_0)=+\infty$   we take that $ \frac{1}{I_{a,b}^+(t_0)}=0$)
Er. (2.1) has a $t_0$-regular solution $x_0(t)$ with  $x_0(t_0) = x_{(0)}$,  and
$$
\frac{x_{(0)} J_{-b}(t)}{1 + x_{(0)} I_{a,b}^+(t_0;t)} \le x_0(t) \le x_{(0)} J_{-b}(t) - I_{a,b}^-(t_0;t), \phantom{aaa} t\ge t_0, \eqno (2.18)
$$
moreover if $x_{(0)} = 0$, then there exists $t_1 \ge t_0$ such that  $x_0(t) = 0, \phantom{a} t\in [t_0;t_1], \linebreak x_0(t) > 0, \phantom{a} t > t_1$.  If  $x_{(0)} > 0$ then  $x_0(t) > 0, \phantom{a} t\ge t_0$.

\noindent
II$^\circ$). The unique  $t_0$-extremal solution $x_*(t)$ of Eq. (2.1) is negative.

\noindent
III$^\circ$). If $I_{a,b}^+(t_0) = +\infty$   or $I_{c,-b}^+(t_0) = -\infty$,  then for each solution  $x(t)$ of Eq. (2.1) with  $x(t_0) \in (x_*(t_0); 0)$ there exists $t_2=t_2(x) \ge t_1 = t_1(x) > t_0$  such that   $x(t) < 0, \phantom{a} t\in [t_0;t_1), \phantom{a} x(t) = 0, \phantom{a} t\in[t_1;t_2]$  and  $x(t)> 0, \phantom{a} t > t_2.$

\noindent
IV$^\circ$). If $I_{a,b}^+(t_0) = +\infty$,  then foe each  $t_0$-normal
solution  $x_N(t)$  of Eq. (2.1) the equality  $\int\limits_{t_0}^{+\infty}a(\tau) x_N(\tau) d\tau =+\infty$. is fulfilled.

\noindent
V$^\circ$). If $I_{c,-b}^+(t_0) = -\infty$, then   $\int\limits_{t_0}^{+\infty}a(\tau) x_*(\tau) d\tau =-\infty$,  where  $x_*(t)$ is the unique $t_0$-extremal solution of Eq. (2.1).

\noindent
VI$^\circ$). If  $I_{a,b}^+(t_0) < + \infty$   and   $I_{-c,-b}^+(t_0) < +\infty$,  then Eq. (2.1) has a negative   $t_0$-normal solution  $x_N^-(t)$ such that for each solution  $x(t)$ of Eq. (2.1) with  $x_0(t_0) \in (x_N^-(t_0);0)$ there exists $t_2=t_2(x) \ge t_1 = t_1(x) > t_0$ such that  $x(t) < 0, \phantom{a} t\in [t_0;t_1), \phantom{a} x(t) = 0, \phantom{a} t\in[t_1;t_2], \phantom{a} x(t)> 0, \phantom{a} t > t_2.$

\noindent
VII$^\circ$). If
$$
I_{a,b}^+(t_0) = +\infty, \phantom{a} \int\limits_{t_0}^{+\infty}|c(\tau)|I_{-b,a}^-(t_0;\tau) d\tau < +\infty, \eqno (2.19)
$$
then  $\int\limits_{t_0}^{+\infty} a(\tau) x_*(\tau) d\tau > - \infty, \phantom{a} \int\limits_{t_0}^{+\infty} c(\tau) u_*(\tau) d\tau =+\infty$, where $u_*(t)$ is the unique   $t_0$-extremal solution of the equation
$$
u'  - c(t) u^2 - b(t) u - a(t) = 0, \phantom{aaa} t\ge t_0.
$$
}

Proof. Set $a_1(t) = a(t) \phantom{a} b_1(t) = b(t), \phantom{a} t\ge t_0, \phantom{a} c_1(t) \equiv 0$.  Then for each $x_{(0)} \ge \frac{-1}{I_{a,b}^+(t_0)}\stackrel{def}{=}~A$ the function  $x_1(t)\equiv \frac{x_{(0)} J_{-b}(t)}{1 + x_{(0)} I_{a,b}^+(t_0;t)}$  is a  $t_0$-regular solution of Eq. (2.5), and the conditions of Theorem 2.1 are fulfilled.
Therefore for each $x_{(0)} \ge A$ Eq. (2.1) has a $t_0$-regular solution $x_0(t)$ with  $x_0(t_0) = x_{(0)}$ and the first of conditions of inequalities (2.18) is satisfied. Set  $a_1(t) = a(t), \phantom{a} b_1(t) = b(t),\ph c_1(t) = c(t), \ph t\ge t_0.$  Then by already proven Eq. (2.5) will have  $t_0$-regular solutions, coinciding  wit the  $t_0$-regular solutions of Eq. (2.1). In the Eq. (2.1) set: $a(t)\equiv 0$. Then $x_2(t)\equiv x_{(0)} J_{-b}(t_0;t) - I_{b,c}^{-}(t_0;t)$ is a $t_0$-regular solution of Eq. (2.1). Obviously in this case the conditions of Theorem 2.1. are satisfied. Therefore the second of the inequalities (2.18) is fulfilled.
Let  $x_{(0)}=0$. Then since $c(t)$ has unbounded support by virtue of (2.17) from the inequality  $c(t)\le 0, \ph t\ge t_0$,  it follows  existence of $t_1 > t_0$  such that $x(t)=0,\ph t\in [t_0;t_1]$  and  $x(t) > 0, \ph t>t_1$.
The assertion  I$^\circ$  is proved. Prove II$^\circ$.  Let $x_0(t)$  be a solution of Eq. (2.1) with  $x_0(t_0) > 0$. By virtue of I$^\circ$ $x_0(t)$ is $t_0$-normal and positive. Therefore From Theorem 2.2 it follows that for each  $t\ge t_0$ the integral  $\nu_{x_0}(t)$ converges. Obviously  $\nu_{x_0}(t) > 0, \ph t\ge t_0$. Then by virtue of the same Theorem 2.2 $x_*(t)\equiv x_0(t) - \frac{1}{\nu_{x_0}(t)}$ is the unique  $t_0$-extremal solution of Eq. (2.1). Show that $x_*(t)\le 0, \ph t\ge t_0.$  By virtue of the first of inequalities (2.18) we have:
$$
\frac{x_0(t) J_{-b}(t;s)}{! + x_0(t) I _{a,b}(t;s)} \le x_0(s), \phh t_0\le t \le s.
$$
Multiplying both sides of this inequality on  $a(s)$  and integrating by  $s$  from  $t$  to  $\tau$
we  get: $\ln [1 + x_0(t) I_{a,b}(t;\tau)] \le \il{t_0}{t} a(s) x_0(s) ds, \ph t_0\le t \le s.$  Then
$$
\nu_{x_0}(t) \le \ilp{t}\frac{a(\tau) J_{-b}(t;\tau) d\tau}{1 + x_0(t) I_{a,b}^+(t;\tau)} = - \frac{1}{x_0(t)}\ilp{t}d\left(\frac{1}{1 + x_0(t) I_{a,b}^+(t;\tau)}\right) =\phantom{aaaaaaaaaaaaaaaaaaaaa}
$$
$$
\phantom{aaaaaaaaaaaaaaaaaaaaaaaaaaaa}=\frac{1}{x_0(t)}\left[1 - \frac{1}{1 + x_0(t) I_{a,b}^+(t;+\infty)}\right] \le \frac{1}{x_0(t)}, \phh t\ge t_0.
$$
From here it follows that  $x_*(t) \le 0, \ph t\ge t_0$. Show that
$$
x_*(t) < 0, \phh t\ge t_0. \eqno (2.20)
$$
Suppose that it is not true. Then since $x_*(t) \le 0, \ph t\ge t_0$, there exists  $t_1 \ge t_0$  such that $x_*(t_1) = 0$. By the first of the inequalities (2.18) from here it follows that  $x_*(t) \ge 0, \ph t\ge t_1$. Hence  $x_*(t) \equiv 0$ on  $[t_1; +\infty)$, which is impossible (since on  $[t_1;+\infty) \ph c(t) \not\equiv 0$.) The obtained contradiction proves (2.20),
and therefore the assertion II$^\circ$ does. Prove III$^\circ$.  Let $x_0(t)$ and  $x_1(t)$  be solutions of Eq.  (2.1) with the initial conditions  $x_0(t_0) > 0, \ph x_1(t_0) \in (x_*(t_0); 0)$. By virtue of Lemma  2.1 $x_0(t)$  and $x_1(t)$ are $t_0$-normal. Therefore by  (2.8) we have
$$
\ilp{t_0} a(\tau) [x_0(\tau) - x_1(\tau)] d\tau < + \infty. \eqno (2.21)
$$
Suppose $I_{a,b}^+(t_0) = +\infty$. Then from the first of inequalities (2.18) it follows:
$$
\ilp{t_0}a(\tau) x_0(\tau) d\tau \ge \ln [1 + x_0(t_1) I_{a,b}^+(t_0)] = +\infty, \eqno (2.22)
$$
Show that there exists  $\widetilde{t}_1 \ge t_0$  such that $x_1(\widetilde{t}_1) = 0$. Suppose that it is not true. Then  $x_1(t) < 0, \ph t\ge t_0$.  Taking into account (2.18) from here we obtain:
\pagebreak
$$
\ilp{t_0} a(\tau)(x_0(\tau) - x_1(\tau)) d\tau \ge \ilp{t_0} a(\tau) x_0(\tau) d\tau = +\infty, \eqno (2.23)
$$
which contradicts (2.21). The obtained contradiction shows that $x_1(\widetilde{t}_1) = 0$ for some  $\widetilde{t}_1 >~t_0$. Since $c(t)$ has unbounded support by virtue of (2.4) from here and from non positivity of  $c(t)$   it follows that  $x_1(t) < 0, \ph t\in [t_0;t_1), \ph x_1(t) = 0, \ph t\in [t_1;t_2]$  and  $x_1(t) >~0, \ph t> t_2$, for some $t_2 \ge t_1 > t_0$. Suppose  $I_{c,-b}^+(t_0) = - \infty$. Consider the equation
$$
u' - c(t) u^2 - b(t) u - a(t) = 0, \phh t\ge t_0. \eqno (2.24)
$$
By II$^\circ$  the unique $t_0$-extremal solution $u_*(t)$ of this equation is negative.
Therefore $\widetilde{x}_*(t) \equiv \frac{1}{u_*(t)}$ is a  $t_0$-regular solution of Eq. (2.1). By already proven, from here and from the equality  $I_{c,-b}^+(t_0) = -\infty$ it follows that each solution  $u(t)$
of Eq. (2.24) with  $u(t_0) \in (u_*(t_0);0)$ vanishes on $[t_0;+\infty)$.  Therefore each solution $x(t)$  of Eq. (2.1) with  $x(t_0) < \widetilde{x}(t_0)$ is not  $t_0$-regular.
By virtue of Lemma 2.1 from here it follows that
$$
u_*(t) = \frac{1}{x_*(t)}, \phh t\ge t_0. \eqno (2.25)
$$
Suppose that some solution  $\widetilde{x}(t)$ of Eq. (2.1) with  $\widetilde{x}(t_0)\in (x_*(t_0); 0)$ is negative.  Then $\widetilde{u}(t)\equiv \frac{1}{\widetilde{x}(t)}, \ph t\ge t_0,$  is a $t_0$-regular solution of Eq. (2.24), and  $\widetilde{u}(t_0) = \frac{1}{\widetilde{x}(t_0)} < \frac{1}{x_*(t_0)}$. From here and from (2.25) it follows that  $\widetilde{u}(t_0) < u_*(t_0)$, which contradicts Lemma~ 2.1.  The obtained contradiction shows that for any solution $x(t)$ of Eq. (2.1) with  $x(t_0) \in (x_*(t_0); 0)$  there exists $t_1= t_1(x) > t_0$ such that  $x(t_1)=0, \ph x(t) < 0, \ph t\in [t_0;t_1).$ Since $c(t)$ has unbounded support by (2.4) from here and from non negativity of $c(t)$  it follows that there exists  $t_2= t_2(x) \ge t_1$  such that  $x(t) = 0, \ph t\in [t_1;t_2], \ph x(t) >~0, \linebreak t > t_2$.  The assertion III$^\circ$ is proved. Prove IV$^\circ$.  Let  $x_+(t)$ be a solution of Eq. (2.1) with  $x_+(t_0) =1$. On the strength of Lemma 2.1 from the assertion I$^\circ$  it follows that $x_+(t)$ is $t_0$-normal. By the first of the inequalities (2.18) we have $x_+(t) \ge J_{-b}(t), \ph t\ge t_0$.  Let  $I_{a,b}^+(t_0) = +\infty$.  Then from the last inequality it follows that
$$
\ilp{t_0}a(\tau) x_+(\tau) d\tau \ge I_{a,b}^+(t_0) = +\infty. \eqno (2.26)
$$
Let  $x_N(t)$ be an arbitrary $t_0$-normal solution of Eq.  (2.1).
By (2.9) we have:
$\ilp{t_0}|x_+(\tau) - x_N(\tau)| d\tau < +\infty$.  From here and from (2.26) we will get:  \linebreak $\ilp{t_0} a(\tau) x_N(\tau) d\tau = \ilp{t_0} a(\tau)(x_N(\tau) - x_+(\tau)) d\tau + \ilp{t_0} a(\tau) x_+(\tau) d\tau = +\infty$. The assertion IV$^\circ$ is proved. Prove V$^\circ$. Since on the strength of  II$^\circ \ph I_*(t) \equiv \il{t_0}{t} a(\tau) x_*(\tau) d\tau$ is a monotonically non increasing function on  $[t_0;+\infty)$, from Lemma 2.3 it follows (after differentiation (2.10)):
$$
a(t)x_0(t)\exp\biggl\{\il{t_0}{t}a(\xi)x_0(\xi) d\xi\biggr\}\ilp t \frac{a(s)x_0(s)}{x_0(t_0)}\exp\biggl\{\il{t_0}{s}\bigl[\frac{c(\xi)}{x_0(\xi)} - a(\xi) x_0(\xi)\bigr]d\xi\biggr\} d s \le \phantom{aaaaaaaaa}
$$
$$
 \phantom{aaaaaaaaaaaaaaaaaaaaaaaaaaa}\le \frac{a(t) x_0(t)}{x_0(t_0)}\exp\biggl\{\il{t_0}{t}\frac{c(\xi)}{x_0(\xi) d\xi}\biggr\}, \phh t\ge t_0. \eqno (2.27)
$$
where  $x_0(t) (> 0, \ph t\ge t_0)$ is a $t_0$-normal solution of Eq. (2.1) (by virtue of Lemma 2.1 from  I$^\circ$   it follows the existence of $x_0(t)$). Since  $u_0(t) \equiv \frac{1}{x_0(t)}$ is a   $t_0$-normal solution of Eq. (2.24) and   $I_{-c,-b}^+(t_0) = +\infty$, by IV$^\circ$ we have:
$$
\ilp{t_0}c(\tau)u_0(\tau) d\tau = \ilp{t_0}\frac{c(\tau)}{x_0(\tau)} d\tau = - \infty. \eqno (2.28)
$$
Since $a(t)$  has unbounded support there exists infinitely large sequence  $t_0 < t_1 < ... < t_m< ...$  such that  $a(t_m) > 0, \ph m=1, 2, ...$ . Then from (2.27) it follows
$$
\exp\biggl\{\il{t_0}{t_m}a(\xi)x_0(\xi) d\xi\biggr\}\ilp {t_m} \frac{a(s)x_0(s)}{x_0(t_0)}\exp\biggl\{\il{t_0}{s}\bigl[\frac{c(\xi)}{x_0(\xi)} - a(\xi) x_0(\xi)\bigr]d\xi\biggr\} d s \le \exp\biggl\{\il{t_0}{t_m}\frac{c(\xi)}{x_0(\xi)} d\xi\biggr\},
$$
$m=1,2, ...$ . Due to Lemma 2.3 From here and from (2.28) it folloes that  $I_*(t_m) \to - \infty$  for  $m \to +\infty$. Hence, $\ilp{t_0} a(\tau) x\*(\tau) d\tau = - \infty$.  The assertion V$^\circ$ is proved. Prove  VI$^\circ$.  Show that Eq. (2.1) has a  $t_0$-normal negative solution.
In Eq. (2.1) make the change: $x = J_{-b}(t_0;t) X, \ph t\ge t_0$.  We will come to the equation
$$
X' + a(t) J_{-b}(t) X^2 + c(t) J_b(t) = 0, \phh t\ge t_0. \eqno (2.29)
$$
Due to conditions of  VI$^\circ$  chose  $t_1 (> t_0)$ so large that  $$
\biggl[\ilp{t_1} a(\tau) J_{-b}(\tau) d\tau\biggr]^{-1} > - \ilp{t_1} c(\tau) J_{b}(\tau) d \tau.
$$
Then
$$
-\bigl[I_{a,b}^+(t_0)\bigr]^{-1} < I_{c,-b}^+(t_0) < 0. \eqno (2.30)
$$
Let then  $X_-(t)$ be a solution to Eq. (2.29) with
$$
X_-(t_1) \in \bigl(-\bigl[I_{a,b}^+(t_0)\bigr]^{-1}; I_{c,-b}^+(t_0)\bigr). \eqno (2.31)
$$
By (2.18) the inequalities
$$
\frac{X_-(t_1)}{1+ X_-(t_1) I_{a,b}^+(t_1;t)} \le X_-(t) \le X_-(t_1) - I_{c,-b}^+(t_1;t),\phh t\ge t_1, \eqno (2.32)
$$
are fulfilled.
From here and from (2.30) and  (2.31) it follows that  $X_-(t)$ is a defined on $[t_1;+\infty)$  negative    $t_1$-normal solution of Eq.  (2.29).  Then  $x_-(t)\equiv X_-(t)J_{-b}(t_1;t)$ is a defined on  $[t_1;+\infty)$ negative $t_1$-normal solution to Eq. (2.1). Show that $x_-(t)$ is continuable on $[t_0;+\infty)$  as a solution to Eq. (2.1). Suppose $x_-(t)$  can not be continued on  $[t_0;+\infty)$  as a solution of Eq. (2.1). Let then  $(t_2;+\infty)$ be the maximum existence interval for $x_-(t)$, where  $t_2 \ge t_0$. By Lemma 2.2  there exists $t_3 > t_2$ such that $x_-(t_3) > 0$. On the strength of the first of the inequalities (2.18) from here it follows that $x_-(t) > 0, \ph t\ge t_3$. The obtained contradiction shows that $x_-(t)$ is continuable on $[t_0;+\infty)$. By virtue of the first of the inequalities (2.18) the supposition that $x_-(t_4) \ge 0$ for some $t_4 \ge t_0$  also leads to the contradiction.
So, $x_-(t) < 0, \ph t\ge t_0$. Since  $x_-(t)$ is $t_1$-normal, by continuable  dependence of solutions of Eq. (2.1) from their  initial values, the solution  $x_-(t)$  also is $t_0$-normal.
According to  I$^\circ$ the solution $x_0(t)$  of Eq. (2.1) with $x_0(t_0) =0$ starting with  some  $t_1 = t_1(x_0) \ge t_6$ becomes positive.  Then by continuable dependence of solutions of Eq. (2.1) from their initial values, all initial values $x_{(0)}$, for which the solutions $x(t)$  of Eq. (2.1) with  $x(t_0) =x_{(0)}$  eventually become positive, form an open set. From here from Lemma 2.1 and from the fact that  $x_-(t)$ is negative it follows that
there exists a negative   $t_0$-normal solution  $X_N^-(t)$ of Eq. (2.1) such that any solution $x(t)$  of Eq.  (2.1) with   $x(t_0) > X_N^-(t_0)$  eventually becomes positive.  By (2.4) from here it follows that for any solution $x(t)$ of Eq. (2.1) with   $x(t_0) \in (x_N^-(t_0);0)$ there exists  $t_2 = t_2(x) \ge t_1=t_1(x) > t_0$ such that  $x(t) < 0, \ph t\in [t_0;t_1),\ph x(t) = 0, \ph t\in [t_1;t_2], \ph x(t) > 0, \ph t >t_2$.  By (2.8)\linebreak  $\ilp{t_0} a(\tau)[x_*(\tau) - x_N(\tau)]d\tau = - \infty$.  Then  $\ilp{t_0} a(\tau) x_*(\tau) d\tau \le \ilp{t_0} a(\tau) [x_*(\tau) - x_N(\tau)] d\tau =~-\infty$. Let $x_+(t)$ be a solution to Eq. (2.1) with  $x_+(t_0) =1$. On the strength of Lemma 2.1 from I$^\circ$  it follows that $x_+(t)$ is $t_0$-normal. Then by (2.9) we have:
$$
0 < \ilp{t_0} a(\tau) x_+(\tau) d\tau \le \ilp{t_0} a(\tau) [x_+(\tau) - x_N^-(\tau)]d\tau < +\infty. \eqno (2.33)
$$
Let $x_N(t)$ be an arbitrary $t_0$-normal solution to Eq. (2.1). Then since \linebreak
$\ilp{t_0}a(\tau)x_N(\tau) d\tau = \ilp{t_0}[x_N(\tau) - x_+(\tau)]d\tau + \ilp{t_0} a(\tau) x_+(\tau)d\tau$
and by virtue of (2.9)\linebreak $\ilp{t_0}a(\tau)|x_N(\tau) - x_+(\tau)|d\tau < +\infty$,  from  (2.33) it follows that the integral  $\ilp{t_0}a(\tau)x_N(\tau)d\tau$ converges. The assertion VI$^\circ$ is proved.   Prove VII$^\circ$. Due to the second of inequalities (2.18) taking into account the inequality   $c(t) \le 0, \ph t\ge t_0$ we have:
$$
\il{t_0}{t}\frac{c(\tau)}{x_0(\tau)} d\tau \ge \frac{1}{x_0(t_0)}\il{t_0}{t} c(\tau) J_b(t)\biggl[1 + x_0(t_0) I_{a,b}^+(t_0;\tau)\biggr] d\tau, \phh t\ge t_0, \eqno (2.34)
$$
where  $x_0(t)$ is a positive $t_0$-normal solution of Eq. (2.1), existence of which follows from Lemma 2.1 and from  I$^\circ$. By virtue of Fubini's theorem from the second relation of  (2.19) it follows that  $I_{-b,-c}^+(t_0;+\infty) < +\infty.$ From here and from the second relation of (2.19) and from  (2.34) it follows that
$$
\ilp{t_0}\frac{c(\tau)}{x_0(\tau)} d\tau > - \infty. \eqno (3.35)
$$
From the first relations of (2.18) and  (2.19) it follows that
$$
\ilp{t_0} a(\tau) x_0(\tau) d\tau = +\infty. \eqno (2.36)
$$
Set  $g(t)\equiv \exp\biggl\{-\il{t_0}{t} a(\xi) x_0(\xi) d\xi\biggr\}$.  Obviously the inverse function $g^{-1}(t)$ of  $g(t)$
exists on $supp \ph a(t)$. Denote:
$$
g_1(t) \equiv \left\{
\begin{array}{l}
g^{-1}(t), \ph t\in supp \ph a(t);\\
\ph\\
t_0,\phh t \not \in supp \ph a(t), \ph t\ge t_0.
\end{array}
\right.
$$
Then taking into account (2.36) the equality (2.16) can be rewritten in the form
$$
\il{t_0}{t} a(\tau)
 x_*(\tau) d\tau =
 - \ln \nu_{x_0}(t_0) + \ln
\Biggl[\frac{\il{0}{g(t)}\exp\biggl\{\il{t_0}{g^{-1}(\zeta)}\frac{c(\xi)}{x_0(\xi)}d\xi\biggr\}d\zeta}{x_0(t_0) g(t)}\Biggr]\ge 1, \phh t\ge t_0.
$$
From here and from (2.35) it follows that
$$
\ilp{t_0} a(\tau) x_*(\tau) d\tau \ge - \ln [x_0(t_0) \nu_{x_0}(t_0)] + \ilp{t_0} \frac{c(\xi)}{x_0(\xi)}d\xi > - \infty.
$$
By virtue of Lemma 2.1  $u_0(t)\equiv \frac{1}{x_0(t)}$ is a $t_0$-normal solution of Eq. (2.24)
(since for  $x_1(t_0) > x_0(t_0)$  the function $u_1(t) \equiv \frac{1}{x_1(t)}$ is a  $t_0$-regular solution of Eq. (2.24), where $x_1(t)$ is a  $t_0$-regular solution to Eq. (2.1)). Then by  (2.8) taking into account (2.35)  we obtain:
$$
\ilp{t_0} c(\tau) u_*(\tau) d\tau = \ilp{t_0} c(\tau)(u_*(\tau) - u_0(\tau)) d\tau + \ilp{t_0}\frac{c(\tau)}{x_0(\tau)} d\tau = +\infty.
$$
The assertion VII$^\circ$ is proved. The theorem is proved.

On the basis of Theorem 2.3 it can be make phase portrait of solutions of Eq. (2.1) in the case  $a(t) \ge 0, \ph c(t) \le 0, \ph t\ge t_0,$  for the following two possible restrictions:

\noindent
a) $I_{a,b}^+(t_0) = +\infty$    or  $I_{c,-b}^+(t_0) = -\infty$ (see  pict. 1) ;

\noindent
b) $I_{a,b}^+(t_0) < +\infty$    and  $I_{-c,-b}^+(t_0) < +\infty$   (see   pict. 2).

\noindent
In pict 1. we see only one negative global solution of Eq (2.1), meanwhile in pict. 2 we see whole slice of negative global solutions of Eq. (2.1).

\begin{picture}(40,100)
\put(-10,-40){\vector(0,1){120}}
\put(-20,20){\vector(1,0){220}}
\put(230,-40){\vector(0,1){120}}
\put(220,20){\vector(1,0){220}}

\put(188,13){$_t$}
\put(-5,75){$_{x(t)}$}
\put(15,13){$_{t_0}$}
\put(22,20){$_\circ$}
\put(158,-16){$_{x_*(t)}$}
\put(-10,-55){$_{pict.1.}$}
\put(20,-55){$_{I_{a,b}^+(t_0)=+\infty}$ $_{or}$ $_{I_{-c,-b}^+(t_0)=+\infty}$}
\put(230,-55){$_{pict.2.}$}
\put(260,-55){$_{I_{a,b}^+(t_0)<+\infty}$ $_{and}$ $_{I_{-c,-b}^+(t_0)<+\infty}$}

\put(428,13){$_t$}
\put(245,13){$_{t_0}$}
\put(388,-26){$_{x_*(t)}$}
\put(388,4){$_{x_N^-(t)}$}
\put(235,75){$_{x(t)}$}
\put(252,20){$_\circ$}

\put(0,-20){\thicklines\qbezier[30](24,-20)(24,40)(24,100)}
\put(44,-20){\thicklines\qbezier[30](24,-20)(24,40)(24,100)}
\put(84,-20){\thicklines\qbezier[30](24,-20)(24,40)(24,100)}

\put(230,-20){\thicklines\qbezier[30](24,-20)(24,40)(24,100)}
\put(274,-20){\thicklines\qbezier[30](24,-20)(24,40)(24,100)}
\put(314,-20){\thicklines\qbezier[30](24,-20)(24,40)(24,100)}

\put(20,-5){\thicklines \gr{0}{10}{10}{14}{10}{6}{10}{15}{10} \grik{-42}{10}{10}{5}{10} \gr{-45}{10}{10}{15}{10}{5}{10}{15}{10}
\grik{-88}{10}{10}{5}{10}}

\put(20,30){\gr{0}{10}{10}{14}{10}{6}{10}{15}{10} \grik{-42}{10}{10}{7}{18} \gr{-45}{10}{18}{23}{18}{12}{18}{23}{18}
\grik{-88}{10}{18}{16}{18}}

\put(20,28){\gr{0}{10}{10}{14}{10}{6}{10}{15}{10} \grik{-42}{10}{10}{7}{18} \gr{-45}{10}{18}{23}{18}{12}{18}{23}{18}
\grik{-88}{10}{18}{16}{18}}

\put(20,26){\gr{0}{10}{10}{14}{10}{6}{10}{15}{10} \grik{-42}{10}{10}{7}{18} \gr{-45}{10}{18}{23}{18}{12}{18}{23}{18}
\grik{-88}{10}{18}{16}{18}}

\put(20,24){\gr{0}{10}{10}{14}{10}{6}{10}{15}{10} \grik{-42}{10}{10}{7}{18} \gr{-45}{10}{18}{23}{18}{12}{18}{23}{18}
\grik{-88}{10}{18}{16}{18}}

\put(20,-2){\gr{0}{10}{10}{13}{10}{7}{10}{16}{13} \grik{-42}{10}{13}{10}{19} \gr{-45}{10}{19}{23}{20}{18}{30}{42}{36}
\grik{-88}{10}{35}{33}{35}}

\put(20,1){\gr{0}{10}{10}{13}{10}{7}{10}{16}{13} \grik{-42}{10}{13}{10}{19} \gr{-45}{10}{19}{23}{20}{18}{30}{42}{36}
\grik{-88}{10}{35}{33}{35}}

\put(20,4){\gr{0}{10}{10}{13}{10}{7}{10}{16}{13} \grik{-42}{10}{13}{10}{19} \gr{-45}{10}{19}{23}{20}{18}{30}{42}{36}
\grik{-88}{10}{35}{33}{35}}

\put(20,-14){\gr{0}{10}{15}{20}{16}{12}{15}{18}{14}}
\put(85, -12){\qbezier(0,2)(5,-2)(11,0)}
\put(87,-12){\qbezier(10,0)(18,-0.1)(19,-25)}

\put(70,10){\qbezier(0,70)(3,50)(10,45)}
\put(75,25){\graph{30}{20}{27}{33}{30}{25}{27}{30}{27}}
\put(159,52){\qbezier(0,0)(15,-5)(25,0)}

\put(110,30){\qbezier(0,50)(1,35)(14,30)}
\put(170,83){\gr{50}{50}{27}{23}{30}{35}{27}{20}{27}}

\put(25,-14){\qbezier(0,0)(7,4)(14,0)}
\put(35,-14){\qbezier(5,0)(7,-7)(20,-3)}

\put(45,-17){\qbezier(10,0)(18,-0.1)(20,-20)}

\put(27,80){\qbezier(0,0)(2,-35)(10,-37)}
\put(33,40){\gr{0}{10}{13}{10}{13}{17}{10}{8}{13}}
\put(120,40){\gr{0}{10}{13}{10}{13}{17}{10}{8}{13}}
\put(100,43){\qbezier(-2,0)(10,4)(23,0)}

\put(250,-15){\thicklines \gr{0}{10}{10}{15}{10}{5}{10}{15}{10} \grik{-42}{10}{10}{5}{10} \gr{-45}{10}{10}{15}{10}{5}{10}{15}{10}
\grik{-88}{10}{10}{5}{10}}

\put(250,-5){\thicklines \gr{0}{10}{10}{15}{10}{6}{10}{15}{10} \grik{-42}{10}{10}{5}{10} \gr{-45}{10}{10}{15}{10}{5}{10}{15}{10}
\grik{-88}{10}{10}{5}{10}}

\put(250,-7){\thinlines \gr{0}{10}{10}{15}{10}{6}{10}{15}{10} \grik{-42}{10}{10}{5}{10} \gr{-45}{10}{10}{15}{10}{5}{10}{15}{10}
\grik{-88}{10}{10}{5}{10}}

\put(250,-9){\thinlines \gr{0}{10}{10}{15}{10}{6}{10}{15}{10} \grik{-42}{10}{10}{5}{10} \gr{-45}{10}{10}{15}{10}{5}{10}{15}{10}
\grik{-88}{10}{10}{5}{10}}

\put(250,-12){\thinlines \gr{0}{10}{10}{15}{10}{6}{10}{15}{10} \grik{-42}{10}{10}{5}{10} \gr{-45}{10}{10}{15}{10}{5}{10}{15}{10}
\grik{-88}{10}{10}{5}{10}}

\put(250,30){\gr{0}{10}{10}{14}{10}{6}{10}{15}{10} \grik{-42}{10}{10}{7}{18} \gr{-45}{10}{18}{23}{18}{12}{18}{23}{18}
\grik{-88}{10}{18}{16}{18}}

\put(250,28){\gr{0}{10}{10}{14}{10}{6}{10}{15}{10} \grik{-42}{10}{10}{7}{18} \gr{-45}{10}{18}{23}{18}{12}{18}{23}{18}
\grik{-88}{10}{18}{16}{18}}

\put(250,32){\gr{0}{10}{10}{14}{10}{6}{10}{15}{10} \grik{-42}{10}{10}{7}{18} \gr{-45}{10}{18}{23}{18}{12}{18}{23}{18}
\grik{-88}{10}{18}{16}{18}}

\put(250,-2){\gr{0}{10}{10}{13}{10}{7}{10}{16}{13} \grik{-42}{10}{13}{10}{14} \gr{-45}{10}{14}{15}{17}{25}{40}{45}{40}
\grik{-88}{10}{40}{33}{38}}

\put(250,1){\gr{0}{10}{10}{13}{10}{7}{10}{16}{13} \grik{-42}{10}{13}{10}{14} \gr{-45}{10}{14}{15}{17}{25}{40}{45}{40}
\grik{-88}{10}{40}{33}{38}}

\put(250,4){\gr{0}{10}{10}{13}{10}{7}{10}{16}{13} \grik{-42}{10}{13}{10}{14} \gr{-45}{10}{14}{15}{17}{25}{40}{45}{40}
\grik{-88}{10}{40}{33}{38}}

\put(250,-24){\gr{0}{10}{15}{20}{16}{12}{15}{18}{14}}
\put(315, -22){\qbezier(0,2)(5,-2)(11,0)}
\put(317,-22){\qbezier(10,0)(18,-0.1)(19,-25)}

\put(300,10){\qbezier(0,70)(3,50)(10,45)}
\put(305,25){\graph{30}{20}{27}{33}{30}{25}{27}{30}{27}}
\put(389,52){\qbezier(0,0)(15,-5)(25,0)}

\put(340,30){\qbezier(0,50)(1,35)(14,30)}
\put(400,83){\gr{50}{50}{27}{23}{30}{35}{27}{20}{27}}

\put(255,-24){\qbezier(0,0)(7,4)(14,0)}
\put(265,-24){\qbezier(5,0)(7,-7)(20,-3)}
\put(275,-27){\qbezier(10,0)(18,-0.1)(20,-20)}

\put(257,83){\qbezier(0,0)(2,-35)(10,-37)}
\put(263,43){\gr{0}{10}{13}{10}{13}{17}{10}{8}{13}}
\put(350,43){\gr{0}{10}{13}{10}{13}{17}{10}{8}{13}}
\put(330,46){\qbezier(-2,0)(10,4)(23,0)}

\end{picture}

\vskip 90pt

{\bf Theorem 2.4.} {\it Suppose  $a(t) \ge 0, \ph c(t) \ge 0, \ph t\ge t_0$, and  Eq. (2.1) has a   $t_0$-regular solution. Then the following assertions are valid.

\noindent
I$^*$  If $I_{a,b}^+(t_0) = +\infty$, then every $t_0$-regular solution of Eq. (2.1) is positive and for its every    $t_0$-normal solution  $x_N(t)$
the equality  $\ilpp a(\tau)x_N(\tau) d\tau = +\infty$ is fulfilled.
Moreover if in addition $\ilpp c(\tau) I_{a,-b}^+(t_0;\tau)  d\tau = +\infty$, then also  $\ilpp a(\tau) x_*(\tau) d \tau = + \infty$,  where  $x_*(t)$ is the unique $t_0$-extremal solution of Eq. (2.1).

\noindent
II$^*$ If $I_{c,-b}^+(t_0) = +\infty$, then for every  $t_0$ -regular solution $x(t)$  of Eq. (2.1) with  $x(t_0) > 0$ there  exist $t_2=t_2(x) \ge t_1 = t_1(x) > t_0$ such that  $x(t) > 0, \ph t\in [t_0;t_1), \ph x(t) = 0, \ph t\in [t_1;t_2], \ph x(t) < 0, \ph t > t_2,$ and if  $x(t_0) = 0 (< 0)$,
there exists $t_1 = t_1(x) \ge t_0$ such that  $x(t) = 0, \ph t\in [t_0;t_1], \ph x(t) < 0, \ph t > t_1$ (then  $x(t) < 0, \ph t\ge t_0$), and  $\ilpp a(\tau)x_*(\tau) d\tau = -\infty$. Moreover if in addition  $\ilpp a(\tau) I_{b,c}^-(t_0;\tau)d\tau = +\infty$, then for every  $t_0$-normal solution  $x_N(t)$ of Eq. (2.1) the equality $\ilpp a(\tau) x_N(\tau) d\tau = -\infty$ is fulfilled.

\noindent
III$^*$ If  $I_{a,b}^+(t_0) < +\infty$  and  $I_{c,-b}^+(t_0) < +\infty$, then there exist  $t_1 \ge t_0$  such that $x_*(t) < 0, \ph t\ge t_1$; the solutions  $x(t)$  of Eq. (2.1) with  $x(t_1) \in (x_*(t_1);0)$    are $t_0$-regular and  $x(t) < 0, \ph t\ge t_1$; there exists a $t_1$-normal positive on $[t_1;+\infty)$   solution  $x_N^+(t)$ of Eq. (2.1) such that fir every  solution $x(t)$ of Eq. (2.1) with   $x(t_1) \in (0;x_N^+(t_1))$  there exist $t_3 = t_3(x) \ge t_2 = t_2(x) > t_1$  such that $x(t) > 0, \ph t\in[t_1;t_2), \ph x(t) = 0, \ph t\in [t_2; t_3], \ph x(t) , 0, \ph t\ < t_3$,  and if   $x(t_1) =0 \ph (x(t_1) < 0)$, there exists $t_2=t_2(x) \ge t_1$ such that  $x(t) = 0, \ph t\in [t_1;t_2], \ph x(t) < 0, \ph t > t_2$  (then  $x(t) < 0, \ph t\ge t_1$); for every
$t_0$-normal solution $x_N(t)$ of Eq. (2.1) the integral $\ilpp a(\tau)x_N(\tau) d\tau$ converges and  $\ilpp a(\tau)x_*(\tau) d\tau = -\infty$.
}

Proof. Let us prove I$^*$.  Let $x_*(t)$  be the $t_0$-extremal solution of Eq. (2.1). Show that  $x_*(t) > 0, \ph t\ge t_0$. Suppose for some  $t_1 > t_0$ the inequality  $x_*(t_1) < 0$ is satisfied.  Let then  $x(t)$  be a solution to Eq. (2.1) with  $x(t_1) \in (x_*(t_1);0)$. By virtue of Lemma 2.1 $x(t)$ is $t_1$-normal. Since  $x(t_1) < 0$    and $c(t) \ge 0, \ph t\ge t_0$, by (2.4) we have $x(t) < 0, \ph t\ge t_1$. From here it follows that
$$
\nu_x(t_1) \ge I_{a,b}^+(t_1). \eqno (2.37)
$$
Suppose  $I_{a,b}^+(t_0)= +\infty$. Then from the easily verifiable equality
$$
I_{a,b}^+(t_0) = I_{a,b}^+(t_0;t) + J_b(t_1) I_{a,b}^+(t_1)   \eqno (2.38)
$$
and from (2.37) it follows that $\nu_x(t_1) = +\infty$. But on the other hand since $x(t)$ is $t_1$-normal by virtue of Theorem 2.2 we have $\nu_x(t_1) , +\infty$.
The obtained contradiction shows that $x_*(t) \ge 0, \ph t\ge t_0$. Show that the equality $x_*(t) =0$ is impossible for all $t\ge t_0$.  Suppose for some  $t_2 \ge t_0$
the equality  $x_*(t_2) = 0$ is satisfied.  Then by (2.4) from the inequality  $c(t) \ge 0, \ph t\ge t_0$       it follows that $x_*(t) \le 0, \ph t\ge t_2$.  Hence, $x_*(t) \equiv 0$ on  $[t_2;+\infty)$, which is impossible  (since on  $[t_2;+\infty)$ we have $c(t) \not \equiv 0$). On the strength of Lemma 2.1 from here it follows that every  $t_0$-regular solution of Eq. (2.1) is positive. Let $x_N(t)$ be a  $t_0$-normal solution of Eq. (2.1). Then since  $x_*(t) > 0, \ph t\ge t_0$, by  (2.8) we have:
$
\ilpp a(\tau) x_N(\tau) d\tau =\linebreak = \ilpp a(\tau)[x_N(\tau) - x_*(\tau)]d\tau + \ilpp a(\tau) x_*(\tau) d\tau \ge \ilpp a(\tau)[x_N(\tau) - x_*(\tau)]d\tau = +\infty.
$
Suppose
$$
\ilpp c(\tau) I_{a,b}(t_0;\tau) d\tau = +\infty, \eqno (2.39)
$$
and  $\ilpp a(\tau)x_*(\tau) d\tau < +\infty$. Show that then
$$
x_*(t) = \ilp{t} J_b(t;\tau) c(\tau)\phi_*(t;\tau)d\tau, \phh t\ge t_0, \eqno (2.40)
$$
where  $\phi_*(t;\tau)\equiv \exp\bigl\{\il{t}{\tau} a(s) x_*(s) d s\bigr\}, \ph \tau \ge t \ge t_0$.                    By  (2.4) we have:
$$
x_*(t) = J_{-h_*}(t_1;t)\biggl[x_*(t_1) - \il{t_1}{t}J_b(t_1;\tau)c(\tau)\phi_*(t;\tau) d\tau\biggr], \phh t\ge t_1\ge t_0, \eqno (2.41)
$$
where  $h_*(t)\equiv a(t)x_*(t) + b(t), \ph t\ge t_0$. From here and from the positivity of $x_*(t)$ it follows that
$$
x_*(t) \ge \ilp{t_1} j_b(t_1;\tau) c(\tau) \phi_*(t_1;\tau) d\tau, \phh t_1 \ge t_0.
$$
Show that the strict inequality
$$
x_*(t) > \ilp{t_1} J_b(t_1;\tau) c(\tau) \phi_*(t_1;\tau) d\tau  \eqno (2.42)
$$
is impossible for all $t_1 \ge t_0$.  Multiplying both sides of (2.41) on  $a(t)\phi_*(t_1;t)$ and integrationg from $t_1$  to  $+\infty$ we  get:
$$
\exp\biggl\{\ilp{t_1}a(\tau)x_*(\tau) d\tau\biggr\} = 1 + \ilp{t_1} J_{-b}(t_1;\tau)\biggl[x_*(t_1) -\il{t_1}{t} J_b(t_1;\tau) c(\tau)\phi_*(t_1;\tau)\biggr] d\tau \ge
$$
$$
1 + I_{a,b}^+(t_1)\biggl[x_*(t_1) - \ilp{t_1}J_b(t_1;\tau)c(\tau)\phi_*(t_1;\tau)d\tau\biggr]. \eqno (2.43)
$$
Suppose for some $t_1 \ge t_0$ the inequality (2.42) is satisfied. Then by (2.38) from the equality  $I_{a,b}^+(t_0)= +\infty$ and from (2.43) it follows that $\ilp{t_1}a(\tau)x_*(\tau)d\tau = +\infty$.                 The obtained contradiction proves (2.40). Multiplying both sides of (2.40) on $a(t)\phi_*(t_0;t)$ and integrating from $t_0$ to  $+\infty$ we will get:
$
\exp\biggl\{\ilpp a(\tau)x_*(\tau) d\tau\biggr\} = 1 +\linebreak \ilpp a(\tau)\phi_*(t_0;\tau)d\tau\ilp{t} J_b(t;\tau) c(\tau)\phi_*(t;\tau) d\tau \ge 1 + \ilpp a(t) d t\ilp{t} J_b(t;\tau) c(\tau) d\tau.
$
By virtue of Fubini's theorem from here and from (2.39) it follows that  $\ilpp a(\tau) x_*(\tau)d\tau = +\infty$.  We came to the contradiction. The assertion  I$^*$ is proved. Let us prove II$^*$.  Let  $x(t)$   be a  $t_0$-regular solution of Eq. (2.1). If  $x(t_0) < 0$, then by (2.4) from inequality  $c(t) \ge0, \ph t\ge t_0$ it follows that  $x(t) < 0, \ph t\ge t_0$. Let  $x(t_0) \ge 0$. Show that in this case it is impossible that
$$
x(t) \ge 0, \phh t\ge t_0. \eqno (2.44)
$$
Suppose  this relation is true. Then since  $a(t) \ge 0, \ph c(t) \ge 0, \ph t\ge t_)$, we have
$$
\il{t_0}{t} J_b(\tau) c(\tau) \exp\biggl\{\il{t_0}{\tau} a(s) x(s) d s\biggr\}d\tau \ge I_{c,-b}^+(t_0;t), \phh t\ge t_0.
$$
By (2.4) from here and from equality  $I_{c,-b}^+(t_0) = +\infty$ it follows that $x(t_1) < 0$ for some  $t_1 > t_0$. We came to the contradiction. Hence by (2.4) if  $x(t_0) > 0,$
then there exists  $t_2 = t_2(x)\ge t_1= t_1(x) > t_0$ such that  $x(t) > 0, \ph t\in [t_0;t_1), \ph x(t) = 0, \ph t\in [t_1;t_2], \ph x(t) < 0, \ph t > t_2$,
and if  $x(t_0) = 0$, then there exists  $t_1 = t_1(x) \ge t_0$  such that  $x(t) = 0, \ph t\in [t_0;t_1]$, and  $x(t) < 0, \ph t\ge t_1$. Let $x_0(t)$  be a  $t_0$-normal solution of Eq. (2.1) with  $x_0(t_0) \ge 0$, and let  $t_1 = t_1(x) \ge t_0$ such that  $x_0(t_1) = 0, \ph x_0(t) < 0, \ph t\>t_1$. Then since by (2.8) $\ilp{t_1} a(\tau)[x_*(\tau) - x_0(\tau)] d\tau = -\infty$, we have
$$
\ilpp a(\tau) x_*(\tau) d\tau = \il{t_0}{t_1}
a(\tau) x_*(\tau) d\tau + \ilp{t_1} a(\tau)[x_*(\tau) - x_0(\tau)] d\tau + \ilp{t_1}a(\tau)x_0(\tau) d\tau \le \phantom{aaaaaa}
$$
$$
\phantom{aaaaaaaaaaaaaaaaaaaaaaaaaaaaaaa}\le \il{t_0}{t} a(\tau) x_*(\tau) d\tau + \ilp{t_1} a(\tau) [x_*(\tau) - x_0(\tau)]d\tau  = -\infty.
$$
Using Theorem  2.1  by analogy of the second of inequalities  (2.18) can be obtained the estimation
$$
x_0(t) \le x_0(t_1) J _{-b}(t_1;t) - I_{b,c}^-(t_1;t) = - I_{b,c}^-(t_1;t), \phh t\ge t_1.
$$
Then
$$
\ilpp a(\tau) x_0(\tau) d\tau \le \il{t_0}{t_1} a(\tau) x_0(\tau) d\tau - \ilp{t_1} a(\tau) I_{a,b}^-(t_1;\tau) d\tau. \eqno (2.45)
$$
Since   $I_{b,c}^-(t_0;t) = I_{b,c}^-(t_0;t_1) J_{-b}(t_1;t) + I_{b,c}^-(t_1;t)$,  we have
$$
\ilpp a(\tau) I_{b,c}^-(t_0;\tau)d\tau = \il{t_0}{t_1} I_{b,c}^-(t_1;\tau) d\tau +  I_{b,c}^-(t_0;t_1)I_{a,b}^+(t_1;+\infty) + \ilp{t_1} a(\tau)I_{b,c}^-(t_1;\tau) d\tau. \eqno (2.46)
$$
Since  $x_0(t) < 0, \ph t > t_1$, by virtue of  I$^*$ we obtain
$$
I_{a,b}^+(t_1) < +\infty. \eqno (2.47)
$$
Let  $\ilpp a(\tau) I_{b,c}^-(t_0;\tau) d\tau = +\infty$. Then from (2.45) -  (2.47) it follows that
$$
\ilpp a(\tau) x_0(\tau) d\tau = - \infty. \eqno (2.48)
$$
Let  $x_N(t)$ be an arbitrary  $t_0$-normal solution of Eq. (2.1). Then since by (2.9) $\ilpp a(\tau)|x_N(\tau) - x_0(\tau)| d\tau < +\infty$, taking into account (2.48) we  have:  $\ilpp a(\tau) x_N(\tau) d\tau = \ilpp a(\tau)[x_N(\tau) - x_0(\tau)]d\tau + \ilpp a(\tau) x_0(\tau) d\tau = -\infty$.                               The assertion II$^*$  is proved. Let us prove  III$^*$. Show that
$$
x_*(t_1) < 0  \eqno (2.49)
$$
for some  $t_1 \ge t_0$. Suppose that it is not true. Then  $x_*(t) \ge 0, \ph t\ge t_)$ and therefore   $\nu_{x_*}(t_0) \le I_{a,b}^+(t_0) < +\infty.$  But on the other hand by (2.7) we have $\nu_{x_*}(t_0) = +\infty$. The obtained contradiction proves (2.49). By (2.4) from (2.49) and from non negativity of  $c(t)$ it follows that
$$
x_*(t) < 0, \ph t\ge t_1. \eqno (2.50)
$$
Hence  $v_*(t)\equiv \frac{1}{x_*(t)}, \ph t\ge t_1$ is a  $t_1$-regular solution of the equation
$$
v' + c(t) v^2 - b(t) v + a(t) = 0, \phh t\ge t_0. \eqno (2.51)
$$
Suppose  $I_{c,-b}^+(t_0) < +\infty$. By (2.38) from here it follows that  $I_{c,-b}^+(t_1) < +\infty$.
Then by already proven   $v_*(t) < 0, \ph t\ge t_2$    for some  $t_2\ge t_1$, where  $v_*(t)$ is the   $t_1$-extremal solution of Eq. (2.51).  From here it follows that  $x_1(t) \equiv - \frac{1}{v_*(t)}$            is an positive solution of Eq. (2.1) defined on  $[t_2; +\infty)$. Since according to (2.50) $x_*(t_2) < 0$, by virtue of Lemma 2.1 $x_1(t)$  is   $t_2$-normal. By virtue of continuously dependence of solutions of Eq. (2.1) from their initial values from here from (2.4) and  (2.50) it follows that there exists $t_1$-normal positive solution  $x_N^+(t)$ of Eq. (2.1) on $[t_2;+\infty)$ having the property: for each solution  $x(t)$   of Eq. (2.1) with  $x(t_2) \in (0;x_N^+(t_2))$  there exists  $t_4 = t_4(x) \ge t_3= t_3(x) > t_2$ such that $x(t) > 0, \ph t\in [t_2;t_3), \ph x(t) = 0, \ph t\in [t_3;t_4], \ph x(t)< 0, \ph t > t_4$; if  $x(t_2) = 0$, then there exists  $t_3 = t_3(x) \ge t_2$ such that $x(t) = 0, \ph t\in [t_2;t_3], \ph x(t) < 0, \ph t> t_3$;  if $x(t_2) \in (x_*(t_2);0)$, then $x(t) < 0, \ph t\ge t_2$.  Let  $x_-(t)$ be a solution of Eq. (2.1) with  $x_-(t_2) \in (x_*(t_2);0)$. Then by virtue of Lemma 2.1  $x_-(t)$ is $ t_2$-normal and, as it was already proved,  $x_-(t) < 0, \ph t\ge t_2$. Therefore taking into account (2.9) we obtain:
$$
0 < \ilpp a(\tau) x_N^+(\tau) d\tau  \le \il {t_0}{t_2}a(\tau) x_N^+(\tau) d\tau + \ilp{t_2}a(\tau)[ x_N^+(\tau) - x_-(\tau)] d\tau < +\infty. \eqno (2.52)
$$
Let $x_N(t)$ be an arbitrary  $t_0$-normal solution of Eq. (2.1). Then \linebreak $\ilpp a(\tau) x_N(\tau) d\tau = \il {t_0}{t_2} a(\tau) x_N(\tau) d\tau + \ilp{t_2}a(\tau)[x_N(\tau) - x_N^+(\tau)] d\tau + \ilp{t_2} a(\tau)x_N^+(\tau) d\tau.$
Since by (2.9) $\ilp{t_2} a(\tau)|x_N(\tau) - x_N^+(\tau)| d\tau < +\infty$,  from the last equality and from (2.52) it follows  convergence of the integral  $\ilpp a(\tau) x_N(\tau) d\tau $.  Since  $x_N^+(t) > 0, \ph t\ge t_2; \ph \ilpp a(\tau) x_*(\tau) d\tau = \il {t_0}{t_2} a(\tau) x_*(\tau) d\tau + \ilp{t_2}a(\tau)[x_*(\tau) - x_N^+(\tau)] d\tau + \ilp{t_2} a(\tau)x_N^+(\tau) d\tau.$       and by  (2.8)  $\ilp{t_2} a(\tau)[x_*(\tau) - x_N^+(\tau)] d\tau =-\infty$,
taking into account (2.52) we  have:  $\ilpp a(\tau)x_*(\tau) d\tau = -\infty$.  The theorem is proved.

{\bf Remark 2.2}. {\it  Existence criteria of  $t_1$-regular solutions of Eq. (2.1) are proved in  [5] and  [12].}

{\bf Remark 2.3.} {\it If  $a(t) > 0, \ph t\ge t_)$,  then existence of  $t_1$-regular solutions of Eq.  (2.1) is equivalent to the non oscillation of the equation
$$
\biggl(\frac{\phi'}{a(t)}\biggr)' - a(t) b(t) \phi' - c(t) \phi = 0, \phh t\ge t_0.
$$
Non oscillatory criteria for the last equation is proved in [13].}

{\bf Corollary 2.2.} {\it Let  $a(t) \ge 0, \ph c(t) \ge 0, \ph t\ge t_0, \ph I_{a,b}^+(t_0) = I_{c,-b}(t_0) =  +\infty$. Then Eq. (2.1) has no $t_1$-regular solutions for all  $t_1 \ge t_0$.}

Proof. Suppose that for some  $t_1 \ge t_)$ Eq. (2.1) has  $t_1$-regular solution $x(t)$. Then by virtue of Theorem 2.4  from the equality $I_{a,b}^+(t_1) = I_{a,b}^+(t_0)$ (see (2.38)) it follows that  $x(t) > 0, \ph t\ge t_1$. Therefore  $v(t) \equiv -\frac{1}{x(t)}, \ph t\ge t_1$, is a negative  $t_1$-regular solution of Eq. (2.51). But on the other hand by virtue of Theorem 2.4 I$^*$ from the equality $I_{c,-b}^+(t_1) = I_{a,b}^+(t_0) = +\infty$ it follows that  $v(t) > 0, \ph t\ge t_1$.  We came to the contradiction. The corollary is proved.

On the basis of Theorem 2.1 and corollary 2.2 we can make the phase portrait of solutions of Eq. (2.1) if $a(t) \ge 0, \ph c(t) \ge 0, \ph t\ge t_0$ for the following four  cases:

\noindent
$\alpha) \ph I_{a,b}^+(t_0) = +\infty$;  and Eq. (2.1) has a $t_1$-regular solution for some  $t_1 \ge t_0$ (see pict.  3);

\begin{picture}(40,100)
\put(-10,-40){\vector(0,1){120}}
\put(-20,-20){\vector(1,0){220}}
\put(230,-40){\vector(0,1){120}}
\put(220,20){\vector(1,0){220}}

\put(188,-27){$_t$}
\put(-5,75){$_{x(t)}$}
\put(15,-27){$_{t_0}$}
\put(22,-20){$_\circ$}
\put(158,-2){$_{x_*(t)}$}
\put(-10,-55){$_{pict.3.}$}
\put(20,-55){$_{I_{a,b}^+(t_0)=+\infty}$}
\put(230,-55){$_{pict.4.}$}
\put(260,-55){$_{I_{c,-b}^+(t_0)=+\infty}$}

\put(428,13){$_t$}
\put(245,13){$_{t_0}$}
\put(388,-26){$_{x_*(t)}$}
\put(235,75){$_{x(t)}$}
\put(252,20){$_\circ$}

\put(0,-20){\thicklines\qbezier[30](24,-20)(24,40)(24,100)}
\put(44,-20){\thicklines\qbezier[30](24,-20)(24,40)(24,100)}
\put(84,-20){\thicklines\qbezier[30](24,-20)(24,40)(24,100)}

\put(230,-20){\thicklines\qbezier[30](24,-20)(24,40)(24,100)}
\put(274,-20){\thicklines\qbezier[30](24,-20)(24,40)(24,100)}
\put(314,-20){\thicklines\qbezier[30](24,-20)(24,40)(24,100)}

\put(20,5){\thicklines \gr{0}{10}{10}{14}{10}{6}{10}{15}{10} \grik{-42}{10}{10}{5}{10} \gr{-45}{10}{10}{15}{10}{5}{10}{15}{10}
\grik{-88}{10}{10}{5}{10}}

\put(20,30){\gr{0}{10}{10}{14}{10}{6}{10}{15}{10} \grik{-42}{10}{10}{7}{18} \gr{-45}{10}{18}{23}{18}{12}{18}{23}{18}
\grik{-88}{10}{18}{16}{18}}

\put(20,26){\gr{0}{10}{10}{14}{10}{6}{10}{15}{10} \grik{-42}{10}{10}{7}{18} \gr{-45}{10}{18}{23}{18}{12}{18}{23}{18}
\grik{-88}{10}{18}{16}{18}}

\put(20,22){\gr{0}{10}{10}{14}{10}{6}{10}{15}{10} \grik{-42}{10}{10}{7}{18} \gr{-45}{10}{18}{23}{18}{12}{18}{23}{18}
\grik{-88}{10}{18}{16}{18}}

\put(20,18){\gr{0}{10}{10}{14}{10}{6}{10}{15}{10} \grik{-42}{10}{10}{7}{18} \gr{-45}{10}{18}{23}{18}{12}{18}{23}{18}
\grik{-88}{10}{18}{16}{18}}

\put(20,14){\gr{0}{10}{10}{14}{10}{6}{10}{15}{10} \grik{-42}{10}{10}{7}{18} \gr{-45}{10}{18}{23}{18}{12}{18}{23}{18}
\grik{-88}{10}{18}{16}{18}}

\put(20,10){\gr{0}{10}{10}{14}{10}{6}{10}{15}{10} \grik{-42}{10}{10}{7}{18} \gr{-45}{10}{18}{23}{18}{12}{18}{23}{18}
\grik{-88}{10}{18}{16}{18}}

\put(20,-14){\gr{0}{10}{15}{20}{16}{12}{15}{18}{14}}
\put(85, -12){\qbezier(0,2)(5,-2)(11,0)}
\put(87,-12){\qbezier(10,0)(18,-0.1)(19,-25)}

\put(70,10){\qbezier(0,70)(3,50)(10,45)}
\put(75,25){\graph{30}{20}{27}{33}{30}{25}{27}{30}{27}}
\put(159,52){\qbezier(0,0)(15,-5)(25,0)}

\put(110,30){\qbezier(0,50)(1,35)(14,30)}
\put(170,83){\gr{50}{50}{27}{23}{30}{35}{27}{20}{27}}

\put(25,-14){\qbezier(0,0)(7,4)(14,0)}
\put(35,-14){\qbezier(5,0)(7,-7)(20,-3)}

\put(45,-17){\qbezier(10,0)(18,-0.1)(20,-20)}

\put(27,80){\qbezier(0,0)(2,-35)(10,-37)}
\put(33,40){\gr{0}{10}{13}{10}{13}{17}{10}{8}{13}}
\put(120,40){\gr{0}{10}{13}{10}{13}{17}{10}{8}{13}}
\put(100,43){\qbezier(-2,0)(10,4)(23,0)}

\put(250,-15){\thicklines \gr{0}{10}{10}{15}{10}{5}{10}{15}{10} \grik{-42}{10}{10}{5}{10} \gr{-45}{10}{10}{15}{10}{5}{10}{15}{10}
\grik{-88}{10}{10}{5}{10}}

\put(250,-7){\gr{0}{10}{30}{35}{30}{26}{20}{15}{20} \grik{-42}{10}{20}{25}{20} \gr{-45}{10}{20}{15}{10}{5}{10}{15}{10}
\grik{-88}{10}{10}{5}{10}}

\put(250,-9){\gr{0}{10}{25}{15}{10}{6}{10}{15}{10} \grik{-42}{10}{10}{5}{10} \gr{-45}{10}{10}{15}{10}{5}{10}{15}{10}
\grik{-88}{10}{10}{5}{10}}

\put(250,-12){\gr{0}{10}{20}{15}{10}{6}{10}{15}{10} \grik{-42}{10}{10}{5}{10} \gr{-45}{10}{10}{15}{10}{5}{10}{15}{10}
\grik{-88}{10}{10}{5}{10}}

\put(250,24){\gr{0}{10}{10}{14}{10}{6}{10}{15}{10} \grik{-42}{10}{10}{7}{1} \gr{-45}{10}{1}{-2}{-14}{-22}{-18}{-14}{-18}
\grik{-88}{10}{-18}{-20}{-18}}

\put(250,26){\gr{0}{10}{14}{14}{10}{6}{10}{15}{10} \grik{-42}{10}{10}{7}{1} \gr{-45}{10}{1}{-2}{-14}{-22}{-18}{-14}{-18}
\grik{-88}{10}{-18}{-20}{-18}}

\put(250,28){\gr{0}{10}{18}{14}{10}{6}{10}{15}{10} \grik{-42}{10}{10}{7}{1} \gr{-45}{10}{1}{-2}{-14}{-22}{-18}{-14}{-18}
\grik{-88}{10}{-18}{-20}{-18}}

\put(250,30){\gr{0}{10}{20}{14}{10}{6}{10}{15}{10} \grik{-42}{10}{10}{7}{1} \gr{-45}{10}{1}{-2}{-14}{-22}{-18}{-14}{-18}
\grik{-88}{10}{-18}{-20}{-18}}

\put(250,32){\gr{0}{10}{25}{14}{10}{6}{10}{15}{10} \grik{-42}{10}{10}{7}{1} \gr{-45}{10}{1}{-2}{-14}{-22}{-18}{-14}{-18}
\grik{-88}{10}{-18}{-20}{-18}}

\put(250,-24){\gr{0}{10}{15}{20}{16}{12}{15}{18}{14}}
\put(315, -22){\qbezier(0,2)(5,-2)(11,0)}
\put(317,-22){\qbezier(10,0)(18,-0.1)(19,-25)}

\put(300,10){\qbezier(0,70)(3,50)(10,35)}
\put(305,25){\graph{21}{10}{7}{1}{-1}{-5}{-14}{-17}{-15}}
\put(389,52){\qbezier(0,-42)(15,-45)(25,-43)}

\put(340,30){\qbezier(0,50)(1,35)(14,20)}
\put(400,83){\gr{50}{50}{17}{0}{-3}{-10}{-17}{-23}{-22}}

\put(255,-24){\qbezier(0,0)(7,4)(14,0)}
\put(265,-24){\qbezier(5,0)(7,-7)(20,-3)}
\put(275,-27){\qbezier(10,0)(18,-0.1)(20,-20)}

\put(257,83){\qbezier(0,0)(2,-35)(10,-37)}
\put(263,43){\gr{0}{10}{13}{0}{3}{7}{4}{2}{-3}}
\put(351,23){\gr{0}{10}{2}{-12}{-5}{-4}{-6}{-8}{-6}}
\put(330,26){\qbezier(-3,4)(10,0)(24,-11)}
\end{picture}

\vskip 70pt

\begin{picture}(40,100)
\put(-10,-40){\vector(0,1){120}}
\put(-20,15){\vector(1,0){220}}
\put(230,-40){\vector(0,1){120}}
\put(220,20){\vector(1,0){220}}

\put(188,8){$_t$}
\put(-5,75){$_{x(t)}$}
\put(15,8){$_{t_0}$}
\put(22,15){$_\circ$}
\put(158,-9){$_{x_*(t)}$}
\put(162,30){$_{x_N^+(t)}$}
\put(-10,-55){$_{pict.5.}$}
\put(20,-55){$_{I_{a,b}^+(t_0)<+\infty}$ $_{and}$ $_{I_{c,-b}^+(t_0)<+\infty}$}
\put(230,-55){$_{pict.6.}$}
\put(260,-55){$_{I_{a,b}^+(t_0) =I_{c,-b}^+(t_0)=+\infty}$}

\put(428,13){$_t$}
\put(245,13){$_{t_0}$}
\put(235,75){$_{x(t)}$}
\put(252,20){$_\circ$}

\put(0,-20){\thicklines\qbezier[30](24,-20)(24,40)(24,100)}
\put(44,-20){\thicklines\qbezier[30](24,-20)(24,40)(24,100)}
\put(84,-20){\thicklines\qbezier[30](24,-20)(24,40)(24,100)}

\put(230,-20){\thicklines\qbezier[30](24,-20)(24,40)(24,100)}
\put(274,-20){\thicklines\qbezier[30](24,-20)(24,40)(24,100)}
\put(314,-20){\thicklines\qbezier[30](24,-20)(24,40)(24,100)}
\put(294,-20){\thicklines\qbezier[30](24,-20)(24,40)(24,100)}
\put(354,-20){\thicklines\qbezier[30](24,-20)(24,40)(24,100)}

\put(20,-2){\thicklines \gr{0}{10}{10}{14}{10}{6}{10}{15}{10} \grik{-42}{10}{10}{5}{10} \gr{-45}{10}{10}{15}{10}{5}{10}{15}{10}
\grik{-88}{10}{10}{5}{10}}

\put(20,25){\thicklines \gr{0}{10}{10}{14}{10}{6}{10}{15}{10} \grik{-42}{10}{10}{5}{10} \gr{-45}{10}{10}{15}{10}{5}{10}{15}{10}
\grik{-88}{10}{10}{5}{10}}

\put(20,0){\gr{0}{10}{10}{14}{10}{6}{10}{15}{10} \grik{-42}{10}{10}{5}{10} \gr{-45}{10}{10}{15}{10}{5}{10}{15}{10}
\grik{-88}{10}{10}{5}{10}}

\put(20,3){\gr{0}{10}{12}{14}{10}{6}{10}{15}{10} \grik{-42}{10}{10}{5}{10} \gr{-45}{10}{10}{15}{10}{5}{10}{15}{10}
\grik{-88}{10}{10}{5}{10}}

\put(20,5){\gr{0}{10}{15}{14}{10}{6}{10}{15}{10} \grik{-42}{10}{10}{5}{10} \gr{-45}{10}{10}{15}{10}{5}{10}{15}{10}
\grik{-88}{10}{10}{5}{10}}

\put(20,28){\gr{0}{10}{10}{14}{10}{6}{10}{15}{10} \grik{-42}{10}{10}{7}{18} \gr{-45}{10}{18}{23}{18}{12}{18}{23}{18}
\grik{-88}{10}{18}{16}{18}}

\put(20,30){\gr{0}{10}{10}{14}{10}{6}{10}{15}{10} \grik{-42}{10}{10}{7}{18} \gr{-45}{10}{18}{23}{18}{12}{18}{23}{18}
\grik{-88}{10}{18}{16}{18}}

\put(20,18){\gr{0}{10}{10}{14}{10}{6}{8}{10}{2} \grik{-42}{10}{2}{-4}{-1} \gr{-45}{10}{-1}{2}{-1}{-8}{-1}{2}{-1}
\grik{-88}{10}{-1}{-4}{-1}}

\put(20,20){\gr{0}{10}{10}{14}{10}{6}{8}{10}{2} \grik{-42}{10}{2}{-4}{-1} \gr{-45}{10}{-1}{2}{-1}{-8}{-1}{2}{-1}
\grik{-88}{10}{-1}{-4}{-1}}

\put(20,22){\gr{0}{10}{10}{14}{10}{6}{8}{10}{2} \grik{-42}{10}{2}{-4}{-1} \gr{-45}{10}{-1}{2}{-1}{-8}{-1}{2}{-1}
\grik{-88}{10}{-1}{-4}{-1}}

\put(20,-14){\gr{0}{10}{15}{20}{16}{12}{15}{18}{14}}
\put(85, -12){\qbezier(0,2)(5,-2)(11,0)}
\put(87,-12){\qbezier(10,0)(18,-0.1)(19,-25)}

\put(70,10){\qbezier(0,70)(3,50)(10,45)}
\put(75,25){\graph{30}{20}{27}{33}{30}{25}{27}{30}{27}}
\put(159,52){\qbezier(0,0)(15,-5)(25,0)}

\put(110,30){\qbezier(0,50)(1,35)(14,30)}
\put(170,83){\gr{50}{50}{27}{23}{30}{35}{27}{20}{27}}

\put(25,-14){\qbezier(0,0)(7,4)(14,0)}
\put(35,-14){\qbezier(5,0)(7,-7)(20,-3)}

\put(45,-17){\qbezier(10,0)(18,-0.1)(20,-20)}

\put(27,80){\qbezier(0,0)(2,-35)(10,-37)}
\put(33,40){\gr{0}{10}{13}{10}{13}{17}{10}{8}{13}}
\put(120,40){\gr{0}{10}{13}{10}{13}{17}{10}{8}{13}}
\put(100,43){\qbezier(-2,0)(10,4)(23,0)}

\put(255,10){\qbezier(0,70)(2,40)(10,30)}
\put(255,10){\qbezier(10,30)(15,20)(20,30)}
\put(255,10){\qbezier(20,30)(25,40)(30,30)}
\put(255,10){\qbezier(30,30)(37,0)(40,-50)}

\put(260,30){\qbezier(0,50)(2,40)(10,30)}
\put(260,30){\qbezier(10,30)(15,20)(20,30)}
\put(260,30){\qbezier(20,30)(25,40)(30,30)}
\put(260,30){\qbezier(30,30)(37,0)(40,-70)}

\put(252,0){\qbezier(2,50)(9,30)(10,30)}
\put(252,0){\qbezier(10,30)(15,23)(20,30)}
\put(252,0){\qbezier(20,30)(25,40)(30,30)}
\put(252,0){\qbezier(30,30)(37,0)(40,-40)}

\put(252,-10){\qbezier(2,50)(9,30)(10,30)}
\put(252,-10){\qbezier(10,30)(15,23)(20,30)}
\put(252,-10){\qbezier(20,30)(25,40)(27,30)}
\put(252,-10){\qbezier(27,30)(33,0)(35,-30)}

\put(300,10){\qbezier(0,70)(2,40)(10,30)}
\put(300,10){\qbezier(10,30)(22,20)(25,10)}
\put(300,10){\qbezier(25,10)(30,0)(35,-50)}

\put(305,10){\qbezier(0,70)(2,40)(10,30)}
\put(305,10){\qbezier(10,30)(22,20)(25,10)}
\put(305,10){\qbezier(25,10)(30,0)(35,-50)}

\put(330,10){\qbezier(0,70)(2,40)(10,30)}
\put(330,10){\qbezier(10,30)(22,20)(25,10)}
\put(330,10){\qbezier(25,10)(30,0)(35,-50)}

\put(340,10){\qbezier(0,70)(2,40)(10,30)}
\put(340,10){\qbezier(10,30)(22,20)(25,10)}
\put(340,10){\qbezier(25,10)(30,0)(35,-50)}

\put(320,10){\qbezier(0,70)(2,40)(10,30)}
\put(320,10){\qbezier(10,30)(22,20)(25,10)}
\put(320,10){\qbezier(25,10)(30,0)(35,-50)}

\put(295,10){\qbezier(0,70)(2,40)(10,30)}
\put(295,10){\qbezier(10,30)(22,20)(25,10)}
\put(295,10){\qbezier(25,10)(30,0)(35,-50)}

\end{picture}

\vskip 100pt

\noindent
$\beta) \ph I_{c,-b}^+(t_0) = +\infty$; and Eq. (2.1) has a $t_1$-regular solution for some  $t_1 \ge t_0$ (see pict. 4);

\noindent
$\gamma) \ph I_{a,b}^+(t_0) < +\infty, \ph I_{c,-b}^+(t_0)< +\infty$; and Eq. (2.1) has a $t_1$-regular solution for some  $t_1 \ge t_0$ (see pict. 5);

\noindent
$\delta) \ph I_{a,b}^+(t_0) = I_{c,-b}^+(t_0) = +\infty$ (see pict. 6).

Let $a(t) > 0, \ph t\ge t_0$, and let $x_0(t)$ be a solution of Eq. (2.1) with  $x_0(t_0) = 0$.  Then by virtue of Theorem 2.3. I$^\circ \ph x_0(t)$ is $t_0$-normal and non negative. Obviously
$$
x_0(t) + \il{t_0}{t} a(\tau)\biggl(x_0(\tau)  + \frac{b(\tau)}{2 a(\tau)}\biggr)^2 d\tau = \il{t_0}{t}\frac{b^2(\tau) - 4 a(\tau) c(\tau)}{4 a(\tau)} d\tau, \phh t\ge t_0.
$$
Then since $x_0(t) \ge 0, \ph t\ge t_0$, we have
$$
\il{t_0}{t} a(\tau)\biggl(x_0(\tau)  + \frac{b(\tau)}{2 a(\tau)}\biggr)^2 d\tau \le \il{t_0}{t}\frac{b^2(\tau) - 4 a(\tau) c(\tau)}{4 a(\tau)} d\tau, \phh t\ge t_0. \eqno (2.53)
$$
According to Cauchy - Schwarz   inequality  we have:
$$
\il{t_0}{t} a(\tau)\biggl(x_0(\tau)  + \frac{b(\tau)}{2 a(\tau)}\biggr) d\tau \le \sqrt{\il{t_0}{t} a(\tau) d\tau}\sqrt{\il{t_0}{t} a(\tau)\biggl(x_0(\tau)  + \frac{b(\tau)}{2 a(\tau)}\biggr)^2 d\tau }
$$
From here and from (2.53) we get:
$$
\il{t_0}{t} a(\tau) x_0(\tau) d\tau \le - \frac{1}{2}\il{t_0}{t} b(\tau) d\tau + \frac{1}{2}\sqrt{\il{t_0}{t} a(\tau)d\tau \biggl[\il{t_0}{t}\frac{b^2(\tau) - 4 a(\tau) c(\tau)}{a(\tau)}\biggr]}d\tau,  \eqno (2.54)
$$
$t\ge t_0$.
Let  $x_*(t)$  be a  $t_0$-extremal solution to Eq. (2.1). Then (see  [11])
$$
a(t) x_*(t) = \frac{\nu_{x_0}(t)}{\nu_{x_0}(t)} - a(t) x_0(t) - b(t), \phh t\ge t_0.
$$
From here and from (2.54) we obtain:
$$
\il{t_0}{t} a(\tau) x_*(\tau) d\tau \ge -\frac{1}{2}\il{t_0}{t}b(\tau) d\tau -\phantom{aaaaaaaaaaaaaaaaaaaaaaaaaaaaaaaaaaaaaaaaaaaaaaaaa}
$$
$$
\phantom{aaaaaaa} -\frac{1}{2} \sqrt{\il{t_0}{t} a(\tau)d\tau \biggl[\il{t_0}{t}\frac{b^2(\tau) - 4 a(\tau) c(\tau)}{a(\tau)}\biggr]}d\tau + \ln\frac{\nu_{x_0}(t)}{\nu_{x_0}(t_0)}, \phh t\ge t_0.  \eqno (2.55)
$$

{\bf Remark 2.2} {\it The estimates (2.54) and  (2.55) are sharp in the sense that for   $a(t) = const, \ph b(t) = const, \ph  c(t) = const$   the estimate  (2.54) becomes an equality up to constant summand and the inequality  (2.55)  becomes an equality.

Suppose  $I_{a,b}^+(t_0) < +\infty, \ph I_{-c,-b}^+(t_0) < +\infty$. Then due to Theorem 2.3. VI$^\circ$  Eq.  (2.1)  has a negative $t_0$-normal solution. Therefore,
$$
\nu_{x_0}(t) = \ilp{t} a(\tau) \exp\biggl\{- \il {t}{\tau}\bigl[2 a(s)\bigl(x_0(\xi) - x_N^-(\xi)\bigr) + 2 a(\xi) x_N^-(\xi) + b(\xi) \bigr]d\xi\biggr\} d\tau \ge \phantom{aaaaaaaaa}
$$
$$
\phantom{aaaaaaaaaaaaaaaaaaaa} \ge \exp\biggl\{- \ilp{t} 2 a(s)\bigl(x_0(\xi)- x_N^-(\xi)\bigr)d\xi\biggr\} I_{a,b}^+(t,+\infty), \phh t\ge t_0.
$$
From here and from (2.55) we get:
$$
\il{t_0}{t} a(\tau) x_*(\tau) d\tau \ge -\frac{1}{2}\il{t_0}{t}b(\tau) d\tau -\phantom{aaaaaaaaaaaaaaaaaaaaaaaaaaaaaaaaaaaaaaaaaaaaaaaaa}
$$
$$
\phantom{aaa} -\frac{1}{2} \sqrt{\il{t_0}{t} a(\tau)d\tau \biggl[\il{t_0}{t}\frac{b^2(\tau) - 4 a(\tau) c(\tau)}{a(\tau)}\biggr]}d\tau + \ln I_{a,b}^+(t) + c, \phh t\ge t_0.  \eqno (2.56)
$$
where
$$
c\equiv - \ilpp a(\tau)\bigl(x_0(\tau) - x_N^-(\tau)\bigr) d\tau - \ln \nu_{x_0}(t_0). \eqno (2.57)
$$

Suppose  $a(t) > 0, \ph c(t) \ge 0, \ph t\ge t_0, \ph I_{a,b}^+(t_0) = +\infty$,  and  Eq. (2.1) has a    $t_0$-regular solution. Obviously
$$
x_1(t) + \il{t_0}{t} a(\tau)\biggl(x_1(\tau) + \frac{b(\tau)}{2a(\tau)}\biggr)^2 d\tau = x_1(t_0) + \il{t_0}{t} \frac{b^2(\tau) - 4 a(\tau) c(\tau)}{4 a(\tau)} d\tau, \phh t\ge t_0.
$$
Then since by Theorem  2.4. I$^* \ph x_1(t) > 0, \ph t\ge t_0$,  we have
$$
\il{t_0}{t} a(\tau)\biggl(x_1(\tau) + \frac{b(\tau)}{2a(\tau)}\biggr)^2 d\tau \le  x_1(t_0) + \il{t_0}{t} \frac{b^2(\tau) - 4 a(\tau) c(\tau)}{4 a(\tau)} d\tau, \phh t\ge t_0.
$$
From here using Cauchy - Schwarz inequality by analogy of (2.54) we get:

$$
\il{t_0}{t} a(\tau) x_0(\tau) d\tau \le - \frac{1}{2}\il{t_0}{t} b(\tau) d\tau + \phantom{aaaaaaaaaaaaaaaaaaaaaaaaaaaaaaaaaaaaaaaaaaaaaaaaa}
$$
$$
\phantom{aaaaaa}+\frac{1}{2}\sqrt{\il{t_0}{t} a(\tau)d\tau \biggl[4x_1(t_0) + \il{t_0}{t}\frac{b^2(\tau) - 4 a(\tau) c(\tau)}{a(\tau)}\biggr]}d\tau, \phh t\ge t_0.  \eqno (2.58)
$$

\vskip 20 pt

\centerline{\bf \S3. The behavior of solutions of the system  (1.1)}

\vskip 20 pt

{\bf 3.1. The general case}.  Denote  $B(t) \equiv a_{11}(t) - a_{22}(t), \ph t\ge t_0$. Consider the Riccati equation
$$
z' + a_{12}(t) z^2 + B(t) z - a_{21}(t) = 0, \phh t\ge t_0. \eqno (3.1)
$$
It is not difficult to show that the solutions  $z(t)$  of this equation existing on some interval е $[t_1;t_2) \ph (t_0 \le t_1 \le t_2 \le +\infty)$ are connected wit solutions  $(\phi(t), \psi(t))$  of the system (1.1) by correlations  (see [6])
$$
\phi(t) = \phi(t_1) \exp\biggl\{\il{t_1}{t}[a_{12}(\tau)z(\tau) + a_{11}(\tau)]d\tau\biggr\}, \ph \phi(t_1) \ne 0. \ph \psi(t) = z(t)\phi(t),  \eqno (3.2)
$$
$t\in [t_1;t_2).$
Denote by  $z_0(t)$  the solution of Eq. (3.1)  with  $z_0(t_0) = i$. This solution exists on  $[t_0;+\infty)$  (i. e. it is a $t_0$-regular solution [see [6]]). Denote: $x_0(t) \equiv Re \hskip 3pt z_0(t), \linebreak y_0(t) \equiv Im \hskip 3pt z_0(t), \ph t\ge t_0$.  Let  $(\phi_\pm(t), \psi_\pm(t))$             be the solutions of the system (1.1) with  $\phi_+(t_0) = 1, \ph \psi_+(t_0) = 0, \ph \phi_-(t_0) = 0, \ph \psi_-(t_0) = 1$. Then the following equalities are valid  (see  [6]):
$$
\phi_+(t) = \frac{J_{S/2}(t)}{\sqrt{y_0(t)}} \cos \biggl(\il {t_0}{t} a_{12}(\tau) y_0(\tau) d \tau\biggr);
$$
$$
\psi_+(t) = \frac{J_{S/2}(t)}{\sqrt{y_0(t)}} \biggr[x_0(t) \cos \biggl(\il {t_0}{t} a_{12}(\tau) y_0(\tau) d \tau\biggr) - y_0(t)\sin \biggl(\il {t_0}{t} a_{12}(\tau) y_0(\tau) d \tau  \biggr)\biggr];
$$
$$
\phi_-(t) = \frac{J_{S/2}(t)}{\sqrt{y_0(t)}} \sin \biggl(\il {t_0}{t} a_{12}(\tau) y_0(\tau) d \tau\biggr);
$$
$$
\psi_-(t) = \frac{J_{S/2}(t)}{\sqrt{y_0(t)}} \biggr[x_0(t) \sin \biggl(\il {t_0}{t} a_{12}(\tau) y_0(\tau) d \tau\biggr) +  y_0(t)\cos \biggl(\il {t_0}{t} a_{12}(\tau) y_0(\tau) d \tau  \biggr)\biggr],
$$
where  $S\equiv a_{11}(t) + a_{22}(t), \ph t\ge t_0$. From here it is not difficult to obtain the relations
$$
\phi_+(t) \cos \biggl(\il {t_0}{t} a_{12}(\tau) y_0(\tau) d \tau\biggr) + \phi_-(t)  \sin \biggl(\il {t_0}{t} a_{12}(\tau) y_0(\tau) d \tau\biggr) = \frac{J_{S/2}(t)}{\sqrt{y_0(t)}}, \eqno (3.3)
$$
$$
\psi_-(t) \cos \biggl(\il {t_0}{t} a_{12}(\tau) y_0(\tau) d \tau\biggr) - \psi_+(t)  \sin \biggl(\il {t_0}{t} a_{12}(\tau) y_0(\tau) d \tau\biggr) = J_{S/2}(t)\sqrt{y_0(t)}, \eqno (3.4)
$$
$$
\psi_-(t) \sin \biggl(\il {t_0}{t} a_{12}(\tau) y_0(\tau) d \tau\biggr) + \psi_+(t)  \cos \biggl(\il {t_0}{t} a_{12}(\tau) y_0(\tau) d \tau\biggr) = J_{S/2}(t)\frac{x_0(t)}{\sqrt{y_0(t)}}, \eqno (3.5)
$$
The real valued solutions  $(\phi(t), \psi(t))$ of the system (1.1) are connected with  $x_0(t)$  and  $y_0(t)$   by relations  (see  [6])
$$
\left\{
\begin{array}{l}
\phi(t) = \mu \frac{J_{S/2}(t)}{\sqrt{y_0(t)}}\sin \biggl(\il {t_0}{t} a_{12}(\tau) y_0(\tau) d \tau + \nu\biggr);\\
\ph\\
\psi(t) = \mu\sqrt{x_0^2(t) + y_0^2(t)}\frac{J_{S/2}(t)}{\sqrt{y_0(t)}}\cos \biggl(\il {t_0}{t} a_{12}(\tau) y_0(\tau) d \tau + \nu - \alpha_0(t)\biggr),
\end{array}
\right.
\eqno (3.6)
$$
$t\ge t_0$,
where $\mu$   and  $\nu$  are arbitrary constants;
$$
\alpha_0(t) \equiv \arcsin \frac{x_0(t)}{\sqrt{x_0^2(t) + y_0^2(t)}} = \arctan \frac{x_0(t)}{y_0(t)}, \phh t\ge t_0. \eqno (3.7)
$$

Let $u(t)$ be a continuous function on  $[t_0;+\infty)$.

{\bf Definition 3.1.} {\it The function $u(t)$ is called oscillatory if it has arbitrary large zeroes, otherwise  $u(t)$ is called non oscillatory.}

{\bf Definition 3.2.} {\it The system (1.1) is called oscillatory (non oscillatory), if for its each non trivial solution $(\phi(t),\psi(t))$   the functions  $\phi(t)$ and  $\psi(t)$  are oscillatory  (non oscillatory).}

{\bf Remark 3.1.} {\it  Some oscillatory and non oscillatory criteria are proved in [6]  (see also~[5]).}

{\bf Definition 3.3.} {\it The system (1.1) is called weak oscillatory  (weak non oscillatory), if for its each non trivial solution $(\phi(t),\psi(t))$ at least one of the functions  $\phi(t)$  and  $\psi(t)$     is oscillatory (is non oscillatory) and there exist two solutions  $(\phi_j(t),\psi_j(t)), \linebreak j=1,2$,  such that $\phi_1(t)$  and  $\psi_1(t)$  are oscillatory  (non oscillatory), and at least one of the functions   $\phi_2(t)$   and  $\psi_2(t)$ is non oscillatory (oscillatory).}

{\bf Definition 3.4.} {\it The system (1.1) is called half oscillatory if for its each non trivial solution $(\phi(t),\psi(t))$   one of the functions   $\phi(t), \ph \psi(t)$   is oscillatory and  other is non \linebreak  oscillatory.}

{\bf Definition 3.5.} {\it The system (1.1) is called singular,
if it has two non trivial solutions  $(\phi_j(t),\psi_j(t)), \ph j=1,2$,  such that  $\phi_1(t)$  and  $\psi_1(t)$  are oscillatory, and  $\phi_2(t)$  and  $\psi_2(t)$  are non oscillatory.}

{\bf Remark 3.2.} {\it It is evident that each system  (1.1) is or else oscillatory or else non oscillatory or else weak oscillatory or else weak non oscillatory or else half oscillatory or else singular.}

Example 3.1. Consider the system
$$
\sist{\phi'= \phh\ph \cos (\lambda t) \psi;}{\psi'= - \cos(\lambda t) \phi, \ph t\ge t_0} \eqno (3.8)
$$
where  $\lambda = const > 0$ is a parameter. The general solution  $(\phi(t),\psi(t))$ to this system is given by formulas:
$$
\phi(t) = c_1\sin \biggl(\frac{1}{\lambda}\sin \lambda t + c_2\biggr), \phh \psi(t) = c_1\cos \biggl(\frac{1}{\lambda}\sin \lambda t + c_2\biggr),
$$
where  $c_1$ and  $c_2$  are arbitrary constants.  It is not difficult to verify that if:

\noindent
1) $0 < \lambda \le \frac{2}{\pi}$,  then the system (3.8) is oscillatory;

\noindent
2) $\frac{2}{\pi} < \lambda \le \frac{4}{\pi}$,  then the system (3.8) is weak oscillatory;

\noindent
3) $\lambda  > \frac{4}{\pi}$, then the system (3.8) is weak non oscillatory.

Example 3.2. Consider the system
$$
\sist{\phi' = \phh \phh \phh \phh \phh \phh  \psi;}{\psi' = \bigl(-\cos^2 t - \sin t\bigr) \phi, \phh \phh t\ge t_0.} \eqno (3.9)
$$
The general solution  $(\phi(t),\psi(t))$  of this system is given by formulas:
$$
\phi(t) = e^{\sin t}\biggl(c_1 + c_2 \il{t_0}{t}e^{-2\sin \tau}d\tau\biggr), \ph \psi(t) = e^{\sin t}\biggl[c_1 \cos t + c_2\biggl\{\cos t \il{t_0}{t}e^{-2\sin \tau}d\tau + e^{-2\sin t}\biggr\}\biggr].
$$
where  $c_1$ and  $c_2$  are arbitrary constants. Obviously  $\phi(t)$  is non oscillatory and   $\psi(t)$ is oscillatory. Hence the system (3.9) is half oscillatory.

Example 3.3. Consider the system
$$
\sist{\phi' = 3 \cos t\hskip 3pt \phi - 2 \cos t\hskip 3pt   \psi;}{\psi' = 4 \cos t\hskip 3pt   \phi - 3 \cos t\hskip 3pt  \psi, \ph t\ge t_0.} \eqno (3.10)
$$
It has the solutions
$$
\left(e^{\sin t},\ph e^{\sin t}\right), \phh \left(e^{\sin t} - e^{-\sin t}, \ph e^{\sin t} - 2e^{-\sin t}\right), \phh t\ge t_0.
$$
Obviously the components of the firs solution are non oscillatory; the firs component of the second solution vanishes in the points  $\pi k \ge t_0, \ph k=0, \pm 1, \pm 2, ... $,  and the nulls of the second component of the second solution are all solutions of the equation  $\sin t = \ln \sqrt{2}$   on  $[t_0;+\infty)$. Therefore the system (3.10) is singular.

{\bf Definition 3.6} {\it The system (1.1) is called   Lyapunov stable (asymptotically stable), if its all solutions are bounded on $[t_0;+\infty)$   (vanish on  $+\infty$).}

{\bf Theorem 3.1.} {\it Let for each solution $(\phi(t), \psi(t))$ of the system (1.1) the function   $J_{-S/2}(t) \phi(t)$ is bounded.  Then there exists a solution  $(\phi_0(t), \psi_0(t))$ of the system  (1.1) such that  $J_{-S/2}(t)\psi_0(t) \not \to 0$  for  $t \to +\infty$. Moreover if in addition  $a_{12}(t)$  does not change sign and $\ilpp|a_{12}(\tau)|d\tau = +\infty$, then the system (1.1) is oscillatory and for each nontrivial solution  $(\phi(t), \psi(t))$ of the system  (1.1)  $J_{-S/2}(t)\psi(t) \not \to 0$  for  $t \to +\infty$.}

Proof. By (3.3) from the conditions of the theorem it follows that
$$
y_0(t) \ge \varepsilon, \phh t\ge t_0, \eqno (3.11)
$$
for some  $\varepsilon > 0$.  Suppose for every solution $(\phi(t), \psi(t))$ of the system
(1.1)  $J_{-S/2}(t)\psi(t) \to 0$  for  $t \to +\infty$. Then according to  (3.4) we have $y_0(t) \to 0$  for $t \to +\infty$,  which contradicts (3.11). The obtained contradiction shows the existence of a solution   $(\phi_0(t), \psi_0(t))$  of the system (2.1) with  $J_{-S/2}(t)\psi_0(t) \not \to 0$    for    $t \to +\infty$.  If in addition  $a_{12}(t)$  does not change sign and  $\ilpp|a_{12}(\tau)|d\tau = +\infty$, then from (3.11) it follows that  $|\ilpp a_{12}(\tau) y_0(\tau)d\tau| = +\infty$. From here and from (3.6) it follows oscillation of the system (1.1). From the last equality from (3.6) and  (3.11) it follows that for every solution  $(\phi(t), \psi(t))$  for the system (1.1)  the relation   $J_{-S/2}(t)\psi(t) \not \to 0$     for  $t \to +\infty$ is fulfilled.  The theorem is proved.

{\bf Theorem 3.2}. {\it Let for each solution  $(\phi(t), \psi(t))$ of the system  (1.1) the relation $J_{-S/2}(t)\phi(t) \to 0$  for  $t \to +\infty$ be satisfied. Then there exists a solution  $(\phi_0(t), \psi_0(t))$ of the system (1.1) such that $J_{-S/2}(t)\psi_0(t)$ is unbounded. Moreover if in addition  $a_{12}(t)$ does not change sign and  $\ilpp|a_{12}(\tau)|d\tau = +\infty$, then the system (1.1) is oscillatory, and for any nontrivial solution $(\phi(t), \psi(t))$ of the system (1.1) the function  $J_{-S/2}(t)\psi(t)$ is unbounded.}

Proof. By (3.3) from the condition of the theorem it follows that
$$
y_0(t) \to +\infty \mb{for} t\to +\infty. \eqno (3.12)
$$
Suppose for every solution  $(\phi(t), \psi(t))$ of the system (1.1) the function  $J_{-S/2}(t)\psi(t)$       is bounded. Then from (3.4) it follows that  $y_0(t)$ is bounded, which contradicts (3.12). Hence for at last one solution  $(\phi_0(t), \psi_0(t))$ of the system (1.1)  the function  $J_{-S/2}(t)\psi_0(t)$ is unbounded.
Suppose $a_{12}(t)$ does not change sign and   $\ilpp|a_{12}(\tau)|d\tau = +\infty$. Then from (3.6) and (3.12) it follows that the system (1.1) is oscillatory and by virtue of the second of equalities  (3.6)  from (3.12) it follows that for any nontrivial solution  $(\phi(t), \psi(t))$ of the system  (1.1) the function   $J_{-S/2}(t)\psi(t)$ is unbounded. The theorem is proved.

{\bf Theorem 3.3 (abut rings)}. {\it Suppose for every solution $\Phi(t) \equiv (\phi(t), \psi(t))$ of the system (1.1) there exists $R_\Phi > 0$ such that  $||\Phi(t)|| \le R_\Phi J_{S/2}(t), \ph t\ge t_0$.           Then for every nontrivial solution $\Phi(t)$ of the system (1.1) there exists $r_\Phi$ such that
$$
||\Phi(t)|| \ge r_\Phi J_{S/2}(t), \phh t\ge t_0. \eqno (3.13)
$$
}

Proof. By (3.3)  - (3.5) from the conditions of the theorem it follows that
$$
\sqrt{y_0(t)} \le M, \phh \frac{1}{\sqrt{y_0(t)}} \le M, \phh \frac{x_0(t)}{y_0(t)} \le M, \phh t\ge t_0. \eqno (3.14)
$$
for some $M= const > 0$. Suppose for some solution  $\Phi_0(t) \equiv (\phi_0(t), \psi_0(t))$  of the system (1.1) the relation (3.13) does not fulfill.  Then there exists  infinitely large sequence  $\{t_n\}_{n-1}^{+\infty}$ such that
$$
J_{-S/2}(t_n)\phi_0(t_n)  \to 0, \ph  J_{-S/2}(t_n)\psi_0(t_n)  \to 0 \mb{for} n \to +\infty. \eqno (3.15)
$$
By (3.6) we have:
$$
J_{-S/2}(t_n)\phi_0(t_n) = \frac{\mu_0}{\sqrt{y_0(t_n)}} \sin (\gamma_n);
$$
$$
J_{-S/2}(t_n)\psi_0(t_n) = \mu_0\sqrt{1 + \frac{x_0(t_0)}{y_0(t_0)}} \sqrt{y_0(t_n)} \cos (\gamma_n - \alpha_0(t_n)),
$$
where  $\gamma_n \equiv \il{t_0}{t_n}y_0(\tau) d\tau + \nu_0, \ph n=1,2, ... ; \ph \nu_0$                                           and  $\mu_0$  are some constants. From here from  (3.14) and  (3.15) it follows:
$$
\sin (\gamma_n) \to 0, \mb{for} n \to +\infty; \eqno (3.16)
$$
$$
\cos (\gamma_n - \alpha_0(t_n)) \to 0, \mb{for} n \to +\infty; \eqno (3.17)
$$
From (3.14) and from (3.7) it follows that there exists $\delta > 0$ such that  $|\cos (\alpha_0(t_n))| >~ \delta, \linebreak n=1,2, ... $ . From here and from (3.16) it follows that  $|\cos(\gamma_n - \alpha_0(t_n))| = 1 \pm \linebreak \pm \sqrt{1 - \sin^2 \gamma_n} \cos (\alpha_0(t_n)) + \sin \gamma_n \sin \alpha_0(t_n) \ge \delta/2$
for all enough large values of  $n$, which contradicts (3.17). The obtained contradiction proves (3.13). The theorem is proved.

{\bf Remark 3.3.} {\it The geometrical meaning of Theorem 3.4 is that if for all solutions  $\Phi(t) \equiv (\phi(t), \psi(t))$ of the system (1.1) the vector functions  $J_{-S/2}(t) \Phi(t)$ are bounded then every of them lies in some ring of radiuses  $0 < r_\Phi <~R_\Phi$.}

By correlation (1.3) between Eq. (1.2), Eq. (1.5) and the system (1.4)  From Theorems 3.1 - 3.4 it follow the following

\vskip 20pt

\centerline{\bf three principles for Eq.  (1.5)}

\vskip 20pt

 {\bf  A)} {\it If all solutions of Eq. (1.5) are bounded then it is oscillatory and for its each nontrivial solution  $\phi(t)$  the relation  $\phi'(t)\not \to 0$ for  $t\to +\infty$ is fulfilled.}

{\bf B)} {\it If all solutions of Eq. (1.5) vanish on  $+\infty$, then the derivative of its every nontrivial solution is unbounded.}

{\bf C)} {\it If Eq. (1.5) is stable by Lyapunov then for its every nontrivial solution  $\phi(t)$    there exist positive numbers  $r_\phi < R_\phi$ such that   $r_\phi \le \sqrt{\phi^2(t) + \phi'(t)^2} \le R_\phi, \ph t\ge t_0$.
}
\vskip 10pt

\begin{picture}(140,100)
\put(70,-40){\vector(0,1){140}}
\put(-10,0){\vector(1,0){210}}
\put(320,-50){\vector(0,1){140}}
\put(220,00){\vector(1,0){220}}
\put(190,5){$\phi$}
\put(75,85){$\phi'$}
\put(430,5){$\phi$}
\put(325,85){$\phi'$}
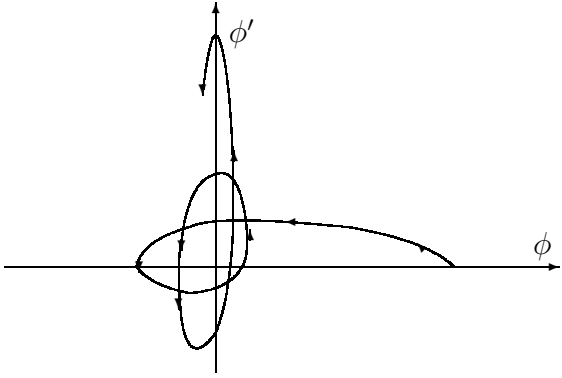
\put(-10,-60){$_{Pict.7.\ph An \ph illustration \ph to \ph the \ph  principle \ph B)}$}
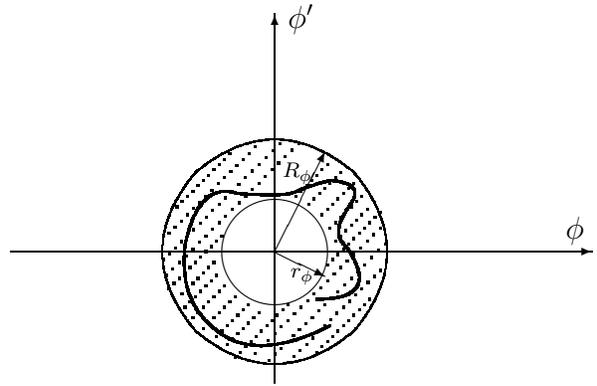
\put(230,-60){$_{Pict.8.\ph An \ph illustration \ph to \ph the \ph  principle \ph C)}$}

\put(100,0){\qbezier(60,0)(50,10)(20,15)}
\put(100,0){\qbezier(20,15)(-30,20)(-40,15)}
\put(100,0){\qbezier(-40,15)(-58,10)(-60,0)}
\put(100,0){\qbezier(-60,0)(-55,-7)(-40,-10)}
\put(100,0){\qbezier(-40,-10)(-25,-9)(-20,0)}
\put(100,0){\qbezier(-20,0)(-17,5)(-19,20)}
\put(100,0){\qbezier(-19,20)(-22,38)(-30,35)}
\put(100,0){\qbezier(-30,35)(-40,33)(-43,10)}
\put(100,0){\qbezier(-43,10)(-45,0)(-43,-20)}
\put(100,0){\qbezier(-43,-20)(-40,-39)(-30,-25)}
\put(100,0){\qbezier(-30,-25)(-20,0)(-25,65)}
\put(100,0){\qbezier(-25,65)(-30,110)(-35,65)}
\put(150,5){\vector(-1,1){4}}
\put(100,17){\vector(-1,0){4}}
\put(57,10){\vector(0,-1){4}}
\put(41,2){\vector(0,-1){4}}
\put(56,-12){\vector(0,-1){4}}
\put(83,10){\vector(0,11){4}}
\put(77,40){\vector(0,11){4}}
\put(65,69){\vector(0,-1){4}}

\put(320,0){\circle{40}}
\put(290,30){\qbezier(0,0)(30,25)(60,0)}
\put(290,30){\qbezier(60,0)(85,-30)(60,-60)}
\put(290,-30){\qbezier(0,0)(30,-25)(60,0)}
\put(290,-30){\qbezier(0,0)(-25,30)(0,60)}
\put(320,0){\vector(2,-1){19}}
\put(320,0){\vector(1,2){19}}
\put(323,30){$_{R_\phi}$}
\put(326,-8){$_{r_\phi}$}
\put(283,17){\thicklines \qbezier[10](0,0)(10,10)(20,20)}
\put(280,7){\thicklines \qbezier[15](0,0)(17.5,17.5)(35,35)}
\put(278,-3){\thicklines \qbezier[18](0,0)(20,20)(44,44)}
\put(279,-9){\thicklines \qbezier[20](0,0)(25,25)(50,50)}
\put(280,-14){\thicklines \qbezier[10](0,0)(10,10)(20,20)}
\put(305,-40){\thicklines \qbezier[10](0,0)(10,10)(20,20)}
\put(340,-5){\thicklines \qbezier[10](0,0)(10,10)(20,20)}
\put(285,-16){\thicklines \qbezier[6](0,0)(6,6)(12,12)}
\put(286,-21){\thicklines \qbezier[6](0,0)(6,6)(12,12)}
\put(289,-25){\thicklines \qbezier[6](0,0)(6,6)(12,12)}
\put(293,-29){\thicklines \qbezier[6](0,0)(6,6)(12,12)}
\put(298,-32){\thicklines \qbezier[6](0,0)(6,6)(12,12)}
\put(302,-35){\thicklines \qbezier[6](0,0)(6,6)(12,12)}

\put(328,21){\thicklines \qbezier[6](0,0)(6,6)(12,12)}
\put(333,18){\thicklines \qbezier[6](0,0)(6,6)(12,12)}
\put(338,14){\thicklines \qbezier[6](0,0)(6,6)(12,12)}
\put(341,9){\thicklines \qbezier[6](0,0)(6,6)(12,12)}
\put(343,4){\thicklines \qbezier[6](0,0)(6,6)(12,12)}

\put(336,-38){\thicklines \qbezier[10](0,0)(10,10)(20,20)}
\put(325,-42){\thicklines \qbezier[15](0,0)(17.5,17.5)(35,35)}
\put(318,-41){\thicklines \qbezier[18](0,0)(20,20)(44,44)}
\put(310,-41){\thicklines \qbezier[20](0,0)(25,25)(50,50)}

\put(336,-18){\thicklines \qbezier(0,0)(25,0)(11,20)}
\put(336,-18){\thicklines \qbezier(11,20)(5,25)(12,35)}
\put(336,-18){\thicklines \qbezier(12,35)(19,44)(5,45)}
\put(336,-18){\thicklines \qbezier(5,45)(0,44)(-10,40)}
\put(336,-18){\thicklines \qbezier(-10,40)(-15,39)(-30,40)}
\put(336,-18){\thicklines \qbezier(-30,40)(-45,45)(-50,20)}
\put(336,-18){\thicklines \qbezier(-50,20)(-52,0)(-40,-10)}
\put(336,-18){\thicklines \qbezier(-40,-10)(-25,-25)(5,-10)}

\end{picture}

\vskip 70pt

Let us compare these principles  wit the following assertion proved in [1]  (see [1], p.  222,
Corollary 6.2.4).

{\bf Proposition 3.1}. {\it  Let  $|r(t)| \le M, \ph t\ge t_0$. Then if all solutions of Eq. (1.5) vanish on $+\infty$, then Eq. (1.5) is asymptotically stable.}

Obviously from any of principles  A) -  C) it follows that the equations (1.5) satisfying the conditions of Proposition 3.1,  form a empty set.

Example 3.4.  Consider the Mathieu equation (see [14])
$$
\phi'' + (\delta +\varepsilon \cos t)\phi = 0, \phh t\ge t_0, \phh \delta, \varepsilon \in R.
$$
From the principle A) it follows that for all pairs  $(\delta,\varepsilon)$  from the  zones of stability this equation is oscillatory, and from the principle  C) it follows that (for this restriction) for its each  nontrivial solution   $\phi(t)$ there exist $R_\phi> r_\phi > 0$ such that  $r_\phi \le \phi^2(t) + \phi'(t)^2 \le R_\phi, \ph t\ge t_0$ (see Pict. 8), which  agrees quite well with the Floquet's theory. Note that some part of mentioned above zones of stability relates to the extremal case of Eq. (1.5), when $\ilp{t_0}r(t) dt = - \infty$.

Example 3.5. Consider the Airy's equation
$$
\phi'' + t \phi = 0, \phh t\ge t_0.
$$
By virtue of L. A.  Gusarov's theorem  (see [15], Theorem 1) all solutions of this equation vanish on $+\infty$.   From the principles  A)  and  B) it follows that this equation is oscillatory and for its every nontrivial solution  $\phi(t)$ the function $\phi'(t)$ is unbounded (see Pict. 7).

{\bf 3.2.  A generalization of Leighton's theorem.} In 1952 Leighton proves (see [16], p. 70, Theorem 2.24), that the equation
$$
(p(t)\phi')' + r(t)\phi =0, \phh t\ge t_0,
$$
is oscillatory if
$$
\ilpp\frac{d\tau}{p(\tau)} = \ilpp r(\tau) d \tau = +\infty.
$$
The next theorem generalizes this result.

{\bf Theorem  3.4.} {\it Let the conditions

\noindent
1) $a_{12}(t) \ge 0, \ph t \ge t_0$;

\noindent
2) $\ilpp a_{12}(\tau) \exp \biggl\{- \il{t_0}{\tau} B(s) d s\biggr\} d \tau = - \ilpp a_{21}(\tau) \exp \biggl\{\il{t_0}{\tau} B(s) d s\biggr\} d \tau = + \infty.$

\noindent
be satisfied. Then the system (1.1) is oscillatory.}

Proof. Suppose the system (1.1) is not oscillatory. Then  from the condition 1) and from the Lemma 4.2  of work [6] it follows that Eq. (3.1)
has a solution on  $[t_1;+\infty)$  for some  $t_1 \ge t_0$. Set  $W(t) \equiv - a_{21}(t) \exp\biggl\{\il{t_1}{t} B(\tau) d\tau \biggr\}, \hskip 2pt t\ge t_1$. In Eq. (3.1) make the substitution
$$
y= z \exp\biggl\{-2\il{t_1}{t} W(\tau)d\tau\biggr\}, \phh t\ge t_1.
$$
We obtain
$$
z' + U(t) z^2 + W(t) =0, \phh t\ge t_1. \eqno (3.18)
$$
where  $U(t)\equiv a_{12}(t)\exp\biggl\{- \il{t_1}{t} B(\tau) d\tau\biggr\}, \ph t\ge t_0.$                            Show that
$$
\ilp{t_1} U(\tau) \exp\biggl\{\il{t_1}{\tau} U(s) d s \il{t_1}{s}W(\zeta) d \zeta \biggr\} d \tau = +\infty. \eqno (3.19)
$$
By 2) we have:
$$
\il{t_1}{t} W(\tau) d\tau = - \il{t_1}{t}a_{21}(\tau)\exp\biggl\{\il{t_1}{\tau} B(s) d s\biggr\} d \tau \ge 0, \phh t\ge t_2,
$$
for some $t_2 \ge t_1$. By 2) from here it follows (3.19). In Eq. (3.18) make the substitution
$$
z = u - \il{t_1}{t} W(\tau) d \tau, \phh t\ge t_1.
$$
We  get:
$$
u' +U(t) u^2 - 2U(t)\il{t_1}{t}W(\tau) d\tau u +U(t)\biggl[\il{t_1}{t}W(\tau) d\tau\biggr]^2 = 0, \phh t\ge t_1. \eqno (3.20)
$$
Since Eq. (3.1) has a solution on  $[t_1;+\infty)$, from the proven above substitutions is seen that Eq. (3.20) also has a solution on  $[t_1;+\infty)$.  By virtue of Lemma 2.4 from here, from (3.19)
and from the inequalities  $a_{12}(t) \ge 0, \ph U(t)\biggl[\il{t_1}{t}W(\tau) d\tau\biggr]^2  \ge 0, \ph t\ge t_1$,  it follows that Eq. (3.20) has a positive solution on  $[t_1;+\infty)$. Then  $z_0(t) \equiv u_0(t) - \il{t_1}{t}W(\tau) d\tau$ is a solution to Eq. (3.18) such that
$$
z_0(t) > \il{t_1}{t}W(\tau) d\tau, \phh t\ge t_1. \eqno (3.21)
$$
From (3.18) it follows that
$$
z_0(t) = z_0(t_0) - \il{t_1}{t}U(\tau) z_0^2(\tau) d \tau - \il{t_1}{t}W(\tau) d\tau, \phh t\ge t_1. \eqno (3.22)
$$
From here and from (3.21) we have:
$$
0 \le \il{t_1}{t} U(\tau) z_0^2(\tau) d\tau < z_0(t_1), \phh t\ge t_1. \eqno (3.23)
$$
$(z_0(t_1) = u_0(t_1) > 0)$.  Taking into account 2) from here we  get:
$$
\biggl[z_0(t_1) - \il{t_1}{t} U(\tau)z_0^2(\tau) d\tau - \il{t_1}{t} W(\tau) d \tau\biggr]^2 \ge 1, \phh t\ge T,
$$
for some $T\ge t_1$. From here and from (3.22) it follows that  $z_0^2(t) \ge 1, \ph t\ge T.$             Hence by 2) we have $\ilp{T} U(\tau) z_0^2(\tau) d\tau \ge \ilp{T} U(\tau) d\tau = +\infty$,  which contradicts (3.23). The theorem is proved.

{\bf 3.3. The case when the system (1.1) has a regular solution}. In the sequel we will assume that the functions $a_{12}(t)$ and  $a_{21}(t)$ have unbounded supports (the case when one of them has bounded support, is trivial).

{\bf Definition 3.7.} {\it A solution  $(\phi(t), \psi(t))$ of the system (1.1) is called   $t_1$-regular if \linebreak $\phi(t) \ne 0, \ph t\ge t_1$}.

Hereafter any  $t_1$-regular solution for any  $t_1\ge t_0$ we will just call  a regular solution.

{\bf Definition 3.8}. {\it The system (1.1) is called regular if it has at least one regular solution.

{\bf Definition 3.9}. {\it The system (1.1) is called strongly regular if each its nontrivial solution is regular}.

{\bf Definition 3.10.}  {\it The system (1.1) is called exotic, if it has  up to arbitrary multiplier the unique regular solution.}

{\bf Definition 3.11.} {\it The regular system (1.1) is called normal if it is not exotic and for its any two regular solutions $(\phi_j(t), \psi_j(t)), \ph j=1,2,$     the function  $\frac{\phi_1(t)}{\phi_2(t)}$ is bounded on $[T;+\infty)$,  where  $T=T(\phi_1,\phi_2)$  such that $\phi_j(t) \ne 0, \ph t\ge T, \ph j=1,2$.}

{\bf Definition 3.12}. {\it The  system (1.1) is  called extremal if it is not exotic  and has a regular solution $(\phi_*(t),  \psi_*(t))$ such that for its any  regular solution   $(\phi(t), \psi(t))$,  linearly independent of $(\phi_*(t),  \psi_*(t))$  the equality  $\lim\limits_{t\to +\infty} \frac{\phi_*(t)}{\phi(t)} = 0$ is fulfilled.}

The unique (up to arbitrary multiplier) solution $(\phi_*(t),  \psi_*(t))$, defined above we will call the minimal solution of the system (1.1)

{\bf Remark 3.4}. {\it Every extremal system is strongly regular. Indeed Let $(\phi_*(t),  \psi_*(t))$ be the minimal solution of the system (1.1) and let $(\phi_1(t),  \psi_1(t))$ be another linearly independent regular solution of the system (1.1). Then for an arbitrary nontrivial solution  $(\phi(t),  \psi(t))$ the equality
$$
\phi(t) = \phi_*(t)\biggl[c_1 + c_2\frac{\phi_*(t)}{\phi_1(t)}\biggr], \phh t\ge t_1,
$$
is fulfilled, where $c_1$ and $c_2$ are some constants, $c_1 \ne 0, \ph t_1 \ge t_0$ such that $\phi_1(t) \ne~ 0, \ph t\ge~t_1$. Since $\lim\limits_{t\to +\infty} \phi_*(t)/\phi_1(t) = 0$ from the last equality it follows that $\phi(t) \ne 0, \ph t\ge t_2$ for some $t_2 \ge t_1$. Therefore the system (1.1) is strongly regular.

{\bf Definition 3.13}. {\it The regular system (1.1) is called super extremal if it is neither exotic, nor normal nor extremal.}

From Theorem 2.2 and Remark 2.1 it follows that the system (1.1) having a regular solution is or else exotic or else normal or else extremal or else super extremal.

On the basis of (3.2) and relations (21) and  (24) of work [11] it is easy to establish that:
to exotic system (1.1) corresponds Eq. (3.1) with the unique regular solution;
to normal system (1.1) corresponds  Eq. (3.1), having only normal solution; to extremal system (1.1) corresponds Eq. (3.1), which has along with the normal solutions  the unique extremal solution; to super extremal system (1.1) corresponds Eq. (3.1), which has along with the normal solutions  two  extremal solutions.

Let  $(\phi_0(t),\psi_0(t))$  and  $(\phi_*(t),\psi_*(t))$ be solutions of the system (1.1) corresponding by formulae (3.2) to a normal  and  a extremal solutions respectively. Then by virtue of the relations  (21) and  (24) we have:
$$
\liminf\limits_{t\to +\infty}\frac{|\phi_*(t)|}{|\phi_0(t)|} = 0, \phh \limsup\limits_{t\to +\infty}\frac{|\phi_*(t)|}{|\phi_0(t)|}< +\infty.
$$

For arbitrary continuous function    $u(t)$ on $[t_0;+\infty)$  denote:
$$
\widetilde{\mu}_u(t_1;t) \equiv \il{t_1}{t} a_{12}(\tau) \exp\biggl\{ - \il{t_1}{\tau}\Bigl[2 a_{12}(s) u(s) + B(s)\Bigr]d s\biggr\} d \tau,\phantom{aaaaaaaaaaaaaaaaaaaaaaaaaa}
$$
$$
\phantom{aaaaaaaaaaaaaaaaaa}\widetilde{\nu}_u(t) \equiv \ilp{t} a_{12}(\tau) \exp\biggl\{ - \il{t}{\tau}\Bigl[2 a_{12}(s) u(s) + B(s)\Bigr]d s\biggr\} d \tau, \phh t_1, t \ge t_0.
$$
Let the system  (1.1) be exotic and let  $x_0(t)\equiv \frac{\psi_0(t)}{\phi_0(t)}, \ph t\ge t_1$,  where $(\phi_0(t),\psi_0(t))$ is the unique (up to arbitrary multiplier)  $t_1$-regular solution of Eq.  (3.1). Then (see [11]) $\limsup\limits_{t \to +\infty}\widetilde{\mu}_{x_0}(t_1;t) = +\infty; \ph \liminf\limits_{t \to +\infty}\widetilde{\mu}_{x_0}(t_1;t) = -\infty$. By (3.2) from here it follows that
$$
\limsup\limits_{t\to +\infty}\il{t_1}{t}\frac{a_{12}(\tau)J_S(t_1;\tau)}{\phi_0^2(\tau)}d\tau = +\infty, \phh \liminf\limits_{t\to +\infty}\il{t_1}{t}\frac{a_{12}(\tau)J_S(t_1;\tau)}{\phi_0^2(\tau)}d\tau = -\infty.
$$
Let the system (1.1) be super extremal and let  $(\phi_*(t),\psi_*(t))$  and   $(\phi_{**}(t),\psi_{**}(t))$   be solutions of the system (1.1), corresponding by formulae (3.2)  to some extremal solutions of Eq. (3.1). Then by relation (24) of work [11] the equalities
$$
\liminf\limits_{t\to +\infty}\frac{|\phi_*(t)|}{|\phi_{**}(t)|} = 0, \phh \limsup\limits_{t\to +\infty}\frac{|\phi_*(t)|}{|\phi_{**}(t)|}=+\infty.
$$
is fulfilled. As we see the exotic and normal systems are the most "poor" in the view of diversity of asymptotic behavior theirs regular solutions. In the same point of view the super extremal system  is the most "rich".

{\bf Theorem 3.5}. {\it  Let for some $t_1 \ge t_0$  there exists a real valued  $t_1$-regular solution  $(\phi_0(t), \psi_0(t))$ of the system (1.1), such that the integral  $\ilp{t_1}\frac{a_{12}(\tau) J_S(\tau)}{\phi_0^2(\tau)} d\tau$ is convergent (conditionally) and   $\ilp{t}\frac{a_{12}(\tau) J_S(\tau)}{\phi_0^2(\tau)} d\tau \ne 0$  for all  $t \ge t_1$. Then:

\noindent
a) the system (1.1) is extremal and for its minimal solution $(\phi_*(t), \psi_*(t))$ the equality
$$
\ilp{t_1}\frac{a_{12}(\tau) J_S(\tau)}{\phi_*^2(\tau)} d\tau = \pm\infty, \eqno(3.18)
$$
is fulfilled;

\noindent
b) for its each  solution $(\phi(t), \psi(t))$ linearly independent of $(\phi_*(t), \psi_*(t))$ the integral \linebreak $\ilp{T}\frac{a_{12}(\tau) J_S(\tau)}{\phi^2(\tau)} d\tau$ converges (conditionally) where  $T=T(\phi;\psi) \ge t_0$  such that  $\phi(t) \ne 0, \linebreak t\ge~T$, and if in addition  $a_{12}(t)$   does not change sign and  $\ilpp |a_{12}(\tau) J_S(\tau)| d\tau = +\infty$, then
$$
\limsup\limits_{t\to+\infty} |\phi(t)| = +\infty; \eqno (3.19)
$$

\noindent
c) for arbitrary real valued solutions  $(\phi_j(t), \psi_j(t)), \ph j=1,2,$ of the system (1.1), linearly independent of  $(\phi_*(t), \psi_*(t))$, there exists finite limit
$$
\lim\limits_{t\to +\infty} \frac{\phi_1(t)}{\phi_2(t)} (\ne 0). \eqno (3.20)
$$
}

Proof. $x_0(t) \equiv \frac{\psi_0(t)}{\phi_0(t)}, \ph t\ge t_1$. By (3.2) we have:  $\widetilde{\nu}_{x_0}(t_1) = \ilp{t_1}\frac{a_{12}(\tau) J_S(\tau)}{\phi_0^2(\tau)} d\tau$.                          From here and from the conditions of the theorem it follows that the integral  $\widetilde{\nu}_{x_0}(t_1)$  converges and $\widetilde{\nu}_{x_0}(t) \ne 0, \ph t\ge t_1$.  Then by virtue of Theorem 2.2 Eq. (3.1) has the unique $t_1$-extremal solution $x_*(t)$. By (3.2) the functions
$$
\phi_*(t) \equiv \exp\biggl\{\il{t_1}{t}\Bigl[a_{12}(\tau) x_*(\tau) + a_{11}(\tau)\Bigr]d\tau\biggr\}, \phh \psi_*(t)\equiv x_*(t)\psi_*(t), \phh  t\ge t_1,
$$
define a real valued $t_1$-regular solution $(\phi_*(t), \psi_*(t))$ of the system (1.1), which can be continued on $[t_0;+\infty)$  as a solution of Eq. (1.1).  Then by  (2.7) from the equality  $\widetilde{\nu}_{x_*}(t_1) = \ilp{t_1}\frac{a_{12}(\tau) J_S(\tau)}{\phi_*^2(\tau)} d\tau$ it follows (3.18). From (3.18) and from the conditions of the theorem it follows that $(\phi_0(t), \psi_0(t))$  and $(\phi_*(t), \psi_*(t))$  are linearly independent. Therefore for arbitrary solution  $(\phi(t), \psi(t))$  of the system (1.1), linearly independent of $(\phi_*(t), \psi_*(t))$ the equality
$$
\phi(t) = \phi_*(t)\biggl[c_0 + c_*\frac{\phi_*(t)}{\phi_0(t)}\biggr], \phh t\ge t_1, \eqno (3.21)
$$
is fulfilled, where  $c_0$   and   $c_*$  are some constants and  $c_0\ne 0$. By virtue of Theorem 2.2 from the equality  $\widetilde{\nu}_{x_0}(t_1) = \ilp{t_1}\frac{a_{12}(\tau) J_S(\tau)}{\phi_0^2(\tau)} d\tau$  and from the conditions of the theorem it follows that   $x_0(t)$ is a  $t_1$-normal solution of Eq. (3.2). Then according to (2.8) from the equality
$$
\frac{\phi_*(t)}{\phi_0(t)} = \frac{\phi_*(t_1)}{\phi_0(t_1)} \exp\biggl\{\il{t_1}{t} a_{12}(\tau)\Bigl[x_0(\tau) - x_*(\tau)\bigr]d\tau\biggr\}, \phh t\ge t_1,
$$
it follows that
$$
\lim\limits_{t\to +\infty}\frac{\phi_*(t)}{\phi_0(t)} = 0. \eqno (3.22)
$$
From here and from (3.21) it follows that  $\phi(t) \ne 0, \ph t\ge t_2$, for some $t_2 \ge t_1$.
The assertion a) is proved. The assertion b)  follows from Theorem 2.2
and from the equality $\widetilde{\nu}_{x}(T) = \ilp{T}\frac{a_{12}(\tau) J_S(T;\tau)}{\phi^2(\tau)} d\tau$,                        where $x(t)\equiv \frac{\psi(t)}{\phi(t)}$ is a $T$-normal solution of Eq. (3.2), and the equality (3.19) immediately follows from the  absolutely convergence of the integral   $\ilp{T}\frac{a_{12}(\tau)J_S(T;\tau)}{\phi^2(\tau)} d\tau$ and from the equality
$$
\ilpp|a_{12}(\tau) J_S(T;\tau)|d\tau = J_{-S}(T)\ilpp|a_{12}(\tau) J_S(\tau)|d\tau.
$$
Let us prove c). Set $x_j(t) \equiv \frac{\psi_j(t)}{\phi_j(t)}, \ph t\ge T, \ph j=1,2,$ where  $T \ge t_0$    such that $\phi_j(t) \ne 0, \ph t\ge T, \ph j=1,2.$ Since  $(\phi_j(t), \psi_j(t)), \ph j=1,2,$ are linearly independent of  $(\phi_*(t), \psi_*(t))$, on the strength of (3.2) and Theorem 2.2 the functions $x_j(t)$ are $T$-normal solutions of Eq. (3.1), and
$$
\frac{\phi_1(t)}{\phi_2(t)} = \frac{\phi_1(T)}{\phi_2(T)}\exp\biggl\{\il{T}{t}a_{12}(\tau)\Bigl[x_1(\tau) - x_2(\tau)\Bigr]d\tau\biggr\}, \phh t\ge T.
$$
By (2.9) from here it follows (3.20). The theorem is proved.

{\bf Remark 3.5}. {\it  Theorem 3.5 is a generalization and supplement of the assertions  (i) and  (ii) of Theorem 6.4 from the book  [17] (see [17], p. 355 ,  Theorem 6.4)}

Denote by  $reg(t_1)$ the set of initial values  $x_{(0)}\in R$, for which the solution  $x(t)$   of Eq. (3.1) with  $x(t_1) = x_{(0)}$  is $t_1$-regular. If $a_{12}(t) \ge 0 \ph (\le 0), \ph t\ge t_0$, and Eq. (3.1) has a real valued $t_1$-regular solution then $reg(t_1) = [x_*(t); +\infty) \ph (reg(t_1) = (- \infty;x_*(t)])$, where  $x_*(t)$ is the unique $t_1$-extremal solution to Eq. (3.1)  (see [11]).

{\bf Corollary 3.1}. {\it If $a_{12}(t)$ does not change sign, then the system (1.1) is or else oscillatory or else regular. In the last case the assertions of Theorem 3.5 are valid, and if:

\noindent
A$) \ph I_{|a_{12}|,-S}^+(t_0) = +\infty$, then the system (1.1) has a solution  $(\phi_0(t), \psi_0(t))$  such that $\limsup\limits_{t \to +\infty}|\phi_0(t)| = +\infty$;

\noindent
B$) \ph I_{|a_{12}|,-S}^+(t_0) < +\infty$, then the system (1.1) has a solution  $(\phi_*(t), \psi_*(t))$  such that $\liminf\limits_{t \to +\infty}|\phi_*(t)| = 0$.
}

Proof. If  $a_{12}(t)$  does not change sign then or else  $\ilpp |a_{12}(\tau)|y_0(\tau) d\tau = +\infty$       or else $\ilpp |a_{12}(\tau)|y_0(\tau) d\tau < +\infty$. By (3.6) in the first case the system  (1.1) is oscillatory. In the second case the system (1.1), and therefore Eq. (3.1) have  $t_1$-regular solutions for some  $t_1 \ge t_0$. Then  $reg(t_1) = [x_*(t_1);+\infty)$  or  $reg(t_1) = (-\infty;
x_*(t_1)]$, where $x_*(t)$ is the unique $t_1$-extremal solution of Eq. (3.1). From here it follows that Eq. (3.1) has a $t_1$-normal solution $x_0(t)$. Then by virtue of Theorem 2.2 the integral  $\widetilde{\nu}_{x_0}(t_1)$ converges and since  $a_{12}(t)$ does not change sign and has unbounded support we have
$$
\widetilde{\nu}_{x_0}(t) \ne 0, \phh t\ge t_1. \eqno (3.23)
$$
By (3.2) we have $\widetilde{\nu}_{x_0}(t) = \ilp{t}\frac{a_{12}(\tau) J_S(t;\tau)}{\phi_0^2(\tau)} d\tau, \ph t\ge t_1$, where $\phi_0(t)\equiv\exp\biggl\{\il{t_1}{t}\Bigl[a_{12}(\tau) x_0(\tau) + a_{11}\Bigr]d\tau\biggr\}, \linebreak \psi_0(t) \equiv x_0(t) \phi_0(t), \ph t\ge t_1$, form the  $t_1$-regular solution  $(\phi_0(t), \psi_0(t))$ of the system (1.1) on $[t_1;+\infty)$, which can be continued on $[t_0;+\infty)$ as a solution of the system  (1.1).  From here and from  (3.23) it follows that all conditions of Theorem  3.5 are fulfilled. Therefore all assertions of Theorem 3.5 are valid.
Suppose  $I_{|a_{12}|,-S}^+(t_0;+\infty) = +\infty$. Show that for the solution $(\phi_0(t), \psi_0(t))$     the component  $\phi_0(t)$ is unbounded. Suppose that  $|\phi_0(t)| \le M, \ph t \ge t_0.$
Then by (3.2) we have $|\widetilde{\nu}_{x_0}(t_1)| = \ilp{t_1}\frac{|a_{12}(\tau)|J_s(t_1;\tau)}{\phi_0^2(\tau)}d\tau \ge \frac{1}{M^2} \ilp{t_1}|a_{12}(\tau)|J_s(t_1;\tau)d\tau = +\infty$. But on the other hand since  $x_0(t)$ is  $t_1$-normal we have $|\widetilde{\nu}_{x_0}(t_1)|< +\infty$.
We came to the contradiction. The assertion  A$)$  is proved. Let us prove B$)$.  Set
$$
\phi_*(t)\equiv \exp\biggl\{\il{t_1}{t}\Bigl[a_{12}(\tau)x_*(\tau) + a_{11}(\tau)\Bigr]d\tau\biggr\}, \phh \psi_*(t)\equiv x_*(t)\phi_*(t), \phh t\ge t_1,
$$
By  (3.2) $(\phi_*(t), \psi_*(t))$ is a solution of the system (1.1) on  $[t_1;+\infty)$, which can be continued on   $[t_0;+\infty)$ as a solution of the system (1.1). Show that   $\liminf\limits_{t \to +\infty}|\phi_*(t)| = 0$. Suppose that it is not true. Then  $|\phi_*(t)| \ge \varepsilon > 0, \ph t\ge T_0,$  for some $T_0 \ge t_0$. Hence by  (3.2) from the condition
$I_{|a_{12}|,-S}^+(t_0) < +\infty$ it follows that  $|\widetilde{\nu}_{x_*}(t_1)| = \ilp{t_1}\frac{|a_{12}(\tau)|J_s(t_1;\tau)}{\phi_*^2(\tau)}d\tau \le \frac{1}{\varepsilon^2} \ilp{t_1}|a_{12}(\tau)|J_s(t_1;\tau)d\tau < +\infty$.
But on the other hand since  $x_*(t)$ is $t_1$-extremal and  $a_{12}(t)$
does not change sign, we have  $\widetilde{\nu}_{x_*}(t_1) =\pm \infty$.
We came to the contradiction. Therefore  $\liminf\limits_{t \to +\infty}|\phi_*(t)| = 0$.
The corollary is proved.

{\bf Corollary 3.2} {\it Let for some  $t_1 \ge t_0$ there exists a  $t_1$-regular solution  $(\phi_*(t), \psi_*(t))$ of the system (1.1) such that $\ilp{t_1}\frac{a_{12}(\tau)J_S(\tau)}{\phi_*^2(\tau)}d\tau = \pm\infty$.  Then the assertions of Theorem 3.5 are valid.}

Proof. Set  $x_*(t)\equiv\frac{\psi_*(t)}{\phi_*(t)}, \ph t\ge t_1$. Obviously  $x_*(t)$ is a      $t_1$-regular solution of Eq. (3.1).  On the basis of (3.2) it is easy to show that  $\widetilde{\nu}_{x_*}(t_1) = \ilp{t_1}\frac{a_{12}(\tau)J_S(\tau)}{\phi_*^2(\tau)}d\tau$. From here and from the condition of corollary it follows that  $\widetilde{\nu}_{x_*}(t_1) = \pm\infty$. On the basis of Corollary  2.1 from this relation we obtain  that Eq. (3.1) ha the unique  $t_1$-extremal solution (which coincides with  $x_*(t)$).   Then by virtue of Theorem 2.2 Eq. (3.1) has a  $t_1$-normal solution  $x_0(t)$  such that the integral  $\widetilde{\nu}_{x_0}(t)$  converges for all  $t\ge t_1$  and  $\widetilde{\nu}_{x_0}(t) \ne 0, \ph t\ge t_1$.  Set
$$
\phi_0(t)\equiv \exp\biggl\{\il{t_1}{t}\bigl[a_{12}(\tau) x_0(\tau) + a_{11}(\tau)\bigr]d\tau\biggr\}, \phh \psi_0(t)\equiv x_0(t)\phi_0(t), \phh t\ge t_1.
$$
By (3.2) $(\phi_0(t), \psi_0(t))$ is a solution of the system (1.1) on  $[t_1;+\infty)$, which can be continued on  $[t_0;+\infty)$  as a solution of the system (1.1) By (3.2) we have $\ilp{t}\frac{a_{12}(\tau)J_S(\tau)}{\phi_0^2(\tau)}d\tau = \widetilde{\nu}_{x_0}(t) \ne 0, \ph t\ge t_1$.  Therefore for the solution $(\phi_0(t), \psi_0(t))$ of the system (1.1) all conditions of Theorem 3.5 are fulfilled. The corollary is proved.

{\bf Theorem 3.6}. {\it Let  $a_j(t)$  and  $c_j(t)$ be continuously differentiable functions on  $[t_0;+\infty)$  such that  $(-1)^j a_j(t) > 0, \ph  (-1)^j c_j(t) > 0, \ph t\ge t_0, \ph j=1,3$, and \hskip 2pt let  $a_1(t) \le a_{12}(t) \le~a_2(t), \linebreak c_1(t) \le - a_{21}(t) \le c_2(t), \ph B(t) \ge \frac{1}{2} \biggl(\frac{a_j'(t)}{a_j(t)} - \frac{c_j(t)}{c_j(t)}\biggr) + 2(-1)^j\sqrt{a_j(t)c_j(t)}, \ph t\ge t_0, \ph j=1,2.$
Then every solution  $(\phi(t), \psi(t))$ of the system (1.1) with  $\frac{\psi(t_0)}{\phi(t_0)} \in \Bigl[-\sqrt{\frac{c_2(t)}{a_2(t)}}; \sqrt{\frac{c_1(t)}{a_1(t)}}\Bigr]$  is \linebreak$t_0$-regular and
$$
-\sqrt{\frac{c_2(t)}{a_2(t)}} \le \frac{\psi(t)}{\phi(t)} \le \sqrt{\frac{c_1(t)}{a_1(t)}}, \phh t\ge t_0. \eqno (3.24)
$$
If $\frac{\psi(t_0)}{\phi(t_0)} \in \Bigl(-\sqrt{\frac{c_2(t)}{a_2(t)}}; \sqrt{\frac{c_1(t)}{a_1(t)}}\Bigr)$, then the function  $I_\phi(t) \equiv \il{t_0}{t}\frac{a_{12}(\tau) J_S(t_0;\tau)}{\phi_0^2(\tau)} d\tau, \ph t\ge t_0$, is bounded.
}

Proof. By virtue of Theorem 4.2 of work  [11] if the conditions of the theorem are fulfilled then every solution  $x(t)$  of Eq. (3.1) with  $x(t_0) \in \Bigl[-\sqrt{\frac{c_2(t)}{a_2(t)}}; \sqrt{\frac{c_1(t)}{a_1(t)}}\Bigr]$ is  $t_0$-regular and $-\sqrt{\frac{c_2(t)}{a_2(t)}} \le  x(t) \le \sqrt{\frac{c_1(t)}{a_1(t)}}, \ph t\ge t_0$. By (3.2) from here it follows that every solution  $(\phi(t), \psi(t))$  of the system (1.1) with  $\frac{\psi(t_0)}{\phi(t_0)} \in \Bigl[-\sqrt{\frac{c_2(t)}{a_2(t)}}; \sqrt{\frac{c_1(t)}{a_1(t)}}\Bigr]$  is $t_0$-regular and the relation (3.24) is fulfilled.  If  $\frac{\psi(t_0)}{\phi(t_0)} \in \Bigl(-\sqrt{\frac{c_2(t)}{a_2(t)}}; \sqrt{\frac{c_1(t)}{a_1(t)}}\Bigr)$,  then  $x(t) \equiv \frac{\psi(t)}{\phi(t)}, \ph t\ge t_0$
is $t_0$-normal. Hence (see [10]) $\widetilde{\mu}_x(t_0;t)$  is bounded. Therefore since by (3.2) $\widetilde{\mu}_x(t_0;t) = I_\phi(t), \ph t\ge t_0$, the function  $I_\phi(t), \ph t\ge t_0$  is bounded. The theorem is proved.

Example 3.6. Consider the system of equations
$$
\sist{\phi' = t^2 \phi + t \sin t \psi;}{\psi' = t^3 \cos t \phi - t^2 \psi, \ph t\ge 1.} \eqno (3.25)
$$
Any solution  $(\phi(t), \psi(t))$  of this system with  $\Bigl|\frac{\psi(t_0)}{\phi(t_0)}\Bigr| \le 1$  is $1$-regular and
$$
\Bigl|\frac{\psi(t)}{\phi(t)}\Bigr| \le t, \phh t \ge 1. \eqno (3.26)
$$
If  $\Bigl|\frac{\psi(t_0)}{\phi(t_0)}\Bigr| < 1$, then the function  $I_\phi(t) \equiv \il{1}{t} \frac{\tau \sin\tau}{\phi^2(\tau)}, \ph t\ge 1,$ is bounded. To prove this fact we should take   $a_j(t) \equiv (-1)^j t, \ph c_j(t) \equiv (-1)^j t^3, \ph t\ge 1, \ph j=1,2,$ and apply Theorem 3.6. From the firs equation of the system (3.25) and from (3.26) is seen that for every solution $(\phi(t), \psi(t))$ of the system (3.25) with   $\Bigl|\frac{\psi(t_0)}{\phi(t_0)}\Bigr| \le 1$  the function $|\phi(t)|$ is monotonically nondecreasing.

{\bf 3.4. The case when $a_{12}(t)$  and  $a_{21}(t)$  do not change their signs}.  In this section we will study the behavior of the system (1.1) in two cases: I. $a_{12}(t) \ge 0, \ph a_{21}(t) \ge 0, \ph t\ge t_0$; II. $a_{12}(t) \ge 0, \ph a_{21}(t) \le 0, \ph t\ge t_0$. The cases  III. $a_{12}(t) \le 0, \ph a_{21}(t) \le 0, \ph t\ge t_0$; IV. $a_{12}(t) \le 0, \ph a_{21}(t) \ge 0, \ph t\ge t_0$  can be reduced to  I  and  II
by simple transformation $\phi \to - \phi$  in the system (1.1) In  (1.1) make the substitutions
$$
\sist{\phi = J_{a_{11}}(t)\rho\sin\theta;}{\psi = J_{a_{22}}(t)\rho\cos\theta.} \eqno (3.27)
$$
We  get:
$$
\sist{\rho' \sin \theta + \theta'\rho  \cos \theta = A_{21}(t) \rho \cos \theta;}{\rho' \cos \theta - \theta'\rho  \sin \theta = A_{12}(t) \rho \sin \theta,} \eqno (3.28)
$$
where  $A_{21}(t)\equiv J_{-B}(t)a_{21}(t), \ph A_{12}(t) \equiv J_B(t)a_{12}(t), \ph t\ge t_0$.                             Multiplying the first equation of the last system on  $\sin \theta$, and the second one on   $\cos \theta$ and summarizing  we obtain the equation
$$
\rho' = [A_{12}(t) + A_{21}(t)] \rho \sin \theta \cos \theta, \phh t\ge t_0. \eqno (3.29)
$$
Now multiplying the first equation of the system (3.28) on   $\cos \theta$, and the second one on  $\sin \theta$  and subtracting from the first obtained the second one we  get:
$$
\theta' \rho = \rho [A_{21}(t) \cos^2 \theta - A_{12}(t) \sin^2 \theta], \phh t\ge t_0. \eqno (3.30)
$$
The system (3.28) (or the system of equations (3.29) and  (3.30)) is equivalent to the system (1.1) in the sense that to any real valued nontrivial solution  $(\phi(t), \psi(t))$ of the system (1.1) corresponds the solution $(\rho(t), \theta(t))$ of the system (3.26) with  $\rho(t) > 0, \ph t\ge t_0$, defined by formulae (3.25), or by formulae:
$$
\rho(t) = \sqrt{J_{-2a_{11}}(t)\phi^2(t) + J_{-2a_{22}}(t) \psi^2(t)}, \phh
\theta(t) = \sist{\arctan \frac{\phi(t)}{\psi(t)} J_{-B}(t), \ph \psi(t) \ne 0;}{\arctan \frac{\psi(t)}{\phi(t)} J_{B}(t), \ph \phi(t) \ne 0,}
$$
(to null solution $(0,0)$ of the system  (1.1)  corresponds the solution  $(0,0)$ of the system (3.28)).
Due to uniqueness  theorem if $\rho (t_0) > 0$,  then
$$
\rho (t)  > 0, \phh t\ge t_0. \eqno (3.31)
$$
Then  Eq. (3.30) can be replaced by
$$
\theta'  = A_{21}(t) \cos^2 \theta - A_{12}(t) \sin^2 \theta, \phh t\ge t_0. \eqno (3.32)
$$

{\bf Definition 3.14}. {\it  We will say that a continuous function  $u(t)$ on  $[t_1;+\infty) \ph (t_1 \ge t_0)$  vanishes no more than one time on  $[t_1; +\infty)$, if it has no nulls on  $[t_1; +\infty)$, or the set of its nulls on $[t_1; +\infty)$ is a segment $[a;b]$ which can be degenerated in a single point.}

{\bf Definition 3.15}. {\it We will say that the system (1.1) satisfies the non conjugation condition on  $[t_1;+\infty) \ph (t_1 \ge t_0)$, if for its each nontrivial solution  $(\phi(t), \psi(t))$  the functions  $\phi(t)$  and  $\psi(t)$  vanish on  $[t_1; +\infty)$ no more than one tome.}

{\bf 3.4.1. The case when $a_{12}(t) \ge 0, \ph a_{21}(t) \ge 0, \ph t\ge t_0$.} In this subsection we will study the behavior of  solutions of the system (1.1) under the restrictions  $a_{12}(t) \ge 0, \linebreak a_{21}(t) \ge~ 0, \ph t\ge t_0$.

{\bf Theorem  3.7.} {\it  If  $a_{12}(t) \ge 0, \ph a_{21}(t) \ge 0, \ph t\ge t_0$,
then the system (1.1) satisfies the non conjugation condition on  $[t_0; +\infty)$.}

Proof. Consider the following family of solutions of the system (1.1) $\Omega (t) \equiv \{ (\phi(t), \psi(t)) : \phi(t_0) = 1\}$. By virtue of Theorem 2.3. II$^\circ$ Eq. (3.1) has the unique  $t_0$-extremal solution $x_*(t)$.  Let  $(\phi(t), \psi(t)) \in \Omega$.  If $\psi(t_0) \ge x_*(t_0)$, then by (3.2) and   Lemma 2.1 the function $x(t) \equiv \frac{\psi(t)}{\phi(t)}$ is a  $t_0$-regular solution of Eq. (3.1).  Therefore if  $\psi(t_0) > 0$,  then $\psi(t) > 0, \ph t\ge t_0$    (since  $\phi(t) > 0, \ph t\ge t_0$);  if $\psi(t_0) = 0$, then there exists $t_1 \ge t_0$  such that $\psi(t) = 0, \ph t\in [t_0;t_1]$  and  $ \psi(t) > 0, \ph t> t_1$. Suppose  $I_{a_{12},B}^+(t_0) = +\infty$ or  $I_{a_{21},-B}^+(t_0) = +\infty$. Then by virtue of Theorem 2.3. III$^\circ$ if  $\psi(t_0) \in (x_*(t_0); 0)$, then there exists $t_2 \ge t_1 > t_0$    such that $\psi(t) < 0, \ph t\in [t_0;t_1), \ph \psi(t) = 0, \ph t \in [t_1;t_2], \ph \psi(t) > 0, \ph t > t_2$ . If $\psi(t_0) = x_*(t_0)$, then by virtue of Theorem 2.3. III$^\circ$  and  (3.2) we have $\psi(t) < 0, \ph t \ge t_0$.  From here it follows that  $u_*(t) \equiv \frac{\phi(t)}{\psi(t)}, \ph t\ge t_0,$ is a   $t_0$-regular solution of the equation
$$
u' + a_{21}(t) u^2 - B(t) u - a_{21}(t) = 0, \phh t\ge t_0. \eqno (3.33)
$$
By Lemma 2.1 from here it follows that for   $\psi(t_0) < x_*(t_0)$ the function $u(t) = \frac{\phi(t)}{\psi(t)}$ is a $t_0$-regular solution of the last equation. Hence for  $\psi(t_0) < 0$  we have  $\psi(t) < 0, \ph t\ge t_0$.    Suppose  $I_{a_{12},B}^+(t_0) < +\infty, \phh I_{a_{21},-B}^+(t_0) < +\infty$. Then by virtue of Theorem  2.3. IV$^\circ$  Eq. (3.1) has a negative $t_0$-normal solution.  By analogy of previous it shows that for  $(\phi(t), \psi(t)) \in \Omega$; for  $\psi(t_0) > 0, \ph \psi(t) > 0, \ph t\ge t_0$: for $\psi(t_0) = 0$  there exists  $t_1 \ge t_0$ such that  $\psi(t) = 0, \ph t\in [t_0;t_1], \ph \psi(t) > 0, \ph t > t_1$;  for $\psi(t_0) \in (x_N^-(t_0); 0)$ there exists $t_2 \ge t_1 > t_0$ such that $\psi(t) < 0, \ph t\in [t_0;t_1), \ph \psi(t) = 0,\ph t\in [t_1;t_2], \ph \psi(t) > 0, \ph t > t_2$;  for  $\psi(t_0) \le x_N^-(t_0)$  the inequality  $\psi(t) < 0, \ph t \ge t_0$ is fulfilled.           So we have shown that for $(\phi(t), \psi(t)) \in\Omega$ the function  $\psi(t)$
vanishes on  $[t_0;+\infty)$ no more than one time.  From here it follows that for arbitrary real valued solution $(\phi(t), \psi(t))$ of the system (1.1) with  $\phi(t_0) \ne 0$ the function  $\psi(t)$ vanishes on $[t_0;+\infty)$ no more than one time. Let $(\phi_0(t), \psi_0(t))$ be a nontrivial solution of the system  (1.1) with   $\phi_0(t_0) = 0$.  Then there exists $t_2 > t_1 \ge t_0$ such that $\phi_0(t) =0, \ph t\in [t_0;t_1]$ and  $\phi_0(t) \ne 0, \ph t\in(t_1;t_2)$. Indeed if  $\phi_0(t) \equiv 0$  on  $[t_0;+\infty)$,      then  $u_0(t) \equiv 0$ is a solution of Eq. (3.33), which is impossible. By uniqueness theorem we have $\psi_0(t) \ne 0, \ph t\in [t_0; t_3]$, for some  $t_3\in (t_1, t_2)$. Then by already proven the function $\psi_0(t)$  vanishes on  $[t_3;+\infty)$ no more than one time.  Therefore $\psi_0(t)$  vanishes on  $[t_0;+\infty)$ no more than one time. We have proved that for each real valued solution  $(\phi(t), \psi(t))$ of the system (1.1) the function $\psi(t)$  vanishes on  $[t_0;+\infty)$ no more than one time. Replacing the system (1.1) by its transpose from here we conclude that for each real valued solution  $(\phi(t), \psi(t))$ of the system (1.1) the function $\phi(t)$  vanishes on $[t_0;+\infty)$ no more than one time.  The theorem is proved.

Let us compare Theorem 3.7 with the following result of Vale - Pussen  (see [18], p. 122;  [1], p. 212).

{\bf Proposition  3.2}. {\it  If in the considering interval  $[t_1;t_2] \ph p(t) \equiv 1, \ph |q(t)| \le M_1$ and  $|r(t)| \le M_2$, then in the considering interval  $[t_1;t_2]$ the distance $h$  between eny two successive nulls of each solution od Eq.  (1.2) satisfies:

\noindent
for  $M_2 > 0$  the inequality
$$
h \ge \frac{\sqrt{9 M_1^2 + 24 M_2} - 3 M_1}{2M_2},
$$
and for  $M_2=0$  the inequality
$$
h \ge \frac{2}{M_2}.
$$
}
By virtue of Theorem 3.7 for  $r(t) \equiv 0$ the system (1.4) satisfies  of non conjugation condition.  В  (1.3) from here it follows that for $r(t) \equiv 0$  every solution of Eq. (1.2) has no more than one null on $[t_0;+\infty)$.  Therefore all equations (1.2), satisfying the conditions of Proposition 3.2 and hawing for $M_2=0$  two successive nulls on $[t_1; t_2]$, form a empty set.

{\bf Theorem 3.8}. {\it Suppose  $a_{12}(t) \ge 0, \ph a_{21}(t) \ge 0, \ph t\ge t_0$                                 Then if:

\noindent
A$_2$) $I_{a_{12}, B}^+(t_0) = +\infty$ or  $I_{a_{21}, -B}^+(t_0) = +\infty$

\noindent
then there exists up to arbitrary multiplier the unique solution  $(\phi_*(t), \psi_*(t))$ of the system (1.1) such that  $\phi_*(t) > 0, \ph \psi_*(t) < 0, \ph t\ge t_)$; for any solution  $(\phi(t), \psi(t))$  of the system (1.1), linearly independent of $(\phi_*(t), \psi_*(t))$ the equalities
$$
\liml{t\to +\infty} \frac{\phi_*(t)}{\phi(t)} = \liml{t\to +\infty} \frac{\psi_*(t)}{\psi(t)} = 0, \eqno (3.34)
$$
are fulfilled and
$$
\ilpp \frac{a_{12}(\tau) J_S(\tau)}{\phi_*^2(\tau)} d \tau = \ilpp \frac{a_{21}(\tau) J_S(\tau)}{\psi_*^2(\tau)} d \tau = + \infty, \eqno (3.35)
$$
$$
\ilp{T} \frac{a_{12}(\tau) J_S(\tau)}{\phi^2(\tau)} d \tau < +\infty, \phh \ilp{T} \frac{a_{21}(\tau) J_S(\tau)}{\psi^2(\tau)} d \tau < + \infty, \eqno (3.36)
$$
where $T=T(\phi;\psi) \ge t_)$ such that $\phi(t) \ne 0, \ph \psi(t) \ne 0, \ph t\ge T$, and

\noindent
a$_1$) if $I_{a_{12}, B}^+(t_0) = +\infty$    ($I_{a_{21}, -B}^+(t_0) = +\infty$), then
$$
\ilpp a_{12}(\tau)\frac{\phi_*(\tau)}{\psi_*(\tau)} d\tau = - \infty, \phh \ilp{T} a_{12}(\tau)\frac{\psi(\tau)}{\phi(\tau)} d\tau = +\infty  \eqno (3.37)
$$
$$
\biggl(\ilpp a_{21}(\tau)\frac{\psi_*(\tau)}{\phi_*(\tau)} d\tau = - \infty, \phh \ilp{T} a_{21}(\tau)\frac{\phi(\tau)}{\psi(\tau)} d\tau = +\infty\biggr);  \eqno (3.38)
$$
b$_1$) if    $I_{a_{12}, B}^+(t_0) = +\infty, \ph  \ilpp a_{21}(\tau) I_{- B, a_{12}}^-(t_0;\tau) d\tau < +\infty \\ \phantom{aaaaaaaaaaaaaaaaaaaaa} (I_{a_{21}, -B}^+(t_0) = +\infty, \ph  \ilpp a_{12}(\tau) I_{B, a_{21}}^-(t_0;\tau) d\tau < +\infty)$, then
$$
\ilpp a_{12}(\tau)\biggl|\frac{\psi_*(\tau)}{\phi_*(\tau)}\biggr| d \tau < +\infty, \phh \ilpp a_{21}(\tau)\frac{\phi_*(\tau)}{\psi_*(\tau)} d \tau = - \infty \eqno (3.39)
$$
$$
\biggl(\ilpp a_{21}(\tau)\biggl|\frac{\phi_*(\tau)}{\psi_*(\tau)}\biggr| d \tau < +\infty, \phh \ilpp a_{12}(\tau)\frac{\psi_*(\tau)}{\phi_*(\tau)} d \tau = - \infty\biggr)  \eqno (3.40)
$$

\noindent
B$_2$) $I_{a_{12}, B}(t_0;+\infty) < +\infty, \phh I_{a_{21}, -B}(t_0;+\infty) < +\infty$
then there exist up to arbitrary multiplier the unique linearly independent solutions  $(\phi_{**}(t), \psi_0(t)),\linebreak (\phi_0(t), \psi_{**}(t))$  of the system  (1.1) such that $\phi_{**}(t) > 0, \ph \psi_{**}(t) > 0, \ph t\ge t_0$; for each solution $(\phi(t), \psi(t))$ of the system (1.1), linearly independent of  $(\phi_{**}(t), \psi_0(t))$  \linebreak(of   $(\phi_0(t), \psi_{**}(t))$) the equality
$$
\liml{t\to +\infty}\frac{\phi_{**}(t)}{\phi(t)} = 0 \phh\biggl(\liml{t\to +\infty}\frac{\psi_{**}(t)}{\psi(t)} = 0 \biggr), \eqno (3.41)
$$
is fulfilled and
$$
\ilpp \frac{a_{12}(\tau) J_S(\tau)}{\phi_{**}^2(\tau)} d \tau = \ilpp \frac{a_{21}(\tau) J_S(\tau)}{\psi_{**}^2(\tau)} d \tau = +\infty \eqno (3.42)
$$
$$
\ilp{T} \frac{a_{12}(\tau) J_S(\tau)}{\phi^2(\tau)} d \tau < +\infty\phh \biggl( \ilp{T} \frac{a_{21}(\tau) J_S(\tau)}{\psi^2(\tau)} d \tau < +\infty\biggr),  \eqno (3.43)
$$
where  $T =T(\phi;\psi) \ge t_0$ such that $\phi(t) \ne 0, \ph \psi(t) \ne 0, \ph t\ge T$.
}

Proof.  Let  $x_*(t)$  be the  $t_0$-extremal solution of Eq. (3.1). Let us prove  A$_2$).  Above it was proved (in proving of Theorem 3.7),  that for fulfillment of the conditions  A$_2$) the function  $u_*(t) \equiv \frac{1}{x_*(t)}$  is the  $t_0$-extremal solution of Eq. (3.33). By  (3.2) the functions
$$
\phi_*(t) \equiv \exp\biggl\{ \il{t_0}{t}\Bigl[a_{12}(\tau) x_*(\tau) + a_{11}(\tau)\Bigr]d\tau\biggr\}, \ph \psi_*(t) \equiv x_*(t) \phi_*(t), \ph t\ge t_0,
$$
form the solution  $(\phi_*(t), \psi_*(t))$  of the system (1.1).  Let  $x_0(t)$ be a solution of Eq. (3.1) with  $x_0(t_0) > x_*(t_0)$. By Lemma 2.1 $x_0(t)$ is $t_0$-normal.  Then by virtue of Theorem 2.2 the integral $\widetilde{\nu}_{x_0}(t_0)$  converges. Obviously  $\widetilde{\nu}_{x_0}(t) > 0, \ph t \ge t_0$   (since $a_{12}(t) \ge 0, \ph t\ge t_0$ and has unbounded support). By (3.2) the functions
$$
\widetilde{\phi}_0(t) \equiv u_*(t) \widetilde{\psi}_*(t), \ph \widetilde{\psi}_*(t) \equiv \exp\biggl\{ \il{t_0}{t}\Bigl[a_{21}(\tau) u_*(\tau) + a_{22}(\tau)\Bigr]d\tau\biggr\}
$$
form the solution $(\widetilde{\phi}_*(t), \widetilde{\psi}_*(t))$ of the system (1.1). Since $u_*(t) = \frac{1}{x_*(t)}$, we have \linebreak $\det\left(\begin{array}{l} \phi_*(t) \ph \psi_*(t)\\
 \widetilde{\phi}_*(t) \ph  \widetilde{\psi}_*(t)\end{array}\right) = 0$.
Then  $(\phi_*(t), \psi_*(t))$ and $(\widetilde{\phi}_*(t), \widetilde{\psi}_*(t))$ are linearly independent. let $u_0(t)$ be a solution of Eq. (3.33) with $u_0(t_0) > u_*(t_0)$. By Lemma 2.1 $u_0(t)$ is $ t_0$-normal. Then by virtue of Theorem 2.2 the integral
$$
\widetilde{\widetilde{\nu}}_{u_0}(t) \equiv \ilp{t} a_{21}(\tau) \exp\biggl\{- \il{t}{\tau}\Bigl[2 a_{21}(\xi) u_0(\xi)  - B(\xi)\Bigr]d\xi\biggr\}d\tau
$$
converges for all  $t\ge t_0$  and  $\widetilde{\widetilde{\nu}}_{u_0}(t) > 0, \ph t\ge t_0$.  Then the assertion  A$_2$) immediately follows from the linearly independence of solutions  $(\phi_*(t), \psi_*(t))$ and $(\widetilde{\phi}_*(t), \widetilde{\psi}_*(t))$
and from Theorem 3.5   (as a up to arbitrary multiplier solution of Eq.  (1.1), appearing  in the formulation of A$_2$)  one can take $(\phi_*(t), \psi_*(t))$).  Let us prove  B$_2$).  By virtue of \linebreak Theorem 2.3. VI$^\circ$ Eq. (3.1) has a $t_0$-normal negative solution $x_N^-(t)$ such that every solution  $x(t)$ of Rq.  (3.1) with  $0\ge x(t_0) > x_N^-(t_0)$ vanishes on  $(t_0;+\infty)$.  From here it follows that  $u_*(t)\equiv \frac{1}{x_N^-{t}}, \ph t\ge t_0$,  is the $t_0$-extremal solution to Eq. (3.33). By  (3.2) the functions
$$
\phi_{**}(t) \equiv \exp\biggl\{ \il{t_0}{t}\Bigl[a_{12}(\tau) x_*(\tau) + a_{11}(\tau)\Bigr]d\tau\biggr\}, \ph \psi_0(t) \equiv x_*(t) \phi_{**}(t),
$$
$$
\phi_0(t) \equiv u_*(t)\psi_{**}(t), \ph \psi_{**}(t) \equiv \exp\biggl\{ \il{t_0}{t}\Bigl[a_{21}(\tau) u_*(\tau) + a_{22}(\tau)\Bigr]d\tau\biggr\} \ph t\ge t_0,
$$
form the solutions  $(\phi_{**}(t),\psi_0(t))$   and   $(\phi_0(t), \psi_{**}(t))$ of the system (1.1). Since \linebreak $\det\left(\begin{array}{l} \phi_{**}(t_0) \ph \psi_0(t_0)\\
 \widetilde{\phi}_0(t_0) \ph  \widetilde{\psi}_{**}(t_0)\end{array}\right) =1 - x_*(t_0) u_*(t_0) = 1 - \frac{x_*(t_0)}{x_N^-(t_0)} \ne 0$, the solutions
 $(\phi_{**}(t),\psi_0(t))$   and   $(\phi_0(t), \psi_{**}(t))$ are linearly independent. Let  $x_0(t)$ be a solution of Eq.  (3.1) with  $x_0(t_0) > x_*(t_0)$.  According to  Lemma 2.1 $x_0(t)$ is $t_0$-normal. Then by virtue of Theorem~ 2.3 the integral $\widetilde{\nu}_{x_0}(t)$ converges for all $t \ge t_0$ and  $\widetilde{\nu}_{x_0}(t) > 0, \ph t\ge t_0$. By analogy if $u_0(t)$ is a solution of Eq. (3. 33) with  $u_0(t_0) > u_*(t_0)$,  then \\ $\widetilde{\widetilde{\nu}}_{u_0}(t) \equiv\ilp{t} a_{21}(\tau) \exp\biggl\{- \il{t}{\tau}\Bigl[2 a_{21}(\xi) u_0(\xi) - B(\xi)\Bigr] d\xi\biggr\}d\tau$ for all $t\ge t_0$  and   $\widetilde{\widetilde{\nu}}_{u_0}(t) \ge~0,\linebreak t\ge t_0$.  From here From Theorem 3.5 and from linearly independence of solutions \linebreak $(\phi_{**}(t),\psi_0(t))$  and  $(\phi_0(t), \psi_{**}(t))$ it follows  B$_2$).   The theorem is proved.

On the basis of Theorems 2.3 and  3.5 and the firs of formulae (3.2) one can make the phase portrait of the family of curves  $\Phi(t)\equiv J_{a_{11}}(t) \phi(t)$, where  $(\phi(t), \psi(t))$  runs  the family  $\Omega$  in the case when  $a_{12}(t) \ge 0, \ph, a_{21}(t) \ge 0, \ph t\ge t_0$ (see pict. 9 and pict.  10)

\vskip 20pt

\begin{picture}(140,100)
\put(30,-40){\vector(0,1){140}}
\put(-10,0){\vector(1,0){210}}
\put(240,-40){\vector(0,1){140}}
\put(220,00){\vector(1,0){210}}
\put(195,-7){$_t$}
\put(60,5){$_{t_0}$}
\put(70,0){\circle*{3}}
\put(30.5,35){\circle*{3}}
\put(180,10){$_{x_*(t)}$}
\put(155,-10){$_{\alpha_*}$}

\put(22,35){$_1$}
\put(35,85){$_{\phi(t)}$}
\put(427,-7){$_t$}
\put(270,5){$_{t_0}$}
\put(280,0){\circle*{3}}
\put(240.5,35){\circle*{3}}
\put(405,6){$_{x_*(t)}$}
\put(402,23){$_{x_N^-(t)}$}
\put(365,-10){$_{\alpha_*}$}
\put(412,-5){$_{\alpha_N^-}$}

\put(232,35){$_1$}
\put(245,85){$_{\phi(t)}$}
\put(-10,-60){$_{Pict.9.\ph \alpha_* = arctan \ph x_*(t_0)}$}
\put(15,-70){$_{I_{a_{12},B}^+(t_0) = +\infty\ph  or \ph  I_{a_{21},-B}^+(t_0) = +\infty }$}

\put(220,-60){$_{Pict.10.\ph \alpha_* = arctan \ph x_*(t_0),\ph \alpha_N^- = arctan \ph x_N^-(t_0)}$}
\put(245,-70){$_{I_{a_{12},B}^+(t_0) < +\infty\ph  and \ph  I_{a_{21},-B}^+(t_0) < +\infty }$}

\put(150,-13){\qbezier(0,0)(5,3)(6,13)}
\put(151,-13){\vector(-1,-1){1}}
\put(157,0){\vector(1,3){.5}}

\put(360,-13){\qbezier(0,0)(5,3)(6,13)}
\put(405,-13){\qbezier(0,2)(5,3)(6,13)}
\put(361,-12){\vector(-1,-1){1}}
\put(367,0){\vector(1,3){.5}}
\put(407,-9){\vector(-1,-1){1}}
\put(412,0){\vector(1,3){.5}}

\put(70,10){\thicklines \qbezier[30](0,-40)(0,25)(0,90)}

\put(280,10){\thicklines \qbezier[30](0,-40)(0,25)(0,90)}

\put(70,10){\thicklines \qbezier[30](-40,25)(40,25)(120,25)}

\put(280,10){\thicklines \qbezier[30](-40,25)(50,25)(140,25)}

\put(110,0){\thicklines \qbezier[30](-40,35)(20,0)(80,-35)}

\put(320,0){\thicklines \qbezier[30](-40,35)(20,0)(80,-35)}

\put(320,10){\thicklines \qbezier[30](-40,25)(30,0)(100,-25)}

\put(70,10){\thicklines \qbezier(0,25)(35,2)(125,-5)}
\put(70,10){\thinlines\qbezier(0,25)(35,-7)(125,-20)}
\put(70,10){\thinlines\qbezier(0,25)(35,-17)(125,-35)}
\put(70,10){\thinlines\qbezier(0,25)(35,-27)(125,-50)}

\put(70,10){\thinlines\qbezier(0,25)(35,3)(85,20) \qbezier(85,20)(100,28)(125,30)}
\put(70,10){\thinlines\qbezier(0,25)(35,6)(85,30) \qbezier(85,30)(100,38)(125,40)}
\put(70,10){\thinlines\qbezier(0,25)(35,8)(85,37) \qbezier(85,37)(100,45)(125,48)}
\put(70,10){\thinlines\qbezier(0,25)(35,12)(85,44) \qbezier(85,44)(100,52)(125,56)}

\put(280,10){\thicklines \qbezier(0,25)(35,2)(135,-8)}
\put(280,10){\thicklines \qbezier(0,25)(45,12)(135,5)}
\put(280,10){\thinlines \qbezier(0,25)(45,5)(135,-2)}
\put(280,10){\thinlines\qbezier(0,25)(35,-7)(125,-20)}
\put(280,10){\thinlines\qbezier(0,25)(35,-17)(125,-35)}
\put(280,10){\thinlines\qbezier(0,25)(35,-27)(125,-50)}
\put(280,10){\thinlines\qbezier(0,25)(55,13)(85,25) \qbezier(85,25)(100,34)(135,35)}
\put(280,10){\thinlines\qbezier(0,25)(55,16)(85,30) \qbezier(85,30)(100,41)(135,45)}
\put(280,10){\thinlines\qbezier(0,25)(55,18)(85,37) \qbezier(85,37)(100,48)(135,53)}
\put(280,10){\thinlines\qbezier(0,25)(55,22)(85,44) \qbezier(85,44)(100,57)(135,61)}

\end{picture}

\vskip 80pt
By analogy it can be  make the phase portrait for the family of curves  $\Psi(t)\equiv J_{a_{22}}(t)\psi(t)$,  where $(\phi(t), \psi(t))$  runs the family of solutions of the system (1.1) with  $\psi(t_0) =1$.

Let  $(\phi(t), \psi(t)) \in \Omega$. If the curve $(\phi(t), \psi(t)) $  is in the first or the third quadrant of phase plane $\Phi\equiv J_{-a_{11}}(t)\phi, \ph  \Psi\equiv J_{-a_{11}}(t)\psi$, then from  (3.39) it follows that the function  $\rho(t) \equiv \sqrt{J_{-2a_{11}}(t)\phi^2(t) +  J_{-2a_{22}}(t)\psi^2(t)}$  monotonically non decreases, and if  $(\phi(t), \psi(t)) $ belongs to the second or the fourth quadrant of phase plane  $\Phi\equiv J_{-a_{11}}(t)\phi, \ph  \Psi\equiv J_{-a_{11}}(t)\psi$,                                  then  $\rho(t)$  monotonically non increases.  Then taking into account Theorem 3.7 one can make the phase portrait of the family of the curves  $\widetilde{\Omega}\equiv \{( J_{-a_{11}}(t)\phi(t),J_{-a_{11}}(t)\psi(t)), \ph (\phi(t), \psi(t)) \in\Omega\}$ for the case  $a_{12}(t) \ge 0, \ph a_{21}(t) \ge 0, \ph t\ge t_0$ $\bigl($see pict.'s  11 -  14, where
$$
\Phi_*(t) \equiv \exp\biggl\{ \il{t_0}{t} a_{12}(\tau) x_*(\tau) d\tau\biggr\}, \phh \Psi_*(t) \equiv x_*(t) J_{a_{11}}(t)\Phi_*(t)\bigr).
$$
If  $I_{a_{12},B}^+(t_0)= +\infty$  and  $I_{a_{21},-B}^+(t_0)= +\infty$, then by virtue
of Theorem .3.V$^\circ$ we have $ \Phi_*(t) \to 0, \ph \Psi_*(t) \to 0$   for  $t\to +\infty$ (see pict. 11). If  $I_{a_{12},B}^+(t_0)= +\infty, \linebreak \ilpp a_{21}(\tau) I_{-B,a_{12}}^-(t_0;\tau) d\tau  < +\infty$,
 then by virtue of Theorem 2.3.~VII$^\circ$ we have \linebreak $\Phi_*(t) \to~\gamma_* > 0, \ph \Psi_*(t) \to 0$  for  $t\to +\infty$ (see pict. 12).

\vskip 100pt
\begin{picture}(1,1)
\put(70,-40){\vector(0,1){140}}
\put(-10,20){\vector(1,0){210}}
\put(280,-40){\vector(0,1){140}}
\put(220,20){\vector(1,0){210}}
\put(190,13){$_\Phi$}
\put(0,-60){$_{Pict.11.}$}
\put(216,-60){$_{Pict.12.}$}
\put(415,13){$_\Phi$}
\put(115,15){$_1$}
\put(325,15){$_1$}
\put(110,20){\circle*{3}}
\put(320,20){\circle*{3}}
\put(75, 90){$_\Psi$}
\put(285, 90){$_\Psi$}
\put(123,0){$_{(\Phi_*(t),\Psi_*(t))}$}
\put(120,0){\thicklines\qbezier[5](0,0)(-10,-1)(-20,-2)}
\put(100,-2.6){\vector(-4,-1){1}}

\put(333,0){$_{(\Phi_*(t),\Psi_*(t))}$}
\put(330,0){\thicklines\qbezier[5](0,0)(-10,-1)(-20,-2)}
\put(310,-2.6){\vector(-4,-1){1}}
\put(33,-58){$_{I_{a_{12},B}^+(t_0)=I_{a_{21},-B}^+(t_0)=+\infty}$}
\put(246,-52){$_{I_{a_{12},B}^+(t_0)=+\infty,  \ilpp a_{21}(\tau)I_{-B,a_{12}}^-(t_0;\tau)d\tau<+\infty.}$}

\put(110,10){\thicklines \qbezier[30](0,-40)(0,25)(0,90)}
\put(320,10){\thicklines \qbezier[30](0,-40)(0,25)(0,90)}

\put(100,10){\thicklines \qbezier(10,-20)(5,-15)(0,-14.5)\qbezier(0,-14.5)(-5,-13.5)(-10,-9.5)}
\put(80,20){\thicklines \qbezier(10,-20)(5,-15)(0,-14.5) \qbezier(0,-14.5)(-5,-13.5)(-8,-3.5)}
\put(72,17){\vector(-1,3){.5}}
\put(100,15){\thinlines \qbezier(10,-20)(5,-15)(0,-14.5)\qbezier(0,-14.5)(-3,-13.5)(-6,-9.5)}
\put(85,25){\thinlines \qbezier(10,-20)(5,-15)(0,-14.5) \qbezier(0,-14.5)(-7.5,-10.5)(-8,-3.5)
\qbezier(-8,-3.5)(-7.5,7.5)(5,15.5)\qbezier(5,15.5)(20,30)(100,36)}
\put(185,61){\vector(4,1){1}}
\put(100,30){\thinlines \qbezier(10,-20)(1,-19)(0,-14.5)\qbezier(0,-14.5)(-1,-13.5)(-1.3,-9.5)}
\put(89,40){\thinlines \qbezier(10,-20)(11,-15)(20,-10) \qbezier(20,-10)(45,1)(100,2)}
\put(120,33.4){\vector(4,1){1}}
\put(89,50){\thinlines  \qbezier(20,-10)(45,1)(100,2)}
\put(120,43.4){\vector(4,1){1}}
\put(89,34){\thinlines  \qbezier(20,-10)(45,1)(100,2)}
\put(120,27.4){\vector(4,1){1}}
\put(89,67){\thinlines  \qbezier(20,-10)(45,-1)(100,2)}
\put(120,60.4){\vector(4,1){1}}
\put (100,4){\thinlines \qbezier (10,-20) (5,-16) (0,-14) \qbezier(0,-14)(-5,-13)(-10,-10)\qbezier(-10,-10)(-20,-2)(-30,-4)\qbezier(-30,-4)(-50,-7)(-90,-45) }
\put(30,-23.5){\vector(-1,-1){1}}
\put (100,-2){\thinlines \qbezier (10,-20) (5,-16) (0,-14) \qbezier(0,-14)(-5,-13)(-10,-10)\qbezier(-10,-10)(-20,-2)(-30,-4)\qbezier(-30,-4)(-50,-7)(-80,-40) }
\put(30,-31.5){\vector(-1,-1){1}}
\put (100,-8){\thinlines \qbezier (10,-20) (5,-16) (0,-14) \qbezier(0,-14)(-5,-13)(-10,-10)\qbezier(-10,-10)(-20,-2)(-30,-4)\qbezier(-30,-4)(-50,-7)(-70,-35) }
\put(34,-37.5){\vector(-2,-3){1}}

\put(310,10){\thicklines \qbezier(10,-20)(5,-15)(0,-14.5)\qbezier(0,-14.5)(-5,-13.5)(-10,-9.5)}
\put(290,20){\thicklines \qbezier(10,-20)(2.99,-18)(3,-1.5)}
\put(293,17){\vector(0,1){.5}}
\put(310,15){\thinlines \qbezier(10,-20)(5,-15)(0,-14.5)\qbezier(0,-14.5)(-3,-13.5)(-6,-9.5)}
\put(295,25){\thinlines \qbezier(10,-20)(5,-15)(6,-4.5)
\qbezier(6,-4.5)(6.5,10.5)(17,19.5)\qbezier(17,19.5)(30,33)(100,36)}
\put(395,61){\vector(4,1){1}}
\put(310,30){\thinlines \qbezier(10,-20)(1,-19)(0,-14.5)\qbezier(0,-14.5)(-1,-13.5)(-1.3,-9.5)}
\put(299,40){\thinlines \qbezier(10,-20)(11,-15)(20,-10) \qbezier(20,-10)(45,1)(100,2)}
\put(330,33.4){\vector(4,1){1}}
\put(299,50){\thinlines  \qbezier(20,-10)(45,1)(100,2)}
\put(330,43.4){\vector(4,1){1}}
\put(299,34){\thinlines  \qbezier(20,-10)(45,1)(100,2)}
\put(330,27.4){\vector(4,1){1}}
\put(299,67){\thinlines  \qbezier(20,-10)(45,-1)(100,2)}
\put(330,60.4){\vector(4,1){1}}
\put (310,4){\thinlines \qbezier (10,-20) (5,-16) (0,-14) \qbezier(0,-14)(-5,-13)(-10,-10)\qbezier(-10,-10)(-20,-2)(-30,-4)\qbezier(-30,-4)(-50,-7)(-90,-45) }
\put(240,-23.5){\vector(-1,-1){1}}
\put (310,-2){\thinlines \qbezier (10,-20) (5,-16) (0,-14) \qbezier(0,-14)(-5,-13)(-10,-10)\qbezier(-10,-10)(-20,-2)(-30,-4)\qbezier(-30,-4)(-50,-7)(-80,-40) }
\put(240,-31.5){\vector(-1,-1){1}}
\put (310,-8){\thinlines \qbezier (10,-20) (5,-16) (0,-14) \qbezier(0,-14)(-5,-13)(-10,-10)\qbezier(-10,-10)(-20,-2)(-30,-4)\qbezier(-30,-4)(-50,-7)(-70,-35) }
\put(244,-37.5){\vector(-2,-3){1}}

\end{picture}

\vskip 80pt
if $I_{a_{21},-B}^+(t_0)= +\infty$ and $\ilpp a_{12}(\tau) I_{B,a_{21}}^-(t_0;\tau) d\tau  < +\infty$, then by virtue of \linebreak Theorem 2.3.~VII$^\circ$ we have $\Phi_*(t) \to 0, \ph \Psi_*(t) \to \gamma_*< 0$               for  $t\to +\infty$ (see pict. 13).
If $I_{a_{12},B}(t_0) < +\infty$  and   $I_{a_{21},-B}(t_0) < +\infty$,
then by virtue of Theorem 2.3.~VII$^\circ$ we have $\Phi_{**}(t) \to 0, \ph \Psi_{**}(t) \to 0$   for   $t\to +\infty$,
where  $\Phi_{**}(t) \equiv J_{-a_{11}}(t)\phi_{**}(t), \ph \Psi_{**}(t) \equiv J_{-a_{22}}(t)\psi_{**}(t), \ph (\phi_{**}(t), \psi_0(t))$ and  $(\phi_0(t), \psi_{**}(t))$  are defined in the formulation of \linebreak Theorem 3.8. B$_2$) as  up to arbitrary multiplier the unique solutions of the system (1.1) (see pict.  14).

\vskip 100pt
\begin{picture}(1,1)
\put(80,-40){\vector(0,1){140}}
\put(-10,20){\vector(1,0){210}}
\put(290,-40){\vector(0,1){140}}
\put(220,20){\vector(1,0){210}}
\put(185,25){$_\Phi$}
\put(113,13){$_{1}$}
\put(110,20){\circle*{3}}
\put(83,90){$_{\Psi}$}

\put(415,25){$_\Phi$}
\put(323,13){$_{1}$}
\put(320,20){\circle*{3}}
\put(293,90){$_{\Psi}$}

\put(0,-60){$_{Pict.13.}$}
\put(123,-2){$_{(\Phi_*(t),\Psi_*(t))}$}
\put(120,-2){\thicklines\qbezier[5](0,0)(-10,-1)(-20,-2)}
\put(100,-4.6){\vector(-4,-1){1}}

\put(231,-60){$_{Pict.14.}$}
\put(343,9){$_{(\Phi_0(t),\Psi_{**}(t))}$}
\put(340,9){\thicklines\qbezier[5](0,0)(-10,-2)(-20,-4)}
\put(317,4){\vector(-4,-1){1}}

\put(333,-30){$_{(\Phi_{**}(t),\Psi_0(t))}$}
\put(330,-30){\thicklines\qbezier[10](0,0)(-10,13)(-20,26)}
\put(312,-5.6){\vector(-1,2){1}}

\put(3,-68){$_{I_{a_{21},-B}^+(t_0)=+\infty,  \ilpp a_{12}(\tau)I_{B,a_{21}}^-(t_0;\tau)d\tau<+\infty.}$}
\put(265,-58){$_{I_{a_{12},B}^+(t_0)<+\infty,\ph I_{a_{21},-B}^+(t_0)<+\infty}$}

\put(110,10){\thicklines \qbezier[30](0,-50)(0,30)(0,90)}
\put(320,10){\thicklines \qbezier[30](0,-50)(0,30)(0,90)}

\put(100,10){\thicklines \qbezier(10,-20)(5,-15)(0,-14.5)\qbezier(0,-14.5)(-5,-13.5)(-10,-9.5)}
\put(80,20){\thicklines \qbezier(10,-20)(2.99,-15)(0,-16)}
\put(80,4.5){\vector(-4,1){.5}}
\put(100,15){\thinlines \qbezier(10,-20)(5,-15)(0,-14.5)\qbezier(0,-14.5)(-3,-13.5)(-6,-9.5)}
\put(85,25){\thinlines \qbezier(10,-20)(5,-15)(6,-4.5)
\qbezier(6,-4.5)(6.5,10.5)(17,19.5)\qbezier(17,19.5)(30,33)(100,36)}
\put(185,61){\vector(4,1){1}}
\put(100,30){\thinlines \qbezier(10,-20)(1,-19)(0,-14.5)\qbezier(0,-14.5)(-1,-13.5)(-1.3,-9.5)}
\put(89,40){\thinlines \qbezier(10,-20)(11,-15)(20,-10) \qbezier(20,-10)(45,1)(100,2)}
\put(120,33.4){\vector(4,1){1}}
\put(89,50){\thinlines  \qbezier(20,-10)(45,1)(100,2)}
\put(120,43.4){\vector(4,1){1}}
\put(89,34){\thinlines  \qbezier(20,-10)(45,1)(100,2)}
\put(120,27.4){\vector(4,1){1}}
\put(89,67){\thinlines  \qbezier(20,-10)(45,-1)(100,2)}
\put(120,60.4){\vector(4,1){1}}
\put (100,4){\thinlines \qbezier (10,-20) (5,-16) (0,-14) \qbezier(0,-14)(-5,-13)(-10,-10)\qbezier(-10,-10)(-20,-2)(-30,-4)\qbezier(-30,-4)(-50,-7)(-90,-45) }
\put(30,-23.5){\vector(-1,-1){1}}
\put (100,-2){\thinlines \qbezier (10,-20) (5,-16) (0,-14) \qbezier(0,-14)(-5,-13)(-10,-10)\qbezier(-10,-10)(-20,-2)(-30,-4)\qbezier(-30,-4)(-50,-7)(-80,-40) }
\put(30,-31.5){\vector(-1,-1){1}}
\put (100,-8){\thinlines \qbezier (10,-20) (5,-16) (0,-14) \qbezier(0,-14)(-5,-13)(-10,-10)\qbezier(-10,-10)(-20,-2)(-30,-4)\qbezier(-30,-4)(-50,-7)(-70,-35) }
\put(34,-37.5){\vector(-2,-3){1}}

\put(310,20){\thicklines \qbezier(10,-20)(5,-15)(0,-14.5)\qbezier(0,-14.5)(-5,-13.5)(-9,-9.5)}
\put(290,30){\thicklines \qbezier(12,-20)(9,-19)(8,-10.5)}
\put(298,20){\vector(0,1){.5}}

\put(310,10){\thicklines \qbezier(10,-20)(5,-15)(0,-14.5)\qbezier(0,-14.5)(-5,-13.5)(-10,-9.5)}
\put(290,20){\thicklines \qbezier(10,-20)(2.99,-15)(0,-16)}
\put(290,4.5){\vector(-4,1){.5}}
\put(310,24){\thinlines \qbezier(10,-20)(5,-15)(2,-14.5)\qbezier(2,-14.5)(-1,-13.5)(-3,-9.5)}
\put(295,34){\thinlines \qbezier(12,-20)(7,-15)(13,-4.5)
\qbezier(13,-4.5)(14.5,.5)(25,6.5)}

\put(310,14){\thinlines \qbezier(10,-20)(5,-19)(-12,-8.5)\qbezier(-12,-8.5)(-17.8,-8.3)(-18,.5)}

\put(295,34){\qbezier(25,13.5)(50,25)(100,27)}
\put(395,61){\vector(4,1){1}}
\put(310,30){\thinlines \qbezier(10,-20)(1,-19)(0,-14.5)\qbezier(0,-14.5)(-1,-13.5)(-1.3,-9.5)}
\put(299,40){\thinlines \qbezier(10,-20)(11,-15)(20,-10) \qbezier(20,-10)(45,1)(100,2)}
\put(330,33.4){\vector(4,1){1}}
\put(299,50){\thinlines  \qbezier(20,-10)(45,1)(100,2)}
\put(330,43.4){\vector(4,1){1}}
\put(299,34){\thinlines  \qbezier(20,-10)(45,1)(100,2)}
\put(330,27.4){\vector(4,1){1}}
\put(299,67){\thinlines  \qbezier(20,-10)(45,-1)(100,2)}
\put(330,60.4){\vector(4,1){1}}
\put (310,4){\thinlines \qbezier (10,-20) (5,-16) (0,-14) \qbezier(0,-14)(-5,-13)(-10,-10)\qbezier(-10,-10)(-20,-2)(-30,-4)\qbezier(-30,-4)(-50,-7)(-90,-45) }
\put(240,-23.5){\vector(-1,-1){1}}
\put (310,-2){\thinlines \qbezier (10,-20) (5,-16) (0,-14) \qbezier(0,-14)(-5,-13)(-10,-10)\qbezier(-10,-10)(-20,-2)(-30,-4)\qbezier(-30,-4)(-50,-7)(-80,-40) }
\put(240,-31.5){\vector(-1,-1){1}}
\put (310,-8){\thinlines \qbezier (10,-20) (5,-16) (0,-14) \qbezier(0,-14)(-5,-13)(-10,-10)\qbezier(-10,-10)(-20,-2)(-30,-4)\qbezier(-30,-4)(-50,-7)(-70,-35) }
\put(244,-37.5){\vector(-2,-3){1}}
\end{picture}

\vskip 90pt

Suppose $(\phi(t), \psi(t)) \in \Omega$. Then by virtue of Theorem 2.3.I$^\circ$  if  $\psi(t_0) \ge \frac{-1}{I_{a_{12},B}^+(t_0)}$ then  $(\phi(t), \psi(t))$ is $t_0$-regular and
$$
\frac{\psi(t_0) J_{-B}(t)}{1 + \psi(t_0) I_{a_{12},B}^+(t_0;t)} \le \frac{\psi(t)}{\phi(t)} \le \psi(t_0) J_{-B}(t) + I_{B,a_{21}}^-(t_0;t), \ph t\ge t_0.
$$
Let us multiply all parts of this inequality on $a_{12}(t)$, add $a_{11}(t)$, integrate from  $t_0$  to  $t$ and take the exponential. We obtain:
$$
J_{a_{11}}(t)\Bigl[1 + \psi(t_0) I_{a_{12},B}^+(t_0;t)\Bigr] \le \phi(t) \le \phantom{aaaaaaaaaaaaaaaaaaaaaaaaaaaaaaaaaaaaaaaaaaaaaaaaa}
$$
$$
\phantom{aaaaa}\le \exp\Bigr\{\psi(t_0) I_{a_{12},B}^+(t_0;t) + \il{t_0}{t}[a_{11}(\tau) + a_{12}(\tau) I_{B, a_{21}}^-(t_0;\tau)]d\tau\Bigr\}, \phh t\ge t_0. \eqno (3.44)
$$
Passing to the transpose of the system (1.1) from here for the solution  $(\phi(t), \psi(t))$
of the system (1.1) with  $\phi(t_0) \ge \frac{-1}{I_{a_{21}, -B}^+(t_0)}, \ph \psi(t_0) =1$  we  get the estimates
$$
J_{a_{22}}(t)\Bigl[1 + \phi(t_0) I_{a_{121},-B}^+(t_0;t)\Bigr] \le \psi(t) \le \phantom{aaaaaaaaaaaaaaaaaaaaaaaaaaaaaaaaaaaaaaaaaaaaaaaaa}
$$
$$
\phantom{aaaaa}\le \exp\Bigr\{\phi(t_0) I_{a_{21},-B}^+(t_0;t) + \il{t_0}{t}[a_{22}(\tau) + a_{21}(\tau) I_{-B, a_{12}}^-(t_0;\tau)]d\tau\Bigr\}, \phh t\ge t_0. \eqno (3.45)
$$
From here and from (3.34) we immediately have

{\bf Theorem 3.9}. {\it  The system (1.1) with $a_{12}(t) \ge 0, \ph a_{21}(t) \ge 0, \ph t\ge t_0$, is Lyapunov  stable if and only if the functions
$J_{a_{11}}(t)I_{a_{12},B}^+(t_0;t), \ph J_{a_{22}}(t)I_{a_{21},-B}^+(t_0;t)$   are bounded on  $[t_0;+\infty)$.}

Example 3.7. Consider the system
$$
\sist{\phi' = \phantom{aaaaaa} \sin^2 t \psi;}{\psi' = \cos^2 t \phi, \phh t\ge t_0.} \eqno (3.46)
$$
for this system we have  $a_{12}(t) = \sin^2(t) \ge 0, \ph a_{21}(t) = \cos^2(t) \ge 0, \ph J_{a_{11}}(t)I_{a_{12},B}^+(t_0;t) =\linebreak = \il{t_0}{t}\sin^2(\tau) d\tau, \ph t\ge t_0$.                                                    Since $\il{t_0}{t}\sin^2(\tau) d\tau$ is unbounded by Theorem 3.9  the system (3.46) is not  Lyapunov  stable.
It is not difficult to verify that the application of Theorem 4.6.2 of Wazevski and mentioned above estimates of Lyapunov and Yu. S. Bogdanov to the system (3.46)  gives no result.

Let $(\phi_0(t), \psi_0(t))$  and $(\phi_1(t), \psi_1(t))$ be the solutions of the system  (1.1)
with   $\phi_0(t_0)=\psi_1(t_0)~=1,\ph \psi_0(t_0) = \phi_1(t_0) =0$. Then from (3.44) and (3.45) we obtain:
$$
\phi_0(t) \le \exp\Bigr\{ \il{t_0}{t}[a_{11}(\tau) + a_{12}(\tau) I_{B, a_{21}}^-(t_0;\tau)]d\tau\Bigr\},\phantom{aaaaaaaaaaaaaaaaaaaaaaaaaaaaaaaaaaaaaa}
$$
$$
\phantom{aaaaaaaaaaaaaaaaaaaaaaaaaa}\psi_0(t) \le \exp\Bigr\{\il{t_0}{t}[a_{22}(\tau) + a_{21}(\tau) I_{-B, a_{12}}^-(t_0;\tau)]d\tau\Bigr\}, \phh t\ge t_0.
$$
By (2.9) and  (3.2) from here for the solution $(\phi(t), \psi(t))$ of the system (1.1) with \linebreak  $\phi(t_0) > 0, \ph \psi(t_0) >~0$ we  get the estimates
$$
0< \phi(t) \le \phi(t_0) c_1 \exp\Bigr\{ \il{t_0}{t}[a_{11}(\tau) + a_{12}(\tau) I_{B, a_{21}}^-(t_0;\tau)]d\tau\Bigr\}, \eqno (3.47)
$$
$$
0 < \psi(t) \le \psi(t_0) c_2 \exp\Bigr\{\il{t_0}{t}[a_{22}(\tau) + a_{21}(\tau) I_{-B, a_{12}}^-(t_0;\tau)]d\tau\Bigr\}, \phh t\ge t_0, \eqno (3.48)
$$
where
$$
c_1\equiv \exp\biggl\{\ilpp a_{12}(\tau)\Bigl(\frac{\psi(\tau)}{\phi(\tau)} - \frac{\psi_0(\tau)}{\phi_0(\tau)}\Bigr) d\tau\biggr\}, \phh c_2\equiv \exp\biggl\{\ilpp a_{21}(\tau)\Bigl(\frac{\phi(\tau)}{\psi(\tau)} - \frac{\phi_1(\tau)}{\psi_1(\tau)}\Bigr) d\tau\biggr\}.
$$
{\bf Remark 3.6}. {\it  For  $I_{a_{12},B}^+(t_0) = +\infty$   (for  $I_{a_{21},-B}^+(t_0) = +\infty$)  the estimate  (3.47)    (the estimate (3.48)) is sharper than the estimate from above from (3.44)      (from  (3.45)).}

Using the estimates (3.47) and  (3.48) for the two particular linearly independent solution  $(\phi_j(t), \psi_j(t)), \ph j =1,2$ with  $\phi_j(t_0) > 0, \ph \psi_j(t_0) > 0, \ph j=1,2$  it is easy to get the following estimates for arbitrary solution  $(\phi(t), \psi(t))$ of the system (1.1) for the case  $a_{12}(t) \ge 0, \ph a_{21}(t) \ge 0, \ph t\ge t_0$:

$$
|\phi(t)| \le C_1(\phi;\psi) \exp\Bigr\{ \il{t_0}{t}[a_{11}(\tau) + a_{12}(\tau) I_{B, a_{21}}^-(t_0;\tau)]d\tau\Bigr\}, \phantom{aaaaaaaaaaaaaaaaaaaa} \eqno (3.49)
$$
$$
\phantom{aaaaaa}|\psi(t)| \le  C_2(\phi;\psi) \exp\Bigr\{\il{t_0}{t}[a_{22}(\tau) + a_{21}(\tau) I_{-B, a_{12}}^-(t_0;\tau)]d\tau\Bigr\}, \phh t\ge t_0, \eqno (3.50)
$$
where $C_j(\phi;\psi), \ph j=1,2$  are some constants depending on $(\phi(t), \psi(t))$.

Example 3.8. Consider the system
$$
\sist{\phi' = \lambda t^\alpha \phi + \mu t ^\beta  \psi;}{\psi' = \nu t ^\gamma \phi + \lambda t^\alpha \psi, \ph t\ge t_0> 0,} \eqno (3.51)
$$
where  $\lambda \in R, \ph \mu > 0, \ph \nu > 0, \ph \alpha > -1, \ph \beta > -1,  \ph  \gamma < -1$                           are some constants and $\beta \ge \alpha, \ph \beta + \gamma + 2 > 0.$ By  (3.49) and  (3.50) for an arbitrary solution  $(\phi(t), \psi(t))$ of the system (1.1) the estimates
$$
|\phi(t)| \le C_1 \exp\biggl\{ \frac{\lambda}{\alpha +1} t ^{\alpha +1} + \frac{\mu\nu t_0^{\gamma+1}}{(\beta + 1)|\gamma + 1|} t ^{\beta + 1}\biggr\}, \eqno (3.52)
$$
$$
|\psi(t)| \le C_2\exp\biggl\{ \frac{\lambda}{\alpha +1} t ^{\alpha +1} + \frac{\mu\nu}{(\beta + 1)(\beta + \gamma + 2)} t ^{\beta + \gamma +  2}\biggr\}, \phh t\ge t_0. \eqno (3.53)
$$
are valid. The estimate of Wazevski (Theorem 4.6.2) for this system is the following:
$$
\sqrt{\phi^2(t) + \psi^2(t)} \le M \exp\biggl\{\frac{\lambda}{\alpha +1} t^{\alpha +1} + \frac{\mu}{2(\beta +1)} t ^{\beta +1}\biggr\}, \phh t\ge t_0,  \eqno (3.54)
$$
where $M=const > 0$. Let as also write for  the  solution $(\phi(t), \psi(t))$  of the   system (3.51) the estimates of Lyapunov [1], p. 132), Yu. S. Bogdanov ([1], p. 133, Theorem  4.61) and S. M. Lozinski ([1], pp. 135 - 137, Theorem 4.6.3).  These estimates are:
$$
\sqrt{\phi^2(t) + \psi^2(t)} \le M_1 \exp\biggl\{\il{t_0}{t}\Bigl[\max\Bigl(\sqrt{\lambda^2\tau^{2\alpha} + \nu ^2 \tau^{2\gamma}}, \sqrt{\lambda^2\tau^{2\alpha} + \mu ^2\tau^{2\beta}}\Bigr)\Bigr] d\tau\biggr\}, \eqno (3.55)
$$
$$
\sqrt{\phi^2(t) + \psi^2(t)} \le M_2\exp\biggl\{\frac{2|\lambda|}{\alpha +1} t^{\alpha +1} + \frac{\mu}{\beta +1} t^{\beta +1}\biggr\}, \eqno (3.56)
$$
$$
\sqrt{\phi^2(t) + \psi^2(t)} \le M_3 \exp\biggl\{\lambda t^{\alpha +1} + \frac{\mu}{\beta +1} t^{\beta +1}\biggr\}, \phh t\ge t_0. \eqno (3.57)
$$
Obviously among the estimates (3.54)  - (3.57) the most sharp is the estimate of Wazevski  (3.54).  For  $\beta + \gamma + 1 > \alpha$  the estimate (3.53) is sharper than (3.54),  and if  $\beta + \gamma + 1 = \alpha$,
then the estimate (3.53) will be sharper of estimate (3.54), provided $\nu < \frac{\alpha +1}{2}.$  For   $\beta \ge \alpha$  and   $t_0^{\gamma +1} < \frac{|\gamma + 1|}{2\nu}$ the estimate (3.52) is sharper than  (3.54). In particular for $\alpha = \beta, \ph \lambda < 0$ on the basis (3.52) and  (3.53) we conclude (taking $t_0 > 0$ enough large),  that the system (3.51) is asymptotically stable, meanwhile under these restrictions and for  $\mu + 2\lambda > 0$ by using (3.54) even impossible establish the Lyapunov stability  of the system  (3.51).
Note that for $\alpha > 1$  or $\beta > 1$  the estimate by freezing method (see [1], p. 141, Theorem 4.6.4)   is not applicable to the system  (3.51).

Let  $a_{12}(t) > 0, \ph t\ge t_0$, and let  $(\phi_0(t), \psi_0(t))$ be a solution of the system (1.1) with $\phi_0(t_0 )=~1,\ph \psi_0(t_0) =~0$.  Then by (3.2)
we have
$$
\phi_0(t) = \exp\biggl\{\il{t_0}{t} \Bigl[a_{11}(\tau) + a_{12}(\tau) x_0(\tau)\Bigr]d\tau\biggr\}, \phh t\ge t_0,
$$
where  $x_0(t)$ is a solution of Eq. (3.1) with  $x_0(t_0) =0$. On the strength of (2.48) from here it follows
$$
0 < \phi_0(t) \le \exp\biggl\{ \frac{1}{2}\il{t_0}{t} S(\tau) d\tau + \frac{1}{2}\sqrt{\il{t_0}{t}a_{12}(\tau) d\tau\il{t_0}{t} \frac{B^2(\tau) + 4a_{12}(\tau)a_{21}(\tau)}{a_{12}(\tau)} d\tau}\biggr\},  \eqno (3.58)
$$
$t\ge t_0.$
Let $(\phi_1(t), \psi_1(t))$ be the solution of the system (1.1)  with   $\phi_1(t_0) =1, \ph \psi_1(t_0)= 1$.           Then by (3.2) we have:
$$
\phi_1(t) = \exp\biggl\{\il{t_0}{t} \Bigl[a_{11}(\tau) + a_{12}(\tau) x_1(\tau)\Bigr]d\tau\biggr\}, \phh t\ge t_0, \eqno (3.59)
$$
where  $x_1(t)$  is a solution of Eq. (3.2) with  $x_1(t_0) = 1$. Since by virtue of Lemma  2.1 and Theorem  2.3.II$^\circ$ the solutions $x_0(t)$ and $x_1(t)$ are $t_0$-normal, taking into account (2.9) from (3.55) and  (3.56) we  get:
$$
0 < \phi_1(t) \le N\exp\biggl\{ \frac{1}{2}\il{t_0}{t} S(\tau) d\tau + \frac{1}{2}\sqrt{\il{t_0}{t}a_{12}(\tau) d\tau\il{t_0}{t} \frac{B^2(\tau) + 4a_{12}(\tau)a_{21}(\tau)}{a_{12}(\tau)} d\tau}\biggr\},  \eqno (3.60)
$$
$t\ge t_0$ where  $N\equiv \exp\biggl\{\ilpp a_{12}(\tau)\Bigr(x_1(\tau) - x_0(\tau)\Bigr)d\tau\biggr\} < +\infty$.  Let $(\phi(t), \psi(t))$ be an arbitrary solution of the system (1.1).  Then since obviously  $(\phi_0(t), \psi_0(t))$ and  $(\phi_1(t), \psi_1(t))$ are linearly independent from the equality  $\phi(t) = \nu_0\phi_0(t) + \nu_1 \phi_1(t), \ph t\ge t_0$, where $\nu_0$  and  $\nu_1$ are some constants, and from  (3.58) and  (3.60) we have:
$$
|\phi(t)| \le M\exp\biggl\{ \frac{1}{2}\il{t_0}{t} S(\tau) d\tau + \frac{1}{2}\sqrt{\il{t_0}{t}a_{12}(\tau) d\tau\il{t_0}{t} \frac{B^2(\tau) + 4a_{12}(\tau)a_{21}(\tau)}{a_{12}(\tau)} d\tau}\biggr\},  \eqno (3.61)
$$
$t\ge t_0$, where  $M=const > 0$. By analogy if $a_{21}(t) > 0, \ph t  \ge t_0$, then
$$
|\psi(t)| \le M_1\exp\biggl\{ \frac{1}{2}\il{t_0}{t} S(\tau) d\tau + \frac{1}{2}\sqrt{\il{t_0}{t}a_{12}(\tau) d\tau\il{t_0}{t} \frac{B^2(\tau) + 4a_{12}(\tau)a_{21}(\tau)}{a_{21}(\tau)} d\tau}\biggr\},  \eqno (3.62)
$$
$t\ge t_0$, where  $M_1=const > 0$.
For the system (3.51) the estimates (3.61) and (3.62) give:
$$
|\phi(t)| \le m \exp\biggl\{\frac{\lambda t^{\alpha + 1}}{\alpha +1} + \sqrt{\frac{\lambda \mu t^{\beta+1}t_0^{\gamma +1}}{(\beta + 1)|\gamma+1|}}\biggr\}, \phantom{aaaaaaaaaaaaaaaaaaaaaaaaaaaaaaaaaaaaaa} \eqno (3.63)
$$
$$
\phantom{aaaaaaaaaaaaaaaaaa}|\psi(t)| \le m_1 \exp\biggl\{\frac{\lambda t^{\alpha + 1}}{\alpha +1} + \sqrt{\frac{\lambda \mu t^{\beta+1}t_0^{\gamma +1}}{(\beta + 1)|\gamma+1|}}\biggr\}, \phh t\ge t_0. \eqno (3.64)
$$
For $0 < \alpha \le \beta$  the estimate (3.63) is sharper of estimate (3.52),  and for  $0 < \alpha \le \beta + \gamma + 2, \ph \beta + 2\gamma + 3 > 0$  the estimate (3.64) is sharper of (3.53).

{\bf 3.4.2. The case when  $a_{12}(t) \ge 0, \ph a_{21}(t) \le 0, \ph t\ge t_0$}. In this subsection we shall assume  $a_{12}(t) \ge 0, \ph a_{21}(t) \le 0, \ph t\ge t_0$.  Then from (3.32) it follows that
$$
\theta(t) \ge 0, \phh t\ge t_0. \eqno (3.65)
$$
Then by (3.27) the radius vector  $\overrightarrow{R}(t) \equiv (J_{-a_{11}}(t)\phi(t), J_{-a_{22}}(t)\psi(t)), \ph t\ge t_0$
rotates monotonically around of origin of the coordinate  plane  $\phi, \psi$ in the negative direction, where $(\phi(t), \psi(t))$  is a solution to the system (1.1) (see pict. 15).
From (3.65) \linebreak
\phantom{aaaaaaaaaaaaaaaaaaaaaaaaaaa} it follows that  $\theta(t)$ can have finite
or negative \linebreak
\phantom{aaaaaaaaaaaaaaaaaaaaaaaaaaa}
infinite limit. By (3.27) in the first case the\linebreak
\phantom{aaaaaaaaaaaaaaaaaaaaaaaaaa}
\hskip 3pt  system (1.1) is non oscillatory, in the second  case\linebreak
\phantom{aaaaaaaaaaaaaaaaaaaaaaaaaaa}
 it is oscillatory. \\
\phantom{aaaaaaaaaaaaaaaaaaaaaaaaaaaaa}
{\bf Theorem 3.10} {\it Let  $a_{12}(t) \ge0, \ph a_{21}(t) \le0, \ph t\ge t_0$,\linebreak
\phantom{aaaaaaaaaaaaaaaaaaaaaaaaaa}
and let the system (1.1) is not oscillatory. Then if:

\vskip -16pt
\begin{picture}(1,1)
\put(50,10){\vector(0,1){70}}
\put(10,40){\vector(1,0){90}}
\put(90,45){$_{\phi}$}
\put(54,70){$_{\psi}$}
\put(0,10){$_{Pict.15.}$}
\put(85,70){$_{\overrightarrow{R}(t)}$}
\put(84,67){\thicklines \qbezier[6](0,0)(-10,-5)(-20,-10)}
\put(66,57.7){\vector(-2,-1){1}}
\put(70,40){\qbezier(0,0)(-3,-14)(-20,-15)\qbezier(-20,-15)(-37,-14)(-40,0)
\qbezier(-40,0)(-37,16)(-20,18)\qbezier(-20,18)(7,19)(10,0)\qbezier(10,0)(7,-25)(-20,-21)
\qbezier(-20,-21)(-34,-18)(-31,-1)\qbezier(-31,-1)(-27,17)(-5,12)}
\put(60,26.5){\vector(-1,-1){1}}
\put(31,39){\vector(0,1){1}}
\put(74.5,26){\vector(-1,-1){1}}
\put(74,52){\vector(1,-1){1}}
\put(44,47.5){\vector(1,1){1}}
\end{picture}

\noindent
 A$_2) \ph I_{a_{12},B}^+(t_0) = +\infty$, then there exists up to arbitrary multiplier the unique linearly independent solutions  $(\phi_*(t), \psi_0(t))$ and  $(\phi_0(t), \psi_*(t))$  of the system (1.1) such that $\phi_0(t) > 0, \ph \psi_*(t) > 0, \ph \phi_0(t), \ph \phi_0(t) < 0, \ph \psi_*(t) > 0, \ph t \ge T$,  for some $T\ge t_0$;        for arbitrary solution $(\phi(t), \psi(t))$  of the system (1.1), linearly independent  of $(\phi_*(t), \psi_0(t))$ \linebreak (of  $(\phi_0(t), \psi_*(t))$) the equality
$$
\liml{t\to+\infty}\frac{\phi_*(t)}{\phi(t)} = 0, \phh \Bigl(\liml{t\to+\infty}\frac{\psi_*(t)}{\psi(t)} = 0,\Bigr), \eqno (3.66)
$$
is fulfilled and
$$
\ilp{T}\frac{a_{12}(\tau) J_S(\tau)}{\phi_*^2(\tau)}d\tau = \ilp{T}\frac{|a_{21}(\tau)| J_S(\tau)}{\psi_*^2(\tau)}d\tau = +\infty, \eqno (3.67)
$$
$$
\ilp{T_1}\frac{a_{12}(\tau) J_S(\tau)}{\phi^2(\tau)}d\tau< +\infty \phh \biggl(\ilp{T_1}\frac{a_{12}(\tau) J_S(\tau)}{\phi^2(\tau)}d\tau< +\infty\biggr); \eqno (3.68)
$$
where $T_1= T_1(\phi;\psi) \ge t_)$ such that   $\phi(t) \ne 0, \ph \psi(t) \ne 0, \ph t\ge T_1$              and \linebreak $\ilp{T_1} a_{12}(\tau)\frac{\psi(\tau)}{\phi(\tau)}d\tau =+\infty$. Moreover if in addition   $\ilpp a_{21}(\tau) I_{a_{12}, -B}^-(t_0;\tau) d\tau = - \infty$,  then $\ilp{T}a_{12}(\tau)\frac{\psi_0(\tau)}{\phi_*(\tau)}d\tau = +\infty;$

\noindent
B$_2) \ph I_{a_{21}, - B}^+(t_0) = - \infty,$  there exists up to arbitrary multiplier the unique solutions \linebreak  $(\phi_{**}(t), \psi_{00}(t))$ and   $(\phi_{00}(t), \psi_{**}(t))$
of the system (1.1) such that  $\phi_{**}(t) > 0, \linebreak \psi_{00}(t) < 0, \ph  \phi_{00}(t) > 0, \ph  \psi_{**}(t)< 0, \ph t\ge T,$  for some  $T\ge t_0$;  for every solution   $(\phi(t), \psi(t))$ of the system (1.1), linearly independent of $(\phi_{**}(t), \psi_{00}(t))$ (of   $(\phi_{00}(t), \psi_{**}(t))$)             the equality
$$
\liml{t\to+\infty}\frac{\phi_{**}(t)}{\phi(t)} = 0, \phh \Bigl(\liml{t\to+\infty}\frac{\psi_{**}(t)}{\psi(t)} = 0,\Bigr), \eqno (3.69)
$$
is fulfilled and
$$
\ilp{T}\frac{a_{12}(\tau) J_S(\tau)}{\phi_{**}^2(\tau)}d\tau = \ilp{T}\frac{|a_{21}(\tau)| J_S(\tau)}{\psi_{**}^2(\tau)}d\tau = +\infty, \eqno (3.70)
$$
$$
\ilp{T_1}\frac{a_{12}(\tau) J_S(\tau)}{\phi^2(\tau)}d\tau< +\infty \phh \biggl(\ilp{T_1}\frac{a_{12}(\tau) J_S(\tau)}{\phi^2(\tau)}d\tau< +\infty\biggr); \eqno (3.71)
$$

where $T_1= T_1(\phi;\psi) \ge t_)$ such that  $\phi(t) \ne 0, \ph \psi(t) \ne 0, \ph t\ge T_1$;
$$
\ilp{T} a_{12}(\tau)\frac{\psi_{00}(\tau)}{\phi_{**}(\tau)}d\tau =-\infty \eqno (3.72).
$$
Moreover if in addition $\ilpp a_{12}(\tau) I_{B,a_{21},}^-(t_0;\tau) d\tau = - \infty$,                                then
$$
\ilp{T_1}a_{12}(\tau)\frac{\psi(\tau)}{\phi(\tau)}d\tau = -\infty. \eqno (3.73)
$$
}

Proof. Since the system (1.1) is not oscillatory from  (3.6) it follows that it has a $T$-regular solution for some $T\ge t_0$.  Then by (3.2) Eq. (3.1) has a $T$-regular solution. Due to Lemma 2.1 from here it follows that Eq. (3.1) has the unique  $T$-extremal solution $x_*(t)$.  Set
$$
\phi_{**}(t) \equiv\exp\biggl\{\il{T}{t}\Bigl[a_{12}(\tau) x_*(\tau) + a_{11}(\tau)\Bigr]d\tau\biggr\}, \phh \psi_{00}(t) \equiv x_*(t)\phi_{**}(t), \phh t\ge T. \eqno (3.74)
$$
By  (3.2) $(\phi_{**}(t).\psi_{00}(t))$ is a solution to the system (1.1) on $[T;+\infty)$, which can be continued on  $[t_0;+\infty)$ as a solution of the system (1.1).  Let us prove  A$_2)$. Suppose  $I_{a_{12},B}^+(t_0) =~ +\infty$. Then by virtue of Theorem 2.4.I$^*$ we have $x_*(t) > 0, \ph t\ge T.$ From here and from (2.74) it follows that  $\phi_{**}(t)> 0, \ph  \psi_{00}(t) > 0, \ph t\ge T$.
Obviously $v_*(t) \equiv - \frac{1}{x_*(t)}$ is the $T$-extremal solution of the equation
$$
z' - a_{21}(t) z^2 - B(t) z + a_{12}(t)=0, \phh t\ge t_0. \eqno (3.75)
$$
Set
$$
\ph_{00}(t) \equiv v_*(t)\psi_{**}(t), \phh \psi_{**}(t) \equiv \exp\biggl\{\il{T}{t}\Bigl[a_{21}(\tau) v_*(\tau) + a_{22}(\tau)\Bigr]d\tau\biggr\}, \phh t\ge T. \eqno (3.76)
$$
By (3.2) $(\phi_{00}(t), \psi_{**}(t))$ is a solution of the system (1.1)  on $[T;+\infty)$, which can be continued on $[t_0;+\infty)$ as a solution of the system  (1.1).  From (3.76) it follows that  $\phi_{00}(t) < 0,\ph \psi_{**}(t) >~0, \ph t\ge T$. Obviously  $(\phi_{**}(t), \psi_{00}(t))$  and    $(\phi_{00}(t), \psi_{**}(t))$  are linearly independent. Since  $a_{12}(t) \ge~0, \ph - a_{21}(t) \ge 0, \ph t\ge t_0$, by virtue of Theorem 2.1 the equations  (3.1) and  (3.22) have the $T$-normal solutions  $x_0(t)$  and  $v_0(t)$ respectively,  for which the integrals $\widetilde{\nu}_{x_0}(t)$   and   $\widetilde{\widetilde{\nu}}_{v_0}(t) \equiv \ilp{t} a_{21}(\tau)\exp\biggl\{- \il{t}{\tau}\Bigl[-2 a_{21}(\xi) v_0(\xi) - B(\xi)\Bigr]d\xi\biggr\}$    converge for all $t \ge T$   and   $\widetilde{\nu}_{x_0}(t) > 0, \ph \widetilde{\widetilde{\nu}}_{v_0}(t) <0, \ph t\ge T.$ Set
$$
\widetilde{\phi}_0(t) \equiv \exp\biggl\{\il{T}{t}\Bigl[a_{12}(\tau) x_0(\tau) + a_{11}(\tau)\Bigr]d\tau\biggr\}, \ph \widetilde{\psi}_0(t) \equiv x_0(t)\widetilde{\phi}_0(t), \phh t\ge T;
$$
$$
\widetilde{\widetilde{\phi}}_0(t) \equiv v_0(t) \widetilde{\widetilde{\phi}}_0(t), \phh \widetilde{\widetilde{\psi}}_0(t) \equiv \exp\biggl\{\il{T}{t}\Bigl[-a_{21}(\tau) v_0(\tau) + a_{22}(\tau)\Bigr]d\tau\biggr\}, \phh t\ge T.
$$
By (3.2) $(\widetilde{\phi}_0(t), \widetilde{\psi}_0(t))$  and  $(\widetilde{\widetilde{\phi}}_0(t), \widetilde{\widetilde{\phi}}_0(t))$ are solutions of the system (1.1) on  $[T;+\infty)$. \linebreak By (3.2) we have
$$
\ilp{t}\frac{a_{12}(\tau) J_S(T;\tau)}{\widetilde{\widetilde{\phi}}_0(\tau)^2} d\tau = \widetilde{\nu}_{x_0}(t) > 0, \ph \ilp{t}\frac{a_{21}(\tau) J_S(T;\tau)}{\widetilde{\widetilde{\psi}}_0(\tau)^2} d\tau = \widetilde{\widetilde{\nu}}_{-v_0}(t)< 0, \ph t\ge t_0.
$$
Then the equalities  (3.66) - (3.68) follow from Theorem 3.5 and from the equality \linebreak $\ilp{T_1}a_{12}(\tau)\frac{\psi(\tau)}{\phi(\tau)} d\tau =+\infty$ (since  $\frac{\psi(t)}{\phi(t)}$ is the  $T_1$-extremal solution of Eq. (3.1))  and \linebreak  $\ilp{T}a_{12}(\tau)\frac{\psi_*(\tau)}{\phi_*(\tau)} d\tau =+\infty$ (since  $\frac{\psi_*(t)}{\phi_*(t)}= x_*(t)$)  from Theorem 2.4.I$^*$. The assertion  A$_2)$  is proved. Let us prove B$_2)$.   Suppose  $I_{a_{21},-B}^+(t_0) = -\infty$.  Then by virtue of Theorem 2.4.II$^*$ the inequality $x_*(t) < 0,\ph t \ge T_1$,  is fulfilled for some $T_1 \ge t_0$. Without loss generality we shall take that  $T_1 =T$.  Show that $v_*(t)\equiv - \frac{1}{x_*(t)}$ is the  $T$-extremal solution of Eq. (3.75). Suppose that it is not so. Let then  $\widetilde{v}_*(t)$ be the  $T$-extremal solution of Eq. (3.75) (existence of $\widetilde{v}_*(t)$   follows from Lemma 2.1). By  Lemma 2.1 we have $\widetilde{v}_*(T) < v_*(T)$. Then   $\widetilde{x}_*(t) \equiv - \frac{1}{\widetilde{v}_*(t)}$ is a solution of eq. (3.1) with  $\widetilde{x}_*(T) < x_*(T)$. Therefore by Lemma 2.1 the solution $\widetilde{x}_*(t)$  is not  $T$-regular and has in an point  $T_2 > T$ a vertical asymptote. From here it follows that  $\widetilde{v}_*(t) =0, \ph t\in [T_2;T_3]$ and  $\widetilde{v}_*(t) \ne 0, \ph t\in (T_3;T_4)$ for some $T_4>T_3\ge T_2$. Since  $x_*(t) < 0, \ph t\ge T$, we have  $\widetilde{x}_*(t) < 0, \ph t\in [T;T_1)$. From here it follows that
$$
\widetilde{v}_*(t) > 0, \phh t\in [T;T_1). \eqno (3.77)
$$
Since by (3.2) $\widetilde{v}_*(t) = \frac{\widetilde{\phi}_*(t)}{\widetilde{\psi}_*(t)}, \ph t\ge T$,              where  $(\widetilde{\phi}_*(t))(\widetilde{\psi}_*(t))$ is a solution of the system  (1.1) with  $\widetilde{\psi}_*(t) \ne 0, \ph t\ge T$, on the strength of (3.27) from (3.24) it follows that  $\widetilde{v}_*(T_5) < 0$  for some $T_5\in(T_3;T_4)$. Then
$$
\widetilde{x}_*(T_5) = - \frac{1}{\widetilde{v}_*(T_5)} > 0, \eqno (3.78)
$$
and since $x_*(T_5) < 0$, we have that  $\widetilde{x}_*(t)$ is a $T_5$-normal solution of Eq.  (3.1). It is easy to show that from the equality  $I_{a_{21}, -B}(t_0) = -\infty$ follows   $I_{a_{21}, -B}(T_5) = -\infty$. By virtue of Theorem 2.4.II$^*$ from here and from (3.78) it follows that  $\widetilde{x}_*(T_6) = 0$ for some  $T_6 > T_5$.
Hence  $\widetilde{v}_*(t)$ has a vertical asymptote in  $T_6$,  i. e.  $\widetilde{v}_*(t)$  is not     $T$-regular. We came to the contradiction which proves that   $v_*(t) \equiv - \frac{1}{x_*(t)}$ is the       $T$-extremal solution of Eq. (3.75). Consider the functions
$$
\phi_{**}(t)\equiv\exp\biggl\{\il{T}{t}\Bigl[a_{12}(\tau)x_*(\tau) + a_{11}(\tau)\Bigr]d\tau\biggr\}, \phh \psi_{00}(t) \equiv x_*(t) \phi_{**}(t), \eqno (3.79)
$$
$$
\psi_{**}(t)\equiv\exp\biggl\{\il{T}{t}\Bigl[-a_{21}(\tau)x_*(\tau) + a_{22}(\tau)\Bigr]d\tau\biggr\}, \phh \phi_{00}(t) \equiv v_*(t) \psi_{**}(t), \phh t\ge T.
$$
By (3.2) $(\phi_{**}(t), \psi_{00}(t))$  and  $(\phi_{00}(t), \psi_{**}(t))$ are solutions of the system  (1.1) on  $[T;+\infty)$, which can be continued on  $[t_0;+\infty)$ as solutions of the system  (1.1).
Obviously $\phi_{**}(t) > 0, \ph \psi_{00}(t) < 0, \ph \phi_{00}(t) > 0, \ph \psi_{**}(t) > 0, \ph t\ge T$.  Therefore  $(\phi_{**}(t), \psi_{00}(t))$ and \linebreak $(\phi_{00}(t), \psi_{**}(t))$ are linearly independent.
Further as in the proof of  A$_2)$ it shows that for some solutions  $(\phi_0(t), \psi_0(t))$
and  $(\phi_1(t), \psi_1(t))$  of the system (1.1)  with  $\phi_0(t) \ne 0, \ph \psi_1(t) \ne 0, \ph t\ge T$, the integrals  $\ilp{t}\frac{a_{12}(\tau) J_S(\tau)}{\phi_0^2(\tau)} d\tau$  and  $\ilp{t}\frac{a_{21}(\tau) J_S(\tau)}{\psi_0^2(\tau)} d\tau$  converge for all  $t \ge T$ and do not vanish on   $[T;+\infty)$. Then the equalities (3.69) - (3.71) are consequences of \linebreak Theorem~3.5, and the equalities (3.72) and  (3.73) follow from Theorem 2.4.II$^*$ and from (3.2) ((3 .79)). The theorem is proved.

{\bf Theorem 3.11}. {\it Suppose   $a_{12}(t) \ge 0, \ph a_{21}(t) \le 0, \ph t\ge t_0, \ph I_{a_{12},B}^+(t_0) < +\infty, \linebreak I_{-a_{21},-B}^+(t_0) < +\infty$.
Then the following assertions are valid:

\noindent
A$_3$) the system (1.1) is non oscillatory;

\noindent
B$_3$) there exist up to arbitrary multiplier the unique solutions  $(\phi_*(t), \psi_0(t))$
and $(\phi_0(t), \psi_*(t))$ of the system (1.1) such that  $\phi_*(t) > 0, \ph \psi_0(t) < 0, \ph \phi_0(t) < 0, \ph \psi_*(t) > 0, \ph t\ge T$ for some  $T > t_0$,  for each solution $(\phi(t), \psi(t))$ of the system (1.1), linearly independent of  $(\phi_*(t), \psi_0(t))$   (of   $(\phi_0(t), \psi_*(t))$) the equality
$$
\liml{t\to+\infty}\frac{\phi_*(t)}{\phi(t)} = 0, \phh \biggl(\liml{t\to+\infty}\frac{\psi_*(t)}{\psi(t)} = 0\biggr). \eqno (3.80)
$$
is fulfilled and
$$
\ilp{T}\frac{a_{12}(\tau) J_S(\tau)}{\phi_*^2(\tau)} d\tau = - \ilp{T}\frac{a_{21}(\tau) J_S(\tau)}{\psi_*^2(\tau)} = +\infty, \eqno (3.81)
$$
$$
\ilp{T_1}\frac{a_{12}(\tau) J_S(\tau)}{\phi^2(\tau)} d\tau < +\infty \phh \biggl(\ilp{T_1}\frac{|a_{21}(\tau)| J_S(\tau)}{\psi^2(\tau)} < +\infty\biggr). \eqno (3.82)
$$
For every solution  $(\phi(t), \psi(t))$ of the system (1.1), linearly independent of  $(\phi_*(t), \psi_0(t))$  the integral  $\ilp{T_1}a_{12}(\tau)\frac{\psi(\tau)}{\phi(\tau)}f\tau \ph \Bigl(\ilp{T_1}a_{21}(\tau)\frac{\phi(\tau)}{\psi(\tau)}d\tau\Bigr)$ converges,
where $T_1 \ge t_0$ such that \linebreak $\phi(t) \ne~0, \ph \psi(t) \ne 0, \ph t\ge T_1$, аnd
$$
\ilp{T}a_{12}(\tau) \frac{\psi_0(\tau)}{\phi_*(\tau)} d\tau = -\infty, \phh \ilp{T}a_{21}(\tau) \frac{\phi_0(\tau)}{\psi_*(\tau)} d\tau = +\infty. \eqno (3.83)
$$
}

Proof. The assertion  A$_3)$ immediately follows from (3.32) and from the equalities
$$
\ilpp A_{12}(\tau) d\tau = I_{a_{12},B}^+(t_0), \phh \ilpp A_{21}(\tau) d\tau = I_{a_{21},-B}^+(t_0).
$$
Let us prove  B$_3)$.  From A$_3)$  and from (3.6) it follows that for some solution  $(\phi_0(t), \psi_0(t))$  of the system (1.1) the inequality $\phi_0(t) \ne 0, \ph t\ge T$, is fulfilled for some  $T \ge t_0$. Therefore $x_0(t) \equiv \frac{\psi_0(t)}{\phi_0(t)}$ is a $T$-regular solution of eq. (3.1). Then since  $a_{12}(t) \ge 0, \ph t\ge t_0$, by virtue of Lemma 2.1 Eq. (3.1) has  the $T$-extremal solution  $x_*(t)$.
By virtue of \linebreak Theorem 2.4.IV$^*$  from the conditions of the theorem it follows that  $x_*(t) < 0, \ph t\ge T_1$,
for some $T_1 \ge t_0$. Without loss of generality we shall take that  $T_1 = T$. Then  $v_0(t) \equiv - \frac{1}{x_*(t)}$ is a $T$-regular solution of Eq. (3.75).  Since  $a_{21}(t)$ does not change sign from here it follows that Eq. (3.75) has  the $T$-extremal solution $v_*(t)$. According to  Theorem  2.4.III$^*$     from the convergence of the integral $I_{a_{21}, -B}^+(t_0)$ it follows that  $v_*(t) < 0, \ph t\ge T_2,$  for some $T_2 \ge t_0$. Without loss of generality we shall take that  $T_2 = T$. By (3.2) the functions
$$
\phi_*(t)\equiv \exp\biggl\{\il{T}{t}\Bigl[a_{12}(\tau)x_*(\tau) + a_{11}(\tau)\Bigr]d\tau\biggr\} \phh \psi_0(t) \equiv x_*(t) \phi_*(t),
$$
$$
\psi_*(t)\equiv \exp\biggl\{\il{T}{t}\Bigl[-a_{21}(\tau)v_*(\tau) + a_{22}(\tau)\Bigr]d\tau\biggr\} \phh \phi_0(t) \equiv v_*(t) \psi_*(t), \phh t\ge T,
$$
form solutions  $(\phi_*(t),\psi_0(t)), \ph  (\phi_*(t),\psi_0(t))$ of the system  (1.1). Obviously   $\phi_*(t) > 0, \linebreak \psi_0(t) < 0, \ph \phi_0(t) < 0, \ph \psi_*(t) > 0, \ph t\ge T$. Prove that these solutions are linearly indepen-\linebreak dent. Suppose that it is not so.  Then  $\phi_0(t) = \lambda\phi_*(t), \ph \psi_*(t) = \lambda\psi_0(t), \ph t\ge T,$  for some $\lambda = const$. From here it follows  $\frac{\psi_0(t)}{\psi_*(t)} = \frac{\psi_*(t)}{\psi_0(t)}$, i. e.  $x_*(t) = \frac{1}{v_*(t)}, \ph t\ge T$.
But on the other hand  $x_*(t) = - \frac{1}{v_0(t)}, \ph t\ge T$. From here it follows that  $v_*(t) = - v_0(t), \ph t\ge T$. Then  $x_0(t)\equiv - \frac{1}{v_*(t)} = \frac{1}{v_0(t)}$ is a positive    $T$-regular solution of Eq. (3.1): $x_*(t) = - x_0(t)$. Since  $x_0(t) > 0,$ and  $x_*(t) < 0, \ph t \ge T$, by virtue of Lemma  2.1 we have that
$x_0(t)$ is $T$-normal. Then on the strength of Theorem 2.4 the integral  $\ilp{T}a_{12}(\tau)x_0(\tau) d\tau$  converges.  But on the other hand by virtue of the same Theorem 2.4.III$^* \ph \ilp{T}a_{12}(\tau)x_0(\tau) d\tau= - \ilp{T}a_{12}(\tau)x_*(\tau) d\tau = +\infty$. We came to the contradiction. Hence  $(\phi_*(t),\psi_0(t))$ and  $(\phi_*(t),\psi_0(t))$ are linearly independent.
Let $x_0(t)$  and  $v_0(t)$ be $ T$-normal solutions of the equations (3.1)  and  (3.75) respectively. Then by virtue of Theorem 2.2  the integrals $\nu_{x_0}(t)$ and  $\widetilde{\widetilde{\nu}}_{v_0}(t)$
converge for all  $t\ge T$  and
$$
\nu_{x_0}(t) > 0, \phh \widetilde{\widetilde{\nu}}_{v_0}(t)< 0, \phh t\ge T. \eqno (3.85)
$$
Set
$$
\phi_0(t)\equiv \exp\biggl\{\il{T}{t}\Bigl[a_{12}(\tau)x_0(\tau) + a_{11}(\tau)\Bigr]d\tau\biggr\} \phh \psi_0(t) \equiv x_0(t) \phi_0(t),
$$
$$
\psi_1(t)\equiv \exp\biggl\{\il{T}{t}\Bigl[-a_{21}(\tau)v_0(\tau) + a_{22}(\tau)\Bigr]d\tau\biggr\} \phh \phi_1(t) \equiv v_0(t) \psi_1(t), \phh t\ge T,
$$
By (3.2) $(\phi_0(t), \psi_0(t)), \ph (\phi_1(t), \psi_1(t))$  are solutions of the system (1.1) on  $[T;+\infty)$, which can be continued on  $[t_0;+\infty)$ as solutions of the system (1.1) By virtue of  (3.2) we have:
$$
I_{\phi_0}(t) \equiv \ilp{t}\frac{a_{12}(\tau)J_S(\tau)}{\phi_0^2(\tau)} d\tau = \nu_{x_0}(t), \phh \widetilde{I}_{\psi_1}(t) \equiv \ilp{t}\frac{a_{21}(\tau)J_S(\tau)}{\psi_1^2(\tau)} d\tau = \widetilde{\widetilde{\nu}}_{v_0}(t) , \ph t\ge T.
$$
From here and from (3.85) it follows that  $I_{\phi_0}(t) \ne 0, \ph  \widetilde{I}_{\psi_1}(t) \ne 0, \ph t\ge T.$  Then the equalities  (3.80), (3.81), (3.82) follow from Theorem 3.5  and the equalities  (3.80), (3.81) and the convergence of the integral  $\ilp{T_1}a_{12}(\tau)\frac{\psi(\tau)}{\phi(\tau)} d\tau \ph \biggl( \ilp{T_1}a_{21}(\tau)\frac{\phi(\tau)}{\psi(\tau)} d\tau \biggr)$
 follow from \linebreak Theorem 2.4.III$^*$ and from  (3.2).  The theorem is proved.

{\bf Theorem 3.12}. {\it Let $a_{12}(t) \ge 0, \ph a_{21}(t) \le 0, \ph t\ge t_0$, and let for some solution $(\phi(t), \psi(t))$  of the system (1.1) and for some  $T \ge t_0$  the inequality $\phi(t) \ne 0, \ph t\ge T$ is fulfilled. Then the system (1.1) satisfies to the non conjugation condition on $[T;+\infty)$.
}

Proof. By (3.6) we have:
$$
\phi(t) = \mu \frac{J_{S/2}(t)}{\sqrt{y_0(t)}} \sin \biggl(\il{t_0}{t} a_{12}(\tau) y_0(\tau) d\tau + \nu\biggr), \phh t\ge t_0,
$$
where  $\mu$   and $\nu$ are some constants. Then  $\il{t_0}{t} a_{12}(\tau)  y_0(\tau)d\tau + \nu  \in (\pi k_0, \pi(k_0 + 1)) \ph t\ge T,$ for some integer $k_0$.  Therefore for arbitrary $\nu_1 \in R$ we have
$$
\il{t_0}{t} a_{12}(\tau)  y_0(\tau)d\tau + \nu_1 \in (\pi k_0 + \nu_1 - \nu; \pi(k_0 + 1) + \nu_1 - \nu), \eqno (3.86)
$$
Let  $(\phi_1(t), \psi_1(t))$ be an arbitrary solution to the system (1.1). Then by  (3.6) we have
$$
\phi_1(t) = \mu_1 \frac{J_{S/2}(t)}{\sqrt{y_0(t)}} \sin \biggl(\il{t_0}{t} a_{12}(\tau) y_0(\tau) d\tau + \nu_1\biggr), \phh t\ge t_0,
$$
and since the function  $f(t) \equiv \il{t_0}{t} a_{12}(\tau) y_0(\tau) d\tau $ monotonically non decreases from (3.86) it follows that  $\phi_1(t)$ vanishes on  $[T;+\infty)$  no more than one time. Let us show that the function  $\psi_1(t)$ vanishes on  $[T;+\infty)$ no more than one time. Suppose that it is not so.  Then for some  $t_k \ge T, \ph k=\overline{1,4}, \ph T < t_1 \le t_2 < t_3 \le t_4$ we have $\psi_1(t) = 0, \ph t\in [t_1;t_2] \cup [t_3;t_4]$   and $\psi_1(t) \ne 0, \ph t\in(t_2;t_3)$.                                                Then since  $a_{12}(t)$  and   $a_{21}(t)$  do not change sign by virtue od Corollary 3.2 from the work  [6] we have  $\phi(t_5) = 0$  for some  $t_5 \in(t_2; t_3)$, which contradicts the condition of the theorem.  Hence $\psi_1(t)$ vanishes on  $[T;+\infty)$  no more than one time. The theorem is proved.

Now we investigate the system (1.1) in the case when it is non oscillatory. Let then  $T \ge t_0$ such that according to Theorem 3.12 the system (1.1) satisfies the non conjugation condition on $[T;+\infty)$, and some its solution   $(\phi_0(t), \psi_0(t))$ is $T$-regular. Without loss of generality we shall take that  $\phi_0(T) =1$. Consider the family of solutions of the system  (1.1): $\Omega_T \equiv \{(\phi(t), \psi(t)) : \ph \phi(t_0) =1\}$. Since $(\phi_0(t), \psi_0(t)) \in \Omega_T$  and  $\phi_0(t) \ne 0, \ph t\ge T$, by    (3.2) we have that $x_0(t) \equiv \frac{\psi_0(t)}{\phi_0(t)}$ is a $T$-regular solution to eq. (3.1). Then according to Lemma 2.1  Eq. (3.1) has the unique  $T$-extremal solution. Suppose   $I_{a_{12},B}^+(t_0) = + \infty$ and   $(\phi(t), \psi(t)) \in \Omega_T, \ph \psi(T) > x_*(T)$. Then by virtue of Theorem  2.4.I$*$ the function  $\Phi_T(t) \equiv J_{-a_{11}}(T;t)\phi(t) \equiv \exp\biggl\{\il{T}{t}a_{12}(\tau) \frac{\psi(\tau)}{\phi(\tau)}\biggr\}$  is monotonically non decreasing  function on $[T;+\infty)$ and if  $\psi(T) > x_*(T)$, then $\Psi_T(+\infty) = +\infty$. Moreover if in addition \linebreak $\ilpp |a_{21}(\tau)| I_{a_{12}, -B}^+(t_0;\tau) = +\infty,$  then  (for $\psi(T) = x_*(T)$)   $\Phi_T(+\infty) = +\infty$. If  $\psi(T) < 0$, then according to Lemma 2.1 the function  $\Phi_T(t)$  vanishes on  $[t_1;t_2], \ph T < t_1\le t_2 < +\infty$, and from (3.27)  it follows (since  $\Phi_T(t)> 0, \ph t\in [T;t_1)$), that  $\Phi_T(t) < 0, \ph t > t_2$  (see pict. 16).
Suppose  $I_{a_{21},-B}^+(t_0) = -\infty$. Then by virtue of Theorem 2.4.II$^*$, if $\psi(T) > x_*(T)$  and  $\psi(T) > 0$, then there exist $t_2 \ge t_1 > T$ such that  $\Phi_T(t)$  monotonically non decreasing on  $[T;t_1)$, is a constant on  $[t_1;t_2]$ and monotonically non increases on  $[t_2;+\infty)$, and  $\Phi_T(t) > 0, \ph t\ge T$.  If $\phi(T) = x_*(T)$, then $\Phi_T(t) > 0, \ph t\ge T$ and   $\Phi_T(+\infty) =0$, and if in addition  $\ilpp a_{12}(\tau) I_{B,|a_{21}|}^-(t_0;\tau) d\tau = +\infty$, then  for  $\psi(T) > x_*(T)$ also we have $\phi_T(+\infty) =0$. If
 $\psi(T) < x_*(T)$,
there exist  $t_2 \ge t_1 > T$ such that  $\Phi_T(t) > 0, \ph t\in [T;t_1), \linebreak \Phi_T(t) =0,\ph t\in [t_1;t_3]$  and  $\Phi_T(t) < 0, \ph t > t_2$ (see pict.  17).

\begin{picture}(100,140)
\put(10,40){\vector(0,1){160}}
\put(-10,90){\vector(1,0){180}}
\put(160,82){$_{t}$}
\put(40,82){$_{T}$}
\put(50,89.5){\circle*{3}}
\put(14,170){$_{\Phi_T(t)}$}
\put(-10,15){$_{Pict.16.\ph \alpha_*= \ph arctan \ph x_*(T)\ph > \ph 0, \ph I_{a_{12},B}^+(t_0) =+\infty.}$}

\put(50,110){\thicklines \qbezier[40](-40,0)(40,0)(120,0)}
\put(50,110){\thicklines \qbezier[40](0,0)(60,20)(120,40)}
\put(50,40){\thicklines \qbezier[60](0,0)(0,80)(0,160)}
\put(77,110){\qbezier (5,10)(10.5,9.8)(11,0)}
\put(5,110){$_1$}
\put(95,115){$_{\alpha_*}$}
\put(50,110){\thicklines \qbezier(0,0)(30,8)(120,60)}
\put(50,110){\thinlines \qbezier(0,0)(30,12)(50,32) \qbezier(50,32)(70,52)(120,75)}
\put(50,110){\thinlines \qbezier(0,0)(30,12)(50,35) \qbezier(50,35)(70,60)(120,85)}
\put(50,110){\thinlines \qbezier(0,0)(30,15)(50,45) \qbezier(50,45)(70,74)(100,85)}
\put(50,110){\thinlines \qbezier(0,0)(40,14)(60,-10) \qbezier(60,-10)(70,-30)(120,-35)}
\put(50,110){\thinlines \qbezier(0,0)(40,13)(60,-17) \qbezier(60,-17)(70,-30)(120,-45)}
\put(50,110){\thinlines \qbezier(0,0)(40,11)(60,-20) \qbezier(60,-20)(70,-35)(120,-55)}

\put(220,40){\vector(0,1){160}}
\put(200,60){\vector(1,0){220}}
\put(410,52){$_{t}$}
\put(250,52){$_{T}$}
\put(260,59.5){\circle*{3}}
\put(290,52){$_{T_1}$}
\put(300,59.5){\circle*{3}}

\put(224,170){$_{\Phi_T(t)}$}
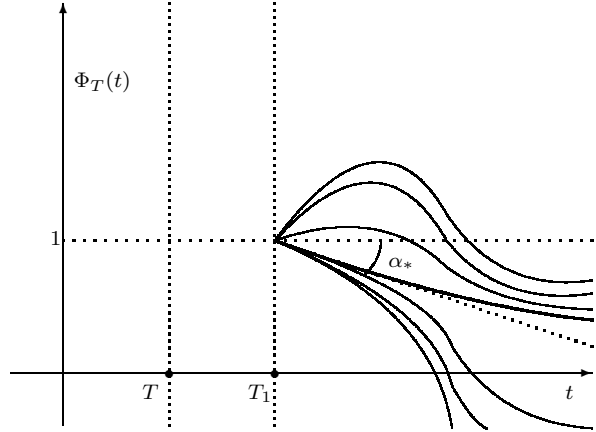
\put(209,15){$_{Pict.17.\ph \alpha_*= \ph arctan \ph x_*(T_1)\ph < \ph 0, \ph I_{a_{21},-B}^+(t_0) =-\infty.}$}

\put(260,110){\thicklines \qbezier[50](-40,0)(60,0)(160,0)}
\put(300,110){\thicklines \qbezier[40](0,0)(60,-20)(120,-40)}

\put(260,40){\thicklines \qbezier[60](0,0)(0,80)(0,160)}
\put(300,40){\thicklines \qbezier[60](0,0)(0,80)(0,160)}

\put(329,87){\qbezier (5,10)(11.5,16.8)(11,23)}
\put(215,110){$_1$}
\put(343,102){$_{\alpha_*}$}
\put(300,110){\thicklines \qbezier(0,0)(70,-25)(120,-30)}
\put(300,110){\thinlines \qbezier(0,0)(40,54)(65,10) \qbezier(65,10)(85,-20)(120,-15)}
\put(300,110){\thinlines \qbezier(0,0)(40,44)(65,0) \qbezier(65,0)(85,-25)(120,-20)}
\put(300,110){\thinlines \qbezier(0,0)(40,14)(65,-10) \qbezier(65,-10)(85,-25)(120,-26)}
\put(300,110){\thinlines \qbezier(0,0)(60,-20)(67,-40) \qbezier(67,-40)(85,-69)(120,-71)}
\put(300,110){\thinlines \qbezier(0,0)(60,-25)(67,-55) \qbezier(67,-55)(75,-69)(80,-71)}
\put(300,110){\thinlines \qbezier(0,0)(60,-25)(67,-71) }
\end{picture}

\noindent
Suppose   $I_{a_{12}, B}^+(t_0) < +\infty,  I_{-a_{21},-B}^+(t_0) < +\infty$.  Then by  Theorem 2.4.III$^*$ there exists \linebreak $t_1 \ge T$ such that
 $x_*(t) < 0, \ph t\ge T_1$.  There  ex-\\ists  $\psi_N^+ > 0$  such that
if $\psi(T_1) \in (0;\phi+N^+)$,\\ then there exist  $t_2 \ge t_1 > T_1$ such that  $\Phi_T(t)$ \\monotonically non decreasing on $[T_1;t_1)$, is \\a constant on  $[t_1;t_2]$  and monotonically non \\decreases for  $t > t_2$, and $\Phi_T(t) > 0, \ph t\ge T_1$. \\If $\psi(T_1) >  \psi_N^+$, then  $\Phi_T(t)$ is monotonically \\non decreasing  and has finite limit,  and for \\$\phi(T_1) = x_*(T_1) \ph \Phi_T(+\infty) = 0$. If  $\psi(T_1) < \\ <x_*(T_1)$,  then   there exist  $t_2\ge t_1 > T_1$ such \\that  $\Phi_T(t) > 0, \ph t\in [T_1;t_1)$, monotonically\\ non increases, $\Phi_T(t) = 0, \ph t\in [t_1;t_2]$ and\\ $\Phi_T(t) < 0, \ph t\ge t_2$   (see pict. 18).

\vskip -15pt

\begin{picture}(1,1)

\put(220,25){\vector(0,1){160}}
\put(200,60){\vector(1,0){220}}
\put(410,52){$_{t}$}
\put(260,52){$_{T}$}
\put(270,59.5){\circle*{3}}
\put(290,52){$_{T_1}$}
\put(300,59.5){\circle*{3}}

\put(224,170){$_{\Phi_T(t)}$}

\put(229,33){$_{Pict.18.}$}
\put(215,25){$_{ \alpha_*= \ph arctan \ph x_*(T_1)\ph < \ph 0,\ph \alpha_+= \ph arctan \ph \psi_N^+(T_1)\ph > \ph 0,}$}
\put(235,15){$_{ I_{a_{12},B}^+(t_0) <+\infty,\ph   I_{-a_{21},-B}^+(t_0) <+\infty.}$}

\put(260,110){\thicklines \qbezier[50](-40,0)(60,0)(160,0)}
\put(300,110){\thicklines \qbezier[40](0,0)(60,-20)(120,-40)}

\put(310,110){\thicklines \qbezier[70](-80,-50)(10,20)(100,70)}

\put(270,40){\thicklines \qbezier[60](0,0)(0,70)(0,140)}
\put(300,40){\thicklines \qbezier[60](0,0)(0,70)(0,140)}

\put(329,87){\qbezier (5,10)(11.5,16.8)(11,23)}
\put(215,110){$_1$}
\put(343,102){$_{\alpha_*}$}
\put(252,67){$_{\alpha_+}$}

\put(237,60){\qbezier (5,10)(10.5,9.8)(11,0)}

\put(300,110){\thicklines \qbezier(0,0)(70,-25)(120,-30)}
\put(300,110){\thinlines \qbezier(0,0)(40,40)(65,46) \qbezier(65,46)(85,53)(120,55)}
\put(300,110){\thinlines \qbezier(0,0)(40,30)(65,36) \qbezier(65,36)(85,43)(120,44)}
\put(300,110){\thinlines \qbezier(0,0)(40,25)(65,27) \qbezier(65,27)(85,31)(120,32)}
\put(300,110){\thicklines \qbezier(0,0)(40,24)(65,10) \qbezier(65,10)(85,-6)(120,-7)}
\put(300,110){\thinlines \qbezier(0,0)(40,14)(65,-10) \qbezier(65,-10)(85,-25)(120,-26)}
\put(300,110){\thinlines \qbezier(0,0)(40,18)(65,0) \qbezier(65,0)(85,-15)(120,-16)}

\put(300,110){\thinlines \qbezier(0,0)(60,-20)(67,-40) \qbezier(67,-40)(85,-69)(120,-71)}
\put(300,110){\thinlines \qbezier(0,0)(60,-25)(67,-55) \qbezier(67,-55)(75,-69)(80,-71)}
\put(300,110){\thinlines \qbezier(0,0)(60,-25)(67,-71) }
\end{picture}
\phantom{aaaaaaaaaaaaaaaaaaaaaaaa} By (2.52) \\for the solution  $(\phi(t), \psi(t)) \in \Omega_T$ with $\psi(T) \ge x_*(T)$  the inequality
$$
0 < \phi(t) \le \exp\Biggl\{\frac{1}{2}\Biggl(\il{T}{t}S(\tau)d\tau +\sqrt{\il{T}{t} a_{12}(\tau)\Bigl[4 \frac{\psi(T)}{\phi(T)} + \il{T}{\tau}\frac{B^2(s) + 4 a_{12}(s)a_{21}(s)}{a_{12}(s)} d s\Bigr] }d\tau\Biggr)\Biggr\},
$$
$t\ge T$ is fulfilled.

\hskip 20 pt

\centerline{\bf References}

\hskip 20 pt

\noindent
1. L. Ya. Adrianova, Vvedenie v teoriu lineinikh sistem differensial'nikh uravnenii (Intro\linebreak \ph duction  to the Theory of Linear Systems of Differential Equations). St. - Peterburg:\linebreak \phantom{a}  Izd. St. - Peterburg. Univ., 1992.

\noindent
2. G. A. Grigorian, On the Stability of Systems of Two First - Order Linear Ordinary\linebreak \phantom{aa} Differential Equations, Differ. Uravn., 2015, vol. 51, no. 3, pp. 283 - 292.

\noindent
3. G. A. Grigorian. Necessary Conditions and a Test for the Stability of a System of Two\linebreak \phantom{aa} Linear Ordinary Differential Equations of the First Order. Differ. Uravn., 2016, \linebreak \phantom{aa} Vol. 52, No. 3, pp. 292 - 300.

\noindent
4. G. A. Grigorian. Some Properties of Solutions of Systems of Two Linear First - Order\linebreak \phantom{aa} Ordinary Differential Equations, Differ. Uravn., 2015, Vol. 51, No. 4, pp. 436 - 444.

\noindent
5. G. A. Grigorian, Global Solvability of Scalar Riccati Equations. Izv. Vissh.\linebreak \phantom{aa} Uchebn. Zaved. Mat., vol. 51, 2015, no. 3, pp. 35 - 48.

\noindent
6.  G. A. Grigorian. Oscillatory Criteria for the Systems of Two First - Order Linear\linebreak \phantom{aa} Ordinary Differential equations. Rocky Mountain Journal of Mathematics, vol. 47,\linebreak \phantom{aa} no. 5, 2017, pp. 1497 - 1524.

\noindent
7. J. D. Mirzov, Asymptotic properties of solutions of nonlinear non autonomus ordinary\linebreak \phantom{aa} differential equations, Brno: Masarik Univ, 2004.

\noindent
8. R. Bellman, Stability theory of differential equations. Moskow, Foreign Literature\linebreak \phantom{aa} Publishers, 1964.

\noindent
9. .A. I. Egorov, Riccati Equations (Fizmatlit, Moscow, 2011) [in Russian].

\noindent
10. G. A. Grigorian.  On Two Comparison Tests for Second-Order Linear  Ordinary\linebreak \phantom{aa} Differential Equations (Russian) Differ. Uravn. 47 (2011), no. 9, 1225 - 1240; trans-\linebreak \phantom{aa} lation in Differ. Equ. 47 (2011), no. 9 1237 - 1252, 34C10.

\noindent
11.  G. A. Grigorian, Properties of solutions of Riccati equation, Journal of Contemporary \\\phantom{aa} Mathematical Analysis, 2007, vol. 42, No 4, pp. 184 - 197.

\noindent
12. G. A. Grigorian, Two Comparison Criteria for Scalar Riccati Equations with\linebreak \phantom{aa} Applications. Russian Mathematics (Iz. VUZ), 56, No. 11, 17 - 30 (2012).

\noindent
13. G. A. Grigorian, Some Properties of Solutions to Second - Order Linear Ordinary\linebreak \phantom{aa} Differential Equations. Trudty Inst. Matem. i Mekh. UrO RAN, 19, No. 1, 69 - 80\linebreak \phantom{aa} (2013).

\noindent
14.  L. Cesary,  Asymptotic  behavior  and stability problems in  ordinary differential equations.\linebreak \phantom{aa} Moskow, ''Mir'', 1964.

\noindent
15. L. A. Gusarov, On Vanishing of Solutions of the Second Order Differential Equations.\linebreak \phantom{aa} Sov. Phys. Dokl. 71 (1) 0 -12 (1950).

\noindent
16.  C. A. Swanson, Comparison and oscillation theory of linear differential equations.\linebreak \phantom{aa}  Academic press. New York and London, 1968.

\noindent
17. Ph. Hartman, Ordinary differential equations, SIAM - Society for industrial and\linebreak \phantom{aaa} applied Mathematics, Classics in Applied Mathematics 38, Philadelphia 2002.

\noindent
18. F. Tricomi, Differential Equations, Izdatelstvo inostrannoj literatury (russian translation\linebreak \phantom{aaa} of the book F. G. Tricomi, Differential Equations, Blackie $\&$ Son Limited).

\end{document}